\documentclass[10pt]{article}
\usepackage{amssymb,amsmath,amsthm,amsfonts,mathrsfs}
\usepackage[pdftex,colorlinks,linktocpage,citecolor=black,linkcolor=black]{hyperref}
\usepackage{verbatim}
\usepackage{mathtools}
\usepackage[usenames]{color}
\usepackage{cite}
\usepackage{graphicx,graphics}
\usepackage{epstopdf,subfigure,tikz,float,wrapfig}
\usepackage{bm}
\usepackage{multirow,multicol}
\usepackage{comment}
\usepackage{booktabs}
\usepackage{threeparttable}
\usepackage[linesnumbered,ruled,algosection]{algorithm2e}
\graphicspath{{./pics/}}
\usepackage[active]{srcltx}
\allowdisplaybreaks

\usepackage{enumitem}

\numberwithin{equation}{section}
\newtheorem{theorem}{Theorem}[section]
\newtheorem{lemma}{Lemma}[section]
\theoremstyle{remark}
\newtheorem{remark}{Remark}[section]

\def\cT {\mathcal T}
\def\bold {\boldsymbol}
\def\Om {\Omega}
\def\J {\mathcal{J}}
\def\F {\mathcal{F}}
\def\U {\mathcal{U}}
\def\to{\rightarrow}
\def\eps {\varepsilon}

\renewcommand{\S}{\bold{\mathcal{S}}}
\newcommand{\E}{\bold{\mathcal{E}}}

\newcommand{\dx}{\,{\rm d}\bold{x}}
\newcommand{\ds}{\,{\rm d}s}
\newcommand{\dd}{\,{\rm d}}

\begin{document}

\title{\Large Convergence Analysis of an Adaptive Nonconforming FEM for Phase-Field Dependent Topology Optimization in Stokes Flow}
\author{Bangti Jin\thanks{Department of Mathematics, The Chinese University of Hong Kong, Shatin, New Territories, Hong Kong, China (b.jin@cuhk.edu.hk, bangti.jin@gmail.com)}
\and Jing Li\thanks{School of Mathematical Sciences, East China Normal University, Shanghai 200241, China. (betterljing@163.com)} \and Yifeng Xu\footnote{corresponding author. Department of Mathematics \& Scientific Computing Key Laboratory of Shanghai Universities, Shanghai Normal University, Shanghai 200234, China. (yfxu@shnu.edu.cn, yfxuma@aliyun.com)} \and Shengfeng Zhu\footnote{Department of Data Mathematics \& Shanghai Key Laboratory of Pure Mathematics and Mathematical Practice, School of Mathematical Sciences, East China Normal University, Shanghai 200241, China. (sfzhu@math.ecnu.edu.cn)}}

\date{}
\maketitle

\begin{abstract}
In this work, we develop an adaptive nonconforming finite element algorithm for the numerical approximation of phase-field parameterized topology optimization governed by the Stokes system. We employ the conforming linear finite element space to approximate the phase field, and the nonconforming linear finite elements (Crouzeix-Raviart elements) and piecewise constants to approximate the velocity field and the pressure field, respectively. We establish the convergence of the adaptive method, i.e., the sequence of minimizers contains a subsequence that converges to a solution of the first-order optimality system, and the associated subsequence of discrete pressure fields also converges. The analysis relies crucially on a new discrete compactness result of nonconforming linear finite elements over a sequence of adaptively generated meshes. We present numerical results for several examples to illustrate the performance of the algorithm, including a comparison with the uniform refinement strategy.



\end{abstract}

\section{Introduction}\label{sect:intro}

In this work, we investigate a numerical approach to topology optimization for fluid mechanics. Let $\Omega\subset\mathbb{R}^d$ ($d=2,3$) be an open bounded domain with a polygonal/polyhedral boundary $\partial\Omega$ with the unit outward normal $\bold{n}$. The Stokes-Darcy/Brinkman flow problem is governed by \cite{NieldBejan:2017}
\begin{align*}
 \left\{\begin{aligned}
   -\mu \bold{\Delta} \bold{u} + \alpha(\rho) \bold{u} + \nabla p&= \bold{f},\quad \mbox{in }\Omega, \\
       \mathrm{div} \bold{u}  &= 0 ,\quad \mbox{in }\Omega,\\
       \bold{u}  &= \bold{g}, \quad \mbox{on }\partial\Omega,
 \end{aligned}\right.
\end{align*}
where $\bold{u}$ is the flow velocity, $p$ is pressure, and $\mu>0$ is fluid viscosity. The body force $\bold{f}$ belongs to $\bold{L}^2(\Omega)$ and the velocity $\bold{g}\in \bold{H}^{1/2}(\partial\Omega)$ satisfies the compatibility condition $\int_{\partial\Omega} \bold{g}\cdot\bold{n}\ds = 0$.  The inverse permeability $\alpha:[0,1]\to[0,\infty]$ is decreasing, surjective and continuously differentiable. The function $\rho:\Omega\to \{0,1\}$ describes the distribution of the fluid in the design domain $\Omega$: the presence and absence of the fluid is indicated by a value of one and zero, respectively.

Topology optimization for fluid mechanics aims at finding an optimal distribution of a fluid within the design domain $\Omega$ that minimizes a problem-specific cost functional and has received a lot attention in both academia and industry \cite{Aleandersen:2020}. It was initiated by the seminal work of Borrvall and Petersson \cite{BorrvallPetersson:2003}, who proposed to minimize the dissipated power plus the potential power of the applied body force $\bold{f}$ in the Stokes-Darcy flow characterized by a regularized inverse variable permeability $\alpha$ and  the artificial porous medium replacing the non-fluid region. Their work has spurred many subsequent developments; see the references  \cite{Gersborg-Hansen:2005,Olesen:2006,GuestPrevost:2006,WikerKlarbringBorrvall:2007,Age:2008,Kreissl:2011,Alonso:2018,Sa:2018, Deng:2018,Alonso:2020,Thore:2021} for an incomplete list. To resolve the issue of the existence of an optimal design \cite{Evgrafov:2005}, Garcke et al \cite{GarckeHecht:2016a,GarckeHecht:2016b,GarckeHechtHinzeKahle:2015} incorporated a total variation (TV) term into the objective, and proposed a phase-field approach using the Ginzburg-Landau energy, replacing the material distribution $\rho$ by a phase field function $\phi:\Omega\to [0,1]$. The domain $\Omega$ can then be decomposed into a fluid region, a non-fluid region occupied by the porous medium and a thin interfacial layer, which are represented by $\phi$ as $\{x\in\Omega: \phi(x) = 1\}$, $\{x\in\Omega:\phi(x) = 0\}$ and $\{x\in\Omega: 0<\phi(x)<1\}$, respectively.

The topology optimization problem in the porous medium-phase field formulation then reads \cite{GarckeHecht:2016a,GarckeHecht:2016b,GarckeHechtHinzeKahle:2015}
\begin{align}
     \qquad \qquad \inf_{\phi \in \U}\mathcal{J}^\eps(\phi) : =  \frac{1}{2}\int_\Omega \alpha_{\eps}(\phi) |\bold{u}|^2 \dx + \int_\Om G(\bold{x},\bold{u},\bold{\nabla}\bold{u}) \dx + \gamma \mathcal{P}_\eps(\phi), \label{min_phase-field}
    \end{align}
subject to the constraint
\begin{equation}
       \bold{u} \in \bold{U}:~	
             \mu \int_\Omega \bold{\nabla} \bold{u} : \bold{\nabla} \bold{v} \dx + \int_\Omega \alpha_\eps(\phi) \bold{u} \cdot \bold{v} \dx  = \int_{\Omega}  \bold{f} \cdot \bold{v} \dx, \quad \forall \bold{v} \in \bold{V},\label{vp_stokes_div-free}
\end{equation}
where $\gamma>0$ is the regularization parameter. The function set $\bold{U}$ and the function space $\bold{V}$ are given by $\bold{U}:=\{\bold{v}\in \bold{H}^1(\Omega)~|~\mathrm{div}\bold{v} = 0,\bold{v}|_{\partial\Omega} = \bold{g} \}$ and $\bold{V}:=\{ \bold{v} \in \bold{H}_0^1(\Omega)~|~\mathrm{div}\bold{v} = 0\}$. The Ginzburg-Landau energy $\mathcal{P}_\epsilon$ and the function $G$ are given, respectively, by
\begin{equation}
 \mathcal{P}_\eps(\phi) = \dfrac{\eps}{2}\int_\Om |\bold{\nabla}\phi|^2 \dx + \dfrac{1}{\eps} \int_\Om f(\phi)  \dx \quad\mbox{and}\quad \label{def:functional}
    G(\bold{x},\bold{u},\bold{\nabla}\bold{u}): = \frac{\mu}{2} |\bold{\nabla}\bold{u}|^2  -  \bold{f} \cdot \bold{u},
 \end{equation}
where $\eps>0$ is a relaxation constant, controlling the width of the interfacial layer. Physically, $G$ represents the dissipative energy and the potential energy of the applied body force $\bold{f}$. $f(\phi)$ is either a double obstacle potential $\frac{1}{2}\phi(1-\phi)$ or a double well potential $\frac{1}{4}\phi^2(1-\phi)^2$. The admissible set $\U$ for the phase-field function $\phi$ is given by
\begin{equation*}
    \U : = \left\{\phi \in H^1(\Omega)~|~\phi \in [0,1]~\text{a.e. in}~\Omega, \int_{\Om}\phi\dx \leq \beta |\Om| \right\} ,
\end{equation*}
where $\beta\in (0,1)$ is the maximal volume fraction of the fluid in $\Omega$. Note that the regularized inverse permeability $\alpha_\eps:[0,1]\to[0,\overline{\alpha}_\eps]$ is
coupled with the phase field parameter $\eps$ so that it stays uniformly bounded by $\overline{\alpha}_\eps>0$.

The numerical treatment of problem \eqref{min_phase-field}--\eqref{vp_stokes_div-free} often employs a gradient flow approach (see, e.g., \cite{GarckeHechtHinzeKahle:2015}). Then the resulting problem is approximated by conforming linear elements for the phase field $\phi$ and lowest-order conforming Taylor-Hood ($P_2$-$P_1$) elements for the velocity $\bold{u}$ and pressure $p$, respectively. There is much interest in reducing the computational cost, especially in the 3d case \cite{Thore:2021}. The recent study  \cite{JinLiXu:2025}  indicates that using nonconforming linear (Crouzeix-Raviart or CR) elements \cite{CrouzeixRaviart:1973} and piecewise constants for the velocity $\bold{u}$ and the pressure $p$, respectively, can achieve this goal. Moreover, adaptive finite element methods (AFEM) using a residual-type \textit{a posteriori} error estimate for the phase field can be used to effectively resolve the interfacial layer between the fluid and porous medium. Thus, a natural idea is to develop an adaptive nonconforming FEM for problem \eqref{min_phase-field}-\eqref{vp_stokes_div-free}. In this work, we shall develop an adaptive algorithm, which consists of successive loops of the cycle
\begin{equation}\label{afem_loop}
      \text{OPTIMIZE} \rightarrow \text{ESTIMATE} \rightarrow \text{MARK} \rightarrow \text{REFINE}.
\end{equation}
Here OPTIMIZE refers to solving the discrete problem \eqref{dismin_phase-field}-\eqref{disvp_stokes_div-free} for the minimizing pair $(\phi_k^\ast, \bold{u}^\ast_k) \in \U_k \times \bold{U}_k$ over the mesh $\cT_k$; ESTIMATE are some computable quantities, termed as \textit{a posteriori} error estimator, to measure the discretization error; MARK determines some elements in $\cT_k$ by a criterion associated with the error estimator; REFINE produces a finer mesh $\cT_{k+1}$ by local refinements of all marked elements and their neighbors for conformity if necessary. The \textit{a posteriori} error estimator has witnessed significant advances \cite{AinsworthOden:2000,Verfurth:2013} for FEM approximations of forward problems, especially elliptic equations; see the {survey papers} \cite{BonitoCanutoNochettoVeeser:2024,CarstensenFeischlPagePraetorius:2014,NochettoSiebertVeeser:2009} for the theory of AFEM. For the analysis of adaptive CR FEMs for the Stokes system, we refer readers to the important works \cite{BeckerMao:2011,CarstensenPeterseimRabus:2013,DariDuranPadra:1995,HuXu:2013}. 

In Section \ref{sect:adaptive}, we derive \textit{a posteriori} error estimators to drive the adaptive algorithm in Algorithm \ref{alg_anfem_topopt-Stokes_phase-field_G1}. The derivation is based on the first-order necessary optimality system \eqref{opt-sys_G1} of problem \eqref{min_phase-field}--\eqref{vp_stokes_div-free}. The experimental study in Section \ref{sect:numerics} indicates that the adaptive  algorithm outperforms the uniform refinement strategy on several benchmark tests. In Section \ref{sect:conv_adaptive}, we establish the convergence of the adaptive algorithm in the sense that a subsequence of adaptively generated discrete minimizers and associated discrete velocity fields converges to a pair satisfying the optimality system \eqref{opt-sys_G1}.  The derivation of an adaptive nonconforming FEM for problem \eqref{min_phase-field}-\eqref{vp_stokes_div-free}, and the convergence analysis represent the main contributions of the work.

Next we outline the analysis strategy. The dependence on a mesh $\cT_k$ is indicated by the refinement level $k\geq 0$ in the subscript.
Let $\{(\phi^\ast_k,\bold{u}^\ast_k)\}_{k\geq 0}$ be an adaptively generated sequence of minimizing pairs. First we prove a null sequence of estimators $\{\eta_{k,2}(\phi_k^\ast,\bold{u}_k^\ast)\}_{k\geq0}$ in Theorem \ref{thm:conv_est_state&costate}, by which we establish the $L^2$ weak convergence of a subsequence $\{\bold{u}^\ast_{k_j}\}_{j\geq 0}$ in Lemma \ref{lem:vel&costate_weak_adaptive}. Then using a variant of the connection operator \cite{Brenner:2003} and the vanishing limit of a computable quantity, measuring the discontinuity of $\bold{u}_k^\ast$ in $\eta_{k,2}$, we prove the $L^2$ strong convergence of $\{\bold{u}^\ast_{k_j}\}_{j\geq 0}$ in Lemma \ref{lem:vel&costate_L2-strong-conv_adaptive}. This discrete compactness property of adaptively generated solutions by the lowest-order CR finite elements generalizes the classic result in \cite{Stummel:1980} (see also \cite[Chapter 10.3]{ShiWang:2013}) over uniformly refined meshes to locally refined meshes and to the best of our knowledge, is a new contribution to the analysis of AFEM. Since the sequence $\{\U_k\}_{k\geq 0}$ is nested, we may define an auxiliary variational problem $\phi_\infty \in \U_\infty:= \overline{\bigcup_{k\geq 0}\U_k}$ in the $H^1(\Omega)$-norm (see \eqref{vp_stokes_div-free_limit}), and study the continuity of $\S_k$ along a convergent subsequence of $\{ \phi_{k}^\ast\}_{k\geq 0}$ with a limit $\overline{\phi}_\infty \in \U_\infty$ in Lemma \ref{lem:sol-map_cont_limit}. Moreover, we prove the convergence of the subsequence $\{\phi_{k_j}^\ast\}_{j\geq0}$ in the $H^1(\Omega)$ norm. These tools and the arguments of \cite{GantnerPraetorius:2022,MorinSiebert:2008,Siebert:2011} allow one to extract a vanishing limit of $\{\eta_{k_j,1}(\phi_{k_j}^\ast,\bold{u}_{k_j}^\ast)\}_{j\geq0}$ in Theorem \ref{thm:conv_est_control}. By the two null subsequences of estimators, we prove that the limit pair of the convergent subsequence $\{(\phi_{k_j}^\ast,\bold{u}_{k_j}^\ast)\}_{j\geq0}$ solves the optimality system \eqref{opt-sys_state_G1} in Theorem \ref{thm:conv_adaptive}. Finally, in Theorem \ref{thm:conv_pre_adaptive}, we prove the convergence of discrete pressure fields $\{p_{k_j}^\ast\}_{j\geq 0}$ in the $L^2(\Omega)$ norm.

Now we situate this work in the context of existing studies. First, the numerical analysis of topology optimization using the Borrvall-Petersson model was investigated with the uniform refinement strategy in \cite{BorrvallPetersson:2003,PapadopoulosSuli:2022,Papadopoulos:2022}. For a piecewise constant FE approximation of the material distribution $\rho$ and an inf-sup stable quadrilateral FE approximation of velocity $\bold{u}$ and pressure $p$, Borrvall and Petersson \cite{BorrvallPetersson:2003} showed the convergence of the discrete velocity $\bold{u}^*_h$ and material distribution $\rho_h^*$ to a minimizer of the model weakly in $H^1(\Omega)^d$ and weakly $\ast$ in $L^\infty(\Omega)$, respectively. For any inf-sup stable conforming pair of FE spaces for the velocity and the pressure, Papadopoulos and Suli \cite{PapadopoulosSuli:2022}  proved that for every isolated minimizer, there exist sequences of FE solutions to the discretized first-order optimality system that strongly converge to the minimizer in $H^1(\Omega)^d$ and $L^s(\Omega)$, $s\in [1,\infty)$, respectively. The analysis was extended in \cite{Papadopoulos:2022} to discontinuous Galerkin methods with an interior penalty. Recently the authors \cite{JinLiXu:2025} presented a convergence analysis of the CR approximation for the phase-field dependent model in \cite{GarckeHecht:2016b}. Second, there has been recent progress in the application and analysis of AFEM for topology optimization, including minimum compliance problem in structural optimization \cite{JinLiXuZhu:2024} and elliptic eigenvalue optimization \cite{LiXuZhu:2023}. For the optimal design of fluid, Garcke et al \cite{GarckeHechtHinzeKahle:2015} proposed an adaptive algorithm for a generalized Cahn-Hilliard equation arising from \eqref{min_phase-field}--\eqref{vp_stokes_div-free}, and showed that the errors of the phase-field and the auxiliary variable in the $H^1(\Omega)$-seminorm are bounded from above by an \textit{a posteriori} error estimate up to a high-order term. However, all these existing works investigate conforming finite elements. In contrast, this work focuses on the convergence analysis based on the optimality system of problem \eqref{min_phase-field}-\eqref{vp_stokes_div-free}, including an \textit{a posteriori} error estimator for the nonconforming  approximation of the velocity by CR elements. The analysis requires new technical tools to overcome the challenge due to the nonconformity of the CR FE space.

The rest of the paper is organized as follows. We present an adaptive algorithm in Section \ref{sect:adaptive}, and establish its convergence in Section \ref{sect:conv_adaptive}. We present numerical results in Section \ref{sect:numerics}.
Throughout, we use standard notation for Sobolev spaces and related (semi-)norms \cite{EvansGariepy:2015}. The notation $c$, with or without subscript, denotes
a generic constant independent of the mesh size and may differ at each occurrence, while $(\cdot,\cdot)_{L^2(G)}$ denotes the $L^2(G)$ scalar product on an
open bounded set $G\subset\mathbb{R}^d$ and the subscript is omitted when $G=\Omega$. For a matrix $\bold{\tau}\in\mathbb{R}^{d\times d}$ and a
 vector $\bold{v}\in \mathbb{R}^d$, $\bold{\tau}\cdot \bold{v}$ and $\bold{\tau}\times\bold{v}$ are defined row-wisely.
The divergence operator $\mathbf{div}$ is applied to $\bold{\tau}$ row-wisely:
$\mathbf{div}\bold{\tau} = (\mathrm{div}\bold{\tau}_{i})_{1\leq i\leq d} =(\sum_{j=1}^d\partial_j\tau_{ij})_{1\leq i\leq d}$.

\section{Adaptive algorithm}\label{sect:adaptive}
In this section we develop an adaptive algorithm for problem \eqref{min_phase-field}-\eqref{vp_stokes_div-free}.

\subsection{Nonconforming FEM approximation}

Let $\cT_0$ be a shape regular triangulation of $\overline{\Omega}$ into closed triangles/tetrahedra
\cite{Ciarlet:2002} and $\mathbb{T}$ be the set of all possible conforming triangulations generated  from $\cT_0$ by
successive bisections \cite{Kossaczky:1995,NochettoSiebertVeeser:2009}. Note that the shape regularity of any
$\mathcal{T}\in\mathbb{T}$ is bounded uniformly by a constant depending only on $\cT_0$ \cite{NochettoSiebertVeeser:2009,Traxler:1997}.
The collection of all edges/faces (respectively interior edges/faces) in each $\cT\in \mathbb{T}$ is denoted by $\F_\cT$ (respectively
$\F_{\cT}(\Omega)$). For each $F\in\mathcal{F}_\cT$, we fix a unit normal vector $\bold{n}_F$, pointing from
$T_+$ to $T_-$ if $T_+, T_-\in \cT$ share a common edge/face $F \in\mathcal{F}_\cT(\Omega)$, and coinciding with the
unit outward normal $\bold{n}$ to $\partial\Omega$ if $F\subset\partial\Omega$, and $\bold{n}_{\partial T}$ denotes the unit outward normal
to the boundary $\partial T$ of each $T\in\cT$. We use $\mathcal{N}_{\cT}(F)$ for the set of elements sharing an $F\in \mathcal{F}_\cT$
and $\mathcal{N}_{\cT}(T)$ for the set of any $T\in \cT$ and its neighboring elements in $\cT$, the union of which is denoted
by $\omega_\cT(T)$.    The local mesh-size $h_T$ ($h_F$) of each $T\in \cT$ ($F\in \F_\cT$) is defined
as $|T|^{1/d}$ ($|F|^{1/(d-1)}$), which yields a piecewise constant mesh-size function $h_\cT:\overline{\Omega}
\to \mathbb{R}^+$ with $h_\cT|_{T} : = |T|^{1/d}$ for each $T\in \cT$. $P_m(T)$ denotes the set of all polynomials of total degree no more than $m$ on $T \in \cT$ for $m\in\mathbb{N}\cup\{0\}$. Now for any  $\cT\in \mathbb{T}$, we define a discrete admissible set $\U_\cT:=S_\cT\cap\U$ (with $S_\cT$ being the $H^1(\Omega)$-conforming linear FE space), the nonconforming linear FE space \cite{CrouzeixRaviart:1973} and the piecewise constant space
\begin{align}\label{C-R_space}
    \bold{X}_\cT &: = \left\{ \bold{v}_\cT \in \bold{L}^2(\Omega)~|~\bold{v}_\cT|_T \in P_1(T)^d~\forall T\in\cT,~([\bold{v}_\cT],1)_{L^2(F)}=\bold{0}~\forall F\in\mathcal{F}_{\cT}(\Omega)\right\},\\
\label{discretespace_pressure}
    Q_{\cT}&:=\{v_\cT \in L_{0}^{2}(\Omega)~|~v_{\cT}|_T \in P_{0}(T)~\forall T\in\mathcal{T}\},
\end{align}
where $[\boldsymbol{v}_\cT]:= \bold{v}_\cT|_{F\subset T_+} - \bold{v}_\cT|_{F\subset T_-}$ denotes the jump across an interior face/edge $F = \partial T_+ \cap \partial T_-$ and takes the value of
$\boldsymbol{v}_\cT$ on $F\subset\partial\Omega$, and $L_0^2(\Omega)$ is a subspace of $L^2(\Omega)$ with zero mean over $\Omega$. A function in $\bold{X}_\cT$ is determined by its values at the centers of faces/edges in $\F_\cT$ and the members in $\bold{X}_\cT$ vanishing at the boundary centers comprise the subspace $\bold{Z}_\cT$. $\bold{\nabla}_\cT$ (respectively $\mathrm{div}_\cT$) denotes the piecewise gradient (respectively divergence) operator over $\cT$. Then $\mathrm{div}_\cT \bold{v}_\cT \in Q_\cT$ for any $\bold{v}_\cT \in \bold{Z}_\cT$. One can define a piecewise divergence free FE subspace
\begin{equation}\label{C-R_space_div-free}
    \bold{V}_\cT : = \left\{\bold{v}_\cT \in \bold{Z}_\cT ~|~( (\mathrm{div}_{\cT}  \bold{v}_\cT ), q_\cT) = 0~\forall q_\cT \in Q_\cT\} = \{\bold{v}_\cT \in \bold{Z}_\cT ~|~\mathrm{div}_{\cT} \bold{v}_\cT = 0\right\},
\end{equation}
and equip $\bold{Z}_{\cT}$ with a mesh-dependent energy norm
$\|\boldsymbol{v}\|^{2}_{1,\cT}:=\sum_{T\in\cT}\int_{T}|\bold{\nabla}\bold{v}|^2\dx=\int_\Omega |\boldsymbol{\nabla}_{\cT}\boldsymbol{v}|^2 \dx,$
which coincides with the usual $\bold{H}^1_0(\Omega)$ norm due to Poincar\'{e} inequality. Over $\bold{Z}_\cT$, there holds a discrete Poincar\'{e} inequality \cite{Brenner:2003, BuffaOrtner:2009}
\begin{equation}\label{disc_Poin_ineq}
    \|\bold{v}_\cT\|_{\bold{L}^2(\Omega)} \leq c_{\text{dp}} \|\bold{v}_\cT\|_{1,\cT},\quad\forall \bold{v}_\cT \in \bold{Z}_\cT,
\end{equation}
with $c_{\text{dp}}$ depending only on the shape regularity of $\cT$.
Also on the Crouzeix-Raviart (CR) FE space $\bold{X}_\cT$, we define an interpolation operator $\bold{\Pi}_\cT:\bold{H}^1(\Omega) \to \bold{X}_\cT$ by \cite{CrouzeixRaviart:1973, OrtnerPraetorius:2011}
\begin{equation}\label{int_operator-cr}
  ( \bold{\Pi}_\cT \bold{v},1)_{L^2(F)}  = ( \bold{v},1)_{L^2(F)}, \quad \forall F \in \mathcal{F}_\cT.
\end{equation}
It satisfies the following approximation property
\begin{equation}\label{int_operator-cr_approx}
  \|\bold{v}-\bold{\Pi}_\cT\bold{v}\|_{\bold{L}^2(T)}\leq c h_T \|\bold{\nabla}\bold{v}\|_{\bold{L}^2(T)},\quad \forall T \in \cT,
\end{equation}
with the positive constant $c$ depending only on the shape regularity of $\cT$ and the stability estimate
\begin{equation}\label{int_operator-cr_stab}
    \|\bold{\nabla}(\bold{\Pi}_\cT \bold{v})\|_{\bold{L}^2(T)} \leq \|\bold{\nabla}\bold{v}\|_{\bold{L}^2(T)}.
\end{equation}
By the assumption on $\bold{g}$ and \cite[Lemma 2.2, Chapter I]{GiraultRaviart:1986}, the boundary condition $\bold{g}$ admits a lift $\bold{w}\in \bold{U}$.
By \eqref{int_operator-cr} and $\mathrm{div}\bold{w}=0$, we have $\mathrm{div}_\cT \bold{\Pi}_\cT\bold{w} = 0$.
The FEM approximation of problem \eqref{vp_stokes_div-free}-\eqref{min_phase-field} reads
\begin{subequations}\label{dismin}
    \begin{align}
         &\qquad\qquad\inf_{\phi_\cT \in \U_\cT}\mathcal{J}^\eps_\cT(\phi_\cT) : =  \tfrac{1}{2}(\alpha_{\eps}(\phi_\cT) ,|\bold{u}_\cT|^2) + ( G(\bold{x},\bold{u}_\cT,\bold{\nabla}_\cT\bold{u}_\cT),1) + \gamma \mathcal{P}_\eps(\phi_\cT), \label{dismin_phase-field}\\
        & \text{s.t.}~\bold{u}_\cT \in \bold{\Pi}_\cT \bold{w} + \bold{V}_\cT:~	
             \mu ( \bold{\nabla}_\cT \bold{u}_\cT, \bold{\nabla}_\cT \bold{v}_\cT) + ( \alpha_\eps(\phi_\cT) \bold{u}_\cT, \bold{v}_\cT)  = ( \bold{f}, \bold{v}_\cT), \quad \forall  \bold{v}_\cT \in \bold{V}_\cT. \label{disvp_stokes_div-free}
    \end{align}
\end{subequations}
The existence of a minimizer $\phi_\cT^*\in \mathcal{U}_\cT$ is ensured \cite[Theorem 3.1]{JinLiXu:2025}.

\subsection{Adaptive phase-field algorithm}
Now we present an adaptive algorithm of form \eqref{afem_loop} for problem \eqref{min_phase-field}-\eqref{vp_stokes_div-free}. We indicate by the subscript $k$ the quantities on the mesh $\mathcal{T}_k$, e.g., $\mathcal{N}_{k}(F)$ for $\mathcal{N}_{\cT_k}(F)$. The basis of the adaptive algorithm is some computable quantities related to the associated optimality system. Let $(\phi^\ast, \bold{u}^\ast)\in \U \times \bold{U}$ be a minimizer and the associated velocity field of problem \eqref{min_phase-field}-\eqref{vp_stokes_div-free}. It may be formulated as \cite{GarckeHecht:2015}
\begin{equation}\label{min_G1_phase-field}
    \inf_{(\phi,\bold{u}) \in \U \times \bold{U}}\mathcal{J}^\eps(\phi,\bold{u}) : =  \tfrac{1}{2}( \alpha_{\eps}(\phi), |\bold{u}|^2) + \tfrac{\mu}{2}  \|\bold{\nabla}\bold{u}\|_{L^2(\Omega)}^2 - (\bold{f},\bold{u})  + \gamma \mathcal{P}_\eps(\phi).
\end{equation}
The optimality system of problem \eqref{min_G1_phase-field} reads
\begin{subequations}\label{opt-sys_G1}
\begin{equation}\label{opt-sys_control_G1}
    (\gamma \eps \bold{\nabla} \phi^\ast , \bold{\nabla} (\phi - \phi^\ast )) + ( \left( \tfrac{\gamma}{\eps} f'(\phi^\ast) + \tfrac{1}{2}\alpha'_\eps(\phi^\ast) |\bold{u}^\ast|^2 \right),  \phi - \phi^\ast) \geq 0, \quad \forall \phi \in \U,
\end{equation}
\begin{equation}\label{opt-sys_state_G1}
    \mu(\bold{\nabla} \bold{u}^\ast, \bold{\nabla} \bold{v}) + ( \alpha_\eps(\phi^\ast) \bold{u}^\ast, \bold{v}) = (\bold{f} ,\bold{v}), \quad \forall \bold{v} \in \bold{V}.
\end{equation}
\end{subequations}
Let $\{\cT_k\}_{k\geq 0}$ be a sequence of meshes over $\overline{\Omega}$ generated by an adaptive algorithm. Then a minimizing pair $(\phi^\ast_k,\bold{u}^\ast_k)\in \U_k \times (\bold{\Pi}_k \bold{w} + \bold{V}_k)$ of the discrete problem \eqref{dismin_phase-field}-\eqref{disvp_stokes_div-free} over the mesh $\cT_k$ satisfies
\begin{subequations}\label{opt-sys_disc_G1}
\begin{equation}\label{opt-sys_disc_control_G1}
    (\gamma \eps \bold{\nabla} \phi^\ast_k,\bold{\nabla} (\phi_k - \phi^\ast_k ))+  ( \left( \tfrac{\gamma}{\eps} f'(\phi^\ast_k) + \tfrac{1}{2}\alpha'_\eps(\phi^\ast_k) |\bold{u}_k^\ast|^2 \right),  \phi_k - \phi^\ast_k) \geq 0, \quad \forall \phi_k \in \U_k,
\end{equation}
\begin{equation}\label{opt-sys_disc_state_G1}
    \mu(\bold{\nabla}_k \bold{u}_k^\ast, \bold{\nabla}_k \bold{v}_k) + ( \alpha_\eps(\phi_k^\ast) \bold{u}^\ast_k, \bold{v}_k) = ( \bold{f},\bold{v}_k), \quad \forall \bold{v}_k \in \bold{V}_k.
\end{equation}
\end{subequations}

Next we derive the a posteriori error estimators. The derivation  uses the following quasi-interpolation operator $\Pi_k: L^1(\Omega) \to S_k$ \cite{JinLiXuZhu:2024,JinXu:2019}
\[
    \Pi_k v : = \sum_{\bold{x}\in\mathcal{V}_k} \frac{1}{|\omega_{\bold{x}}|} \int_{\omega_{\bold{x}}} v \dx \varphi_{\bold{x}},
\]
where $\mathcal{V}_k$ is the set of all vertices associated with $\cT_k$, $\{\varphi_{\bold{x}}\}_{\bold{x}\in \mathcal{V}_k}$ the set of all nodal basis functions in $S_k$, $\omega_{\bold{x}}$ the support of $\varphi_{\bold{x}}$ with respect to $\bold{x}\in \mathcal{V}_k$ and $|\omega_{\bold{x}}|$ denotes the Lebesgue measure of $\omega_{\bold{x}}$. $\Pi_k$ preserves the box constraint and the volume constraint in $\U$, i.e.,  $\Pi_k \phi \in \U_k$ if $\phi\in \U$ \cite{JinLiXuZhu:2024}.
For any $\phi \in \U$, by taking $\phi_k = \Pi_k \phi$ in \eqref{opt-sys_disc_control_G1}, we get
\begin{align}
    &\quad (\gamma \eps \bold{\nabla} \phi^\ast_k, \bold{\nabla} (\phi - \phi^\ast_k )) + ( \tfrac{\gamma}{\eps} f'(\phi^\ast_k) + \tfrac{1}{2}\alpha'_\eps(\phi^\ast_k) |\bold{u}_k^\ast|^2, \phi - \phi^\ast_k) \nonumber \\
    & = (\gamma \eps \bold{\nabla} \phi^\ast_k , \bold{\nabla} (\phi - \Pi_k\phi )) +  ( \tfrac{\gamma}{\eps} f'(\phi^\ast_k) + \tfrac{1}{2}\alpha'_\eps(\phi^\ast_k) |\bold{u}^\ast_k|^2, \phi - \Pi_k\phi) \nonumber \\
    & \quad +  (\gamma \eps \bold{\nabla} \phi^\ast_k , \bold{\nabla} ( \Pi_k\phi - \phi^\ast_k )) + ( \tfrac{\gamma}{\eps} f'(\phi^\ast_k) + \tfrac{1}{2}\alpha'_\eps(\phi^\ast_k) |\bold{u}^\ast_k|^2, \Pi_k\phi - \phi^\ast_k) \nonumber \\
    &\ge (\gamma \eps \bold{\nabla} \phi^\ast_k, \bold{\nabla} (\phi - \Pi_k\phi )) + ( \tfrac{\gamma}{\eps} f'(\phi^\ast_k) + \tfrac{1}{2}\alpha'_\eps(\phi^\ast_k) |\bold{u}^\ast_k|^2,\phi - \Pi_k\phi). \label{estimate-pf01_control_G1}
\end{align}
Elementwise integration by parts and the Cauchy-Schwarz inequality lead to
\begin{align}
    &\quad \bigg|(\gamma \eps \bold{\nabla} \phi^\ast_k, \bold{\nabla} (\phi - \Pi_k\phi )) + (\left( \tfrac{\gamma}{\eps} f'(\phi^\ast_k) + \tfrac{1}{2}\alpha'_\eps(\phi^\ast_k) |\bold{u}^\ast_k|^2 \right) ,\phi - \Pi_k\phi) \bigg| \nonumber \\
    & = \bigg| \sum_{T \in \cT_k}\left( ( \tfrac{\gamma}{\eps} f'(\phi^\ast_k) + \tfrac{1}{2}\alpha'_\eps(\phi^\ast_k) |\bold{u}^\ast_k|^2 ,\phi - \Pi_k\phi)_{L^2(T)} + (\gamma \eps \bold{\nabla} \phi^\ast_k \cdot \bold{n}_{\partial T},\phi - \Pi_k\phi)_{L^2(\partial T)}  \right)\bigg| \nonumber \\
    & = \bigg| \sum_{T\in \cT_k} ( \tfrac{\gamma}{\eps} f'(\phi^\ast_k) + \tfrac{1}{2}\alpha'_\eps(\phi^\ast_k) |\bold{u}^\ast_k|^2 ,\phi - \Pi_k\phi)_{L^2(T)} + \sum_{F \in  \mathcal{F}_k}  (\gamma \eps [\bold{\nabla}\phi_k^\ast]\cdot\bold{n}_F, \phi - \Pi_k\phi)_{L^2(F)}  \bigg | \nonumber \\
    & \leq \sum_{T\in \cT_k}\bigg( \left\| \tfrac{\gamma}{\eps} f'(\phi^\ast_k) + \tfrac{1}{2}\alpha'_\eps(\phi^\ast_k) |\bold{u}^\ast_k|^2 \right\|_{L^2(T)}\|\phi - \Pi_k\phi\|_{L^2(T)} + \sum_{F\subset \partial T} \gamma\eps\|[\bold{\nabla}\phi_k^\ast]\cdot\bold{n}_F\|_{L^2(F)} \|\phi - \Pi_k\phi\|_{L^2(F)} \bigg). \label{estimate-pf02_control_G1}
\end{align}
By defining a residual functional $\mathcal{R}_{1}(\phi_k^\ast,\bold{u}_k^\ast)$ associated with $(\phi_k^\ast,\bold{u}_k^\ast)$ over $H^1(\Omega)$ by
\begin{equation}\label{def:res1_G1}
    \langle \mathcal{R}_{1}(\phi_k^\ast,\bold{u}_k^\ast),\phi \rangle : = (\gamma \eps \bold{\nabla} \phi^\ast_k, \bold{\nabla} \phi) +  (\left( \tfrac{\gamma}{\eps} f'(\phi^\ast_k) + \tfrac{1}{2}\alpha'_\eps(\phi^\ast_k) |\bold{u}_k^\ast|^2 \right),  \phi), \quad \forall \phi \in H^1(\Omega),
\end{equation}
and using error estimates for $\Pi_k$ \cite[Lemma 5.3]{JinXu:2019}, we finally arrive at
\begin{align}\label{estimate-pf03_control_G1}
    \left| \langle \mathcal{R}_{1}(\phi_k^\ast,\bold{u}_k^\ast),\phi - \Pi_k\phi\rangle \right| & \leq c \sum_{T\in \cT_k} \left(h_T \|R_{T,1}(\phi_k^\ast,\bold{u}^\ast_k)\|_{L^2(T)} + \sum_{F\subset\partial T}h_F^{1/2} \|J_{F,1}(\phi_k^\ast)\|_{L^2(F)}\right)   \|\bold{\nabla}\phi\|_{L^2(\omega_k(T))} \nonumber \\
    &  \leq c \sum_{T\in\cT_k} \eta_{k,1}(\phi_k^\ast,\bold{u}^\ast_k;T) \|\bold{\nabla}\phi\|_{L^2(\omega_k(T))} \leq c \eta_{k,1}(\phi_k^\ast,\bold{u}^\ast_k;\cT_k)\|\bold{\nabla}\phi\|_{L^2(\Omega)},\quad \forall \phi \in \U,
\end{align}
where $R_{T,1}(\phi_k^\ast,\bold{u}^\ast_k)$ and $J_{F,1}(\phi_k^\ast)$ denote the element residual on each element $T\in \cT_k$ and the jump residual on each face $F\in\mathcal{F}_{k}$, given respectively by
\[
    R_{T,1}(\phi_k^\ast,\bold{u}^\ast_k) : = \frac{\gamma}{\eps} f'(\phi^\ast_k) + \frac{1}{2}\alpha'_\eps(\phi^\ast_k) |\bold{u}_k^\ast|^2\quad\mbox{and} \quad J_{F,1}(\phi_k^\ast):=\left\{\begin{array}{l}
        \gamma\eps[\bold{\nabla}\phi_k^\ast]\cdot\bold{n}_F, \quad \forall F\in \mathcal{F}_k(\Omega), \\ [2ex]
        \gamma\eps\bold{\nabla}\phi_k^\ast\cdot\bold{n}, \quad \forall F\in \mathcal{F}_k\setminus\mathcal{F}_k(\Omega).
    \end{array} \right.
\]
Then one can define a local error indicator and a global error estimator respectively by
\begin{align*}
    \eta_{k,1}(\phi_k^\ast,\bold{u}_k^\ast;T)&:=\bigg(h_T^2\|R_{T,1}(\phi_k^\ast,\bold{u}^\ast_k)\|_{L^2(T)}^2 + \sum_{F\subset\partial T}h_F \|J_{F,1}(\phi_k^\ast)\|^2_{L^2(F)}\bigg)^{1/2},\\
    \eta_{k,1}(\phi_k^\ast,\bold{u}_k^\ast;\mathcal{M})&=\bigg(\sum_{T\in \mathcal{M}}\eta_{k,1}^2(\phi_k^\ast,\bold{u}_k^\ast;T)\bigg)^{1/2} ,\quad \forall \mathcal{M}\subseteq \cT_k.
\end{align*}
In view of \eqref{estimate-pf03_control_G1}, the computable quantity $\eta_{k,1}$ can bound the residual functional $\mathcal{R}_{1}$ in \eqref{opt-sys_control_G1}. To obtain an estimator for \eqref{opt-sys_state_G1}, we define a second residual functional $\mathcal{R}_{2}(\phi_k^\ast,\bold{u}_k^\ast)$ associated with $(\phi_k^\ast,\bold{u}_k^\ast)$ over $\bold{H}_0^1(\Omega)$ by
\begin{equation}\label{def:res2_G1}
    \langle \mathcal{R}_{2}(\phi_k^\ast,\bold{u}_k^\ast),\bold{v} \rangle : = \mu( \bold{\nabla}_k \bold{u}_k^\ast ,\bold{\nabla} \bold{v} ) + ( \alpha_\eps(\phi^\ast_k) \bold{u}_k^\ast, \bold{v} ) - ( \bold{f}, \bold{v}) ,\quad \forall \bold{v} \in \bold{H}_0^1(\Omega).
\end{equation}
In view of \eqref{int_operator-cr}, we have $\mathrm{div}_k\bold{\Pi}_k\bold{v} = 0$ for any $\bold{v} \in \bold{V}$. Inserting $\bold{v}_k = \bold{\Pi}_k \bold{v} $ in \eqref{opt-sys_disc_state_G1} and arguing similarly yield
\begin{align*}
     \langle \mathcal{R}_{2}(\phi_k^\ast,\bold{u}_k^\ast),\bold{v} \rangle &= \mu( \bold{\nabla}_k \bold{u}_k^\ast , \bold{\nabla}_k ( \bold{v} - \bold{\Pi}_k \bold{v} )) + ( \alpha_\eps(\phi^\ast_k) \bold{u}_k^\ast, \bold{v} - \bold{\Pi}_k \bold{v} ) - (\bold{f} , \bold{v} - \bold{\Pi}_k \bold{v} ) \\
     &= \sum_{T\in\cT_k}\left( (\alpha_\eps(\phi_k^\ast) \bold{u}_k^\ast - \bold{f},  \bold{v} - \bold{\Pi}_k \bold{v} )_{L^2(T)} + ( \bold{v} - \bold{\Pi}_k \bold{v} , \mu\bold{\nabla}\bold{u}^\ast_k \cdot \bold{n}_{\partial T})_{L^2(\partial T)} \right)\\
     &= \sum_{T\in\cT_k} (\alpha_\eps(\phi_k^\ast) \bold{u}_k^\ast - \bold{f},  \bold{v} - \bold{\Pi}_k \bold{v} )_{L^2(T)}, \quad \forall \bold{v} \in \bold{V},
\end{align*}
where the last line follows from the fact that $\mu\bold{\nabla}\bold{u}^\ast_k \cdot \bold{n}_{\partial T}$ is a constant vector on $F\subset \partial T$ and the identity  \eqref{int_operator-cr}. Then we deduce from \eqref{int_operator-cr_approx} that
\begin{align}
    \left| \langle \mathcal{R}_{2}(\phi_k^\ast,\bold{u}_k^\ast),\bold{v} \rangle \right|
    &\leq c \sum_{T\in \cT_k}h_T\|\alpha_\eps(\phi_k^\ast) \bold{u}_k^\ast - \bold{f}\|_{\bold{L}^2(T)}\|\bold{\nabla}\bold{v}\|_{\bold{L}^2(T)} \nonumber \\
    & \leq c \bigg(\sum_{T\in \cT_k}h_T^2\|\alpha_\eps(\phi_k^\ast) \bold{u}_k^\ast - \bold{f}\|_{\bold{L}^2(T)}^2\bigg)^{1/2} \|\bold{\nabla}\bold{v}\|_{\bold{L}^2(\Omega)},\quad \forall \bold{v}\in \bold{V}. \label{estimate-pf01_state_G1}
\end{align}
Moreover, \eqref{opt-sys_disc_state_G1} is  an exterior approximation of problem \eqref{opt-sys_state_G1} since $\bold{X}_k \nsubseteq \bold{H}^1(\Omega)$. Naturally, a computable quantity is expected to measure the discontinuity of $\bold{u}^\ast_k$ across an interelement $F\in \mathcal{F}_k(\Omega)$ and the mismatch with the Dirichlet boundary condition. For any $ \bold{\varphi} \in C^\infty(\overline{\Omega})^{d\times d}$, by elementwise integration by parts, we obtain
\begin{align}\label{estimate-pf02_state_G1}
      (\bold{u}_k^\ast, \mathbf{div} \bold{\varphi} )
      & =  \sum_{T\in\cT_{k}} \bigg( - (\bold{\nabla} \bold{u}_k^\ast , \bold{\varphi} )_{L^2(T)} + \sum_{F\subset \partial T}  ( \bold{u}_k^\ast , \bold{\varphi} \cdot \bold{n}_{\partial T} )_{L^2(F)}  \bigg)\nonumber\\
      & = - \sum_{T\in\cT_k}( \bold{\nabla} \bold{u}_k^\ast , \bold{\varphi})_{L^2(T)} + \sum_{F\in\mathcal{F}_k}( [\bold{u}_k^\ast], \bold{\varphi}\cdot\bold{n}_F )_{L^2(F)}.
\end{align}
The desired measure may be quantified by bounding each summand in the second term of \eqref{estimate-pf02_state_G1}. In view of the identities $$\int_F [\bold{u}_k^\ast] \ds = \bold{0},\quad \forall F\in \mathcal{F}_k(\Omega) \quad \mbox{and}\quad  \int_F \bold{u}_k^\ast \dd s = \int_{F} \bold{\Pi}_k \bold{w} \dd s = \int_F \bold{w}  \dd s = \int_F \bold{g} \dd s,\quad \forall F\in \mathcal {F}_k \setminus \mathcal{F}_k(\Omega),$$
by letting $\bold{\varphi}_F = \int_F \bold{\varphi} \dd s / | F |$ for $F\in \mathcal{F}_k $ and $\bold{\varphi}_T = \int_T \bold{\varphi} \dd x / | T |$ for $T \in \mathcal{T}_k $, we get for any $\bold{\varphi}\in C^\infty(\overline{\Omega})^{d\times d}$,
\begin{align}
    &\quad \bigg| \sum_{F\in\mathcal{F}_k} ([\bold{u}_k^\ast], \bold{\varphi}\cdot\bold{n}_F)_{L^2(F)} - (\bold{g}, \bold{\varphi} \cdot \bold{n})_{L^2(\partial\Omega)}  \bigg| \nonumber\\
    &= \bigg| \sum_{F \in \mathcal{F}_k(\Omega)}( [\bold{u}_k^\ast], \bold{\varphi}\cdot\bold{n}_F )_{L^2(F)} + \sum_{F\in \mathcal{F}_k\setminus\mathcal{F}_k(\Omega)} (\bold{u}_k^\ast - \bold{g}, \bold{\varphi} \cdot \bold{n})_{L^2(F)}  \bigg| \nonumber \\
    &= \bigg| \sum_{F\in\mathcal{F}_k(\Omega)}([\bold{u}_k^\ast] , (\bold{\varphi} - \bold{\varphi}_{F}) \cdot \bold{n}_F)_{L^2(F)} + \sum_{F\in\mathcal{F}_k\setminus\mathcal{F}_k(\Omega)}(\bold{u}_k^\ast - \bold{g} ,(\bold{\varphi} - \bold{\varphi}_F) \cdot \bold{n})_{L^2(F)} \bigg| \nonumber \\
    & \leq \sum_{F\in\mathcal{F}_k(\Omega)} \|[\bold{u}_k^\ast]\|_{\bold{L}^2(F)} \|\bold{\varphi} - \bold{\varphi}_{F}\|_{\mathbb{L}^2(F)} + \sum_{F\in\mathcal{F}_k \setminus \mathcal{F}_k(\Omega)} \|\bold{u}_k^\ast -\bold{g}\|_{\bold{L}^2(F)} \|\bold{\varphi} - \bold{\varphi}_{F}\|_{\mathbb{L}^2(F)}  \nonumber \\
    & \leq \frac{1}{2} \sum_{F\in\mathcal{F}_k(\Omega)} \sum_{T\in\mathcal{N}_k(F)} \|[\bold{u}_k^\ast]\|_{\bold{L}^2(F)} \|\bold{\varphi} - \bold{\varphi}_{T}\|_{\mathbb{L}^2(F)} +     \sum_{F\in\mathcal{F}_k \setminus \mathcal{F}_k(\Omega)} \|\bold{u}_k^\ast -\bold{g}\|_{\bold{L}^2(F)} \|\bold{\varphi} - \bold{\varphi}_{T}\|_{\mathbb{L}^2(F)}  \nonumber \\
    &\leq c \left( \sum_{F\in\mathcal{F}_k(\Omega)} h_F \|[\bold{u}_k^\ast]\|_{\bold{L}^2(F)}^2 + \sum_{T\in \mathcal{F}_k \setminus \mathcal{F}_k(\Omega)} h_F \|\bold{u}_k^\ast -\bold{g}\|_{\bold{L}^2(F)}^2 \right)^{1/2} \|\bold{\nabla}\bold{\varphi}\|_{\mathbb{L}^2(\Omega)}.\label{estimate-pf03_state_G1}
\end{align}
Likewise, we define the element residual on each $T\in \cT_k$ and the jump residual on $F\in \mathcal{F}_k$ respectively by
\[
    R_{T,2}(\phi_k^\ast,\bold{u}_k^\ast) : = \alpha_\eps(\phi_k^\ast) \bold{u}_k^\ast - \bold{f}\quad \mbox{and} \quad J_{F,2}(\bold{u}_k^\ast) : = \left\{ \begin{array}{l}
    [\bold{u}_k^\ast] , \quad \forall F\in \mathcal{F}_k(\Omega),\\ [2ex]
    \bold{u}_k^\ast - \bold{g} , \quad F\in \mathcal{F}_k \setminus \mathcal{F}_k(\Omega).
     \end{array} \right.
\]
From \eqref{estimate-pf01_state_G1} and \eqref{estimate-pf02_state_G1}, the second global estimator $\eta_{k,2}(\phi_k^\ast,\bold{u}_k^\ast;\cT_k)$ for \eqref{opt-sys_state_G1} is defined by
\begin{align*}
    \eta_{k,2}(\phi_k^\ast,\bold{u}^\ast;T)&:= \bigg( h_T^2\|R_{T,2}(\phi_k^\ast,\bold{u}^\ast_k)\|_{\bold{L}^2(T)}^2 + \sum_{F\subset\partial T \cap \Omega}\tfrac{1}{2}h_F \|J_{F,2}(\bold{u}_k^\ast)\|_{\bold{L}^2(F)}^2 + \sum_{F\subset\partial T \cap \partial \Omega}h_F \|J_{F,2}(\bold{u}_k^\ast)\|_{\bold{L}^2(F)}^2\bigg)^{1/2},\\
        \eta_{k,2}(\phi_k^\ast,\bold{u}_k^\ast;\mathcal{M}) &: =  \bigg(\sum_{T\in \mathcal{M}}\eta^2_{k,2}(\phi_k^\ast,\bold{u}_k^\ast;T)\bigg)^{1/2},\quad  \forall \mathcal{M}\subseteq\cT_k.
\end{align*}
Now we can present an adaptive algorithm for problem \eqref{min_phase-field}-\eqref{vp_stokes_div-free}. The set $\mathcal{M}$ and related FE functions will be dropped in the estimators if $\mathcal{M}=\cT_k$ or no confusion arises.

\begin{algorithm}
\caption{AFEM for problem \eqref{min_phase-field} subject to \eqref{vp_stokes_div-free}.}\label{alg_anfem_topopt-Stokes_phase-field_G1}
  \LinesNumbered
  \KwIn{Specify an initial mesh $\cT_{0}$, fix $\eps>0$ and $\gamma>0$ and set the maximum number $K$ of refinement steps.}
  \For {$k=0:K-1$}{
  {(\textsf{OPTIMIZE})} Solve problem \eqref{dismin_phase-field}-\eqref{disvp_stokes_div-free} over $\cT_{k}$ for $(\phi_k^\ast,\bold{u}_{k}^{\ast})\in \U_k \times\bold{U}_{k} $.

  {(\textsf{ESTIMATE})} Compute two error estimators $\eta_{k,1}(\phi_{k}^{\ast},\bold{u}_{k}^{\ast})$ and
  $\eta_{k,2}(\phi_{k}^{\ast},\bold{u}_{k}^{\ast})$.

  {(\textsf{MARK})} Mark two subsets $\mathcal{M}_k^j\subseteq\cT_k$ ($j=1,2$), each containing at least one element $T_{k}^j$ ($j =1,2$) satisfying
    \begin{equation}\label{marking_G1}		
    \eta_{k,j}(T_k^{j})=\max_{T\in\cT_k}\eta_{k,j}(T).
	\end{equation}
    Then $\mathcal{M}_k:=\mathcal{M}_{k}^1\bigcup\mathcal{M}_k^2$.

    {(\textsf{REFINE})} Refine each element $T\in \mathcal{M}_{k}$ by bisection to get $\cT_{k+1}$.

   }

   \KwOut{($\phi_k^\ast ,\bold{u}_k^\ast)$.}
\end{algorithm}

Setting a maximum number of adaptive loops in Algorithms \ref{alg_anfem_topopt-Stokes_phase-field_G1} allows terminating the algorithm in a finite number of finite steps. More frequently, either a prescribed tolerance for $\sum_{j=1}^2\eta_{k,j}^2$ after the module \textsf{ESTIMATE} or a given upper bound on the number of vertices associated with $\cT_k$ after the module \textsf{REFINE} can be used to  stop the iteration in practical implementation. Moreover, the module \textsf{MARK} employs a separate-collective marking scheme. In the first phase, a single subset is determined by a general yet reasonable assumption with respect to the associated estimator, which selects at least one element holding the largest error indicator and is fulfilled by the maximum strategy, the equi-distribution strategy, the
modified equi-distribution strategy and D\"{o}rfler's strategy \cite{NochettoSiebertVeeser:2009, Siebert:2011}. Then the union of the resulting subsets yields the marked set for further refinements.

\section{Convergence of adaptive algorithm}\label{sect:conv_adaptive}

First, we  split the sequence $\{\cT_k\}_{k\geq 0}$ generated by Algorithm \ref{alg_anfem_topopt-Stokes_phase-field_G1} into two classes \cite{Siebert:2011}
\[
    \mathcal{T}_{k}^{+}:=\bigcap_{l\geq k}\mathcal{T}_{l},\quad
    \mathcal{T}_{k}^{0}:=\mathcal{T}_{k}\setminus\mathcal{T}_{k}^{+},\quad
    \Omega_{k}^{+}:=\bigcup_{T\in\mathcal{T}^{+}_{k}}T,\quad
    \Omega_{k}^{0}:=\bigcup_{T\in\mathcal{T}^{0}_{k}}T.
\]
The set $\mathcal{T}_{k}^{+}$ consists of all elements not refined after the $k$-th iteration and
all elements in $\mathcal{T}_{k}^{0}$ are refined at least once after the $k$-th iteration. Moreover, since $\{ \cT_k\}_{k\geq 0}$ is a sequence of locally refined meshes, the mesh-size function $h_k$ has the following property \cite{MorinSiebert:2008,Siebert:2011}
\begin{equation}\label{mesh-size->zero}
        \lim_{k\rightarrow\infty}\|h_{k}\chi^{0}_{k}\|_{L^\infty(\Omega)}=0,
\end{equation}
where $\chi^{0}_{k}$ is the characteristic function of the set $\Omega_{k}^{0}$. For the boundary condition $\bold{u} = \bold{g}$ on $\partial\Omega$, we use an interpolant $\bold{\Pi}_k\bold{w}$ in \eqref{disvp_stokes_div-free} to approximate the lifting $\bold{w}$ of $\bold{g}$. The property $\lim_{k\to\infty}\|\bold{w} - \bold{\Pi}_k\bold{w}\|_{\bold{L}^2(\Omega)} + \|\bold{w} - \bold{\Pi}_k\bold{w}\|_{1,k} =0 $
generally fails over a sequence of adaptively generated meshes, since AFEM is intended to approximate specific function(s), i.e., a minimizing pair $(\phi^\ast, \bold{u}^\ast)$ of \eqref{min_phase-field}-\eqref{vp_stokes_div-free}. So in the analysis, we take $\bold{g} = \bold{0}$.

The convergence of Algorithm \ref{alg_anfem_topopt-Stokes_phase-field_G1} consists of two results:
$\{\eta_{k,1} + \eta_{k,2} \}_{k\geq 0}$ vanishes in the limit up to a subsequence, and the subsequence of discrete minimizing pairs $\{(\phi_{k_j}^\ast,\bold{u}_{k_j}^\ast)\}_{j\geq 0}$ converge to a pair solving \eqref{opt-sys_G1}. These are given in Sections \ref{subsect:est->zero} and \ref{subsect:conv_phase&velocity}. We first prove $\eta_{k,2}\to 0$ in Theorem \ref{thm:conv_est_state&costate}, by which, the discrete compactness of adaptively generated sequences $\{\bold{u}_{k}^\ast\}_{k\geq 0}$ and
the strong continuity of a discrete solution map $\S_{k_j}: \phi_{k_j}^\ast \mapsto \bold{u}^{\ast}_{k_j}$ with respect to the $H^1(\Omega)$-topology of $\{\phi_{k_j}^\ast\}_{j\geq0}$ are proved in Lemmas \ref{lem:vel&costate_weak_adaptive}-\ref{lem:sol-map_cont_limit}. Then
the limit pair of $\{(\phi_{k_j}^\ast,\bold{u}_{k_j}^\ast)\}_{j\geq0}$ is shown to solve \eqref{opt-sys_state_G1}. By combining the technical results in Lemmas \ref{lem:max-est_vi->zero}-\ref{lem:est_reduction_refined} and the ideas in \cite{GantnerPraetorius:2022,MorinNochettoSiebert:2000,Siebert:2011}, we get a null subsequence of the estimators $\{\eta_{k_j,1}\}_{j\geq0}$ in Theorem \ref{thm:conv_est_control}, which is sufficient for \eqref{opt-sys_disc_control_G1} to hold in Theorem \ref{thm:conv_adaptive}. Moreover, Algorithm \ref{alg_anfem_topopt-Stokes_phase-field_G1} yields a sequence of discrete pressure fields $p_{k}^\ast$. We prove that it has an $L^2(\Omega)$ convergent subsequence in Theorem \ref{thm:conv_pre_adaptive}.

With $\bold{g}=\bold{0}$, $\bold{U} = \bold{V}$ in \eqref{vp_stokes_div-free} and $\bold{V}_k + \bold{\Pi}_k \bold{w}$ reduces to $\bold{V}_k$ in \eqref{disvp_stokes_div-free} over $\cT_k$. This, the discrete Poincar\'{e} inequality \eqref{disc_Poin_ineq} and the non-negativity of $\alpha_\eps$ lead to
\begin{equation}\label{stab_discvp-hom_adaptive}
    \|\bold{u}_k\|_{\bold{L}^2(\Omega)} + \|\bold{u}_k\|_{1,k} \leq c \|\bold{f}\|_{\bold{L}^2(\Omega)}.
\end{equation}
The discrete compactness of $\{\bold{u}_k^\ast\}_{k\geq0}$ depends on a vectorial version of the connection operator adapted for the zero boundary condition between the nonconforming $P_1$ FE space and the conforming $P_d$ FE space \cite{Brenner:2003}.
By letting
$$\bold{X}_k^c: = \{\bold{v}\in
\bold{H}^{1}(\Omega)~|~\bold{v}|_{T}\in P_{d}(T)^d,~\forall T\in\mathcal{T}_k\}\quad \mbox{and}\quad \bold{Z}_k^c:= \bold{X}_k^c\cap \bold{H}_0^1(\Omega),$$
we define a connection operator $\E_k: \bold{Z}_k \mapsto \bold{Z}_k^c$ by
\begin{equation}\label{def:con_op_hom}
   \E_{k} \bold{v}(p)= \left\{\begin{array}{ll}
        \dfrac{1}{\#\mathcal{N}_{p}}\displaystyle{\sum_{T\in\mathcal{N}_{p}}}\bold{v}|_{T}(p),&\quad\hbox{if $p$ is an interior node},\\
        \bold{0},&\quad\hbox{if $p$ is a boundary node},
    \end{array}\right.
\end{equation}
where $\mathcal{N}_{p}$ is the set of $T\in\mathcal{T}_k$ that share a common node $p$ and $\#\mathcal{N}_p$ is the cardinality of $\mathcal{N}_{p}$. Also we use a vectorial version of the nodal interpolation operator \cite{Brenner:2003} denoted by  $\bold{\mathcal{I}}_{k}: \bold{Z}_k^c \rightarrow \bold{Z}_{k}$ with
\begin{equation}\label{def:int_nod_hom}
     \bold{\mathcal{I}}_{k}\bold{v}(p)=\bold{v}(p),\quad\hbox{$p$ is the center of a face/edge in $\mathcal{F}_k$}.
\end{equation}
By the argument of \cite[Lemma 3.2 and Corollary 3.3]{Brenner:2003}, the following stability estimate holds
\begin{equation}\label{con_operator-stab_hom}
    |\E_k\bold{v}|_{\bold{H}^1(\Omega)}\leq c\|\bold{v}\|_{1,k},\quad\forall \bold{v}\in \bold{Z}_k,
\end{equation}
with $c$ depending only the shape regularity of $\mathcal{T}_k$.

We give two useful preliminary results. The first is a discrete version of Sobolev inequality on $\bold{Z}_k$. The argument follows the idea of connecting $\bold{Z}_k$ to its conforming counterpart $\bold{Z}_k^c$ and also using $\bold{\mathcal{I}}_k$.
\begin{lemma}\label{lem:disc_Sob_ineq_hom}
    There holds for $2\leq p < \infty$ when $d=2$ or $2 \leq p \leq 6$ when $d=3$
    \begin{equation}\label{disc_Sob_ineq_hom}
        \|\bold{v}\|_{\bold{L}^p(\Om)}\leq c_{\mathrm{ds}}  \|\bold{v}\|_{1,k}, \quad\forall\bold{v}\in \bold{Z}_k,
    \end{equation}
    with the constant $c_{\mathrm{ds}}$ only depending on the shape regularity of $\cT_k$.
\end{lemma}
\begin{proof}
It follows directly from \eqref{def:con_op_hom} and \eqref{def:int_nod_hom} that $\bold{\mathcal{I}}_k \E_{k} \bold{v} = \bold{v}$ for any $\bold{v} \in \bold{Z}_k$. The standard interpolation error estimate \cite{BrennerScott:2008, Ciarlet:2002} and the inverse estimate give that for any $T\in \cT_k$,
\begin{align*}
    \| \bold{v}- \E_k\bold{v} \|^{p}_{\bold{L}^p(T)} & = \| \E_k \bold{v} - \bold{\mathcal{I}}_{k} \E_k\bold{v} \|^p_{\bold{L}^p(T)} \leq
    c h_{T}^{2p - \frac{dp}{2} + d} |\E_k \bold{v}|_{\bold{H}^2(T)}^p \leq
    c h_{T}^{p - \frac{dp}{2} + d}|\E_k \bold{v}|_{\bold{H}^1(T)}^p \\
    &\leq c h_{T}^{p - \frac{dp}{2} + d}h_{T}^{- p + \frac{dp}{2} - d}\|\E_k \bold{v}\|_{\bold{L}^p(T)}^p=
    c \|\E_k \bold{v}\|_{\bold{L}^p(T)}^p.
\end{align*}
Summing up the estimate over $T\in \mathcal{T}_k$ and the triangle inequality give
\[
    \|\bold{v}\|_{\bold{L}^p(\Om)}\leq \|\bold{v} - \E_k \bold{v}\|_{\bold{L}^p(\Om)}+\|\E_k \bold{v}\|_{\bold{L}^p(\Om)}
    \leq c\|\E_k \bold{v}\|_{\bold{L}^p(\Om)},\quad\forall \bold{v}\in \bold{Z}_k,
\]
which, together with Sobolev embedding theorem, Poincar\'{e} inequality and the estimate \eqref{con_operator-stab_hom}, yields
\[
    \|\bold{v}\|_{\bold{L}^p(\Om)}\leq c\|\E_k \bold{v}\|_{\bold{L}^p(\Om)}\leq c\|\E_k \bold{v}\|_{\bold{H}^1(\Omega)}\leq c \|\bold{v}\|_{1,k}.
\]
Thus we obtain the desired estimate \eqref{disc_Sob_ineq_hom}.
\end{proof}

Next we present an estimate for $\E_k$ in terms of the jump term in the estimator $\eta_{k,2}$.
\begin{lemma}\label{lem:est_noncon}
For the operator $\E_k$, there exists $c$ depending only on the shape regularity of $\mathcal{T}_k$ such that
\begin{equation}\label{est_noncon}
    \|\bold{v}-\E_k\bold{v}\|_{\bold{L}^{2}(\Omega)}^{2}\leq c\sum_{F\in\mathcal{F}_{k}}h_F\|[\bold{v}]\|_{\bold{L}^{2}(F)}^{2},\quad\forall \bold{v}\in \bold{Z}_k.
    \end{equation}
\end{lemma}
\begin{proof}
The proof follows the idea of \cite[(3.3), Lemma 3.2]{Brenner:2003}. Let $\mathcal{C}_T$ be the set of centers of the faces/edges of any $T\in\mathcal{T}_k$ and
$\mathcal{N}_T$ be the set of other nodes of $T$.
For any interior node $p\in\mathcal{N}_T\cap\Omega$ and any $T'\in\mathcal{N}_{p}$, there exists a
sequence of elements $T_{1}=T,\ldots,T_{m}=T'$ in $\mathcal{N}_p$ with each consecutive pair
$T_{i}$ and $T_{i+1}$ sharing a face, where the integers $\#\mathcal{N}_p$ and $m$ are both uniformly bounded by a
constant depending on the shape regularity of $\mathcal{T}_k$. Then by the inverse estimate, we have
\[
    |\bold{v}|_{T}(p)-\bold{v}|_{T'}(p)|^{2}\leq c\sum_{i=1}^{m-1}|\bold{v}|_{T_i}(p)-\bold{v}|_{T_i+1}(p)|^{2}
    \leq c\sum_{F\in\mathcal{V}(p)}h^{1-d}_{F}\|[\bold{v}]\|_{\bold{L}^{2}(F)}^{2},
\]
where $\mathcal{V}(p)$ is the set of faces/edges $F$ with $p$ as a common node on $\partial F$. Thus there holds
\[
    \left|\bold{v}|_{T}(p)-\E_{k}\bold{v}(p)\right|^{2}
    =\bigg|\bold{v}|_{T}(p)-\frac{1}{\#\mathcal{N}_{p}}\sum_{T'\in\mathcal{N}_{p}}\bold{v}|_{T'}(p)\bigg|^{2}
    \leq c\sum_{F\in\mathcal{V}(p)}h^{1-d}_{F}\|[\bold{v}]\|_{\bold{L}^{2}(F)}^{2}.
\]
For any boundary node $p\in\mathcal{N}_T\cap\partial\Om$, it can be shown similarly that for $p\in F=\partial T\cap\partial\Omega$,
\[
    |\bold{v}|_{T}(p)-\E_{k}\bold{v}(p)|^{2}=|\bold{v}|_{T}(p)|^{2}\leq ch^{1-d}_{F}\|[\bold{v}]\|_{\bold{L}^{2}(F)}^{2}.
\]
On any $T\in\mathcal{T}_k$, since there is no contribution from $\mathcal{C}_T$, by the norm equivalence, the definition of $\E_{k}$ and the local quasi-uniformity of $\mathcal{T}_k$, we get
    \begin{align}
        \|\bold{v}-\E_{k}\bold{v}\|_{\bold{L}^{2}(T)}^{2}
        &\leq ch_{T}^{d}\sum_{p\in\mathcal{N}_T\cap\overline{\Omega}}|\bold{v}(p)- \E_{k}\bold{v}(p)|^{2}\nonumber\\
        &=ch_{T}^{d}\sum_{p\in\mathcal{N}_T\cap\Omega}\bigg|\bold{v}|_{T}(p)-\frac{1}{|\omega_{p}|}\sum_{T'\in\omega_{p}}\bold{v}|_{T'}(p)\bigg|^{2}
        +h_{T}^{d}\sum_{p\in\mathcal{N}_T\cap\partial\Omega}|\bold{v}|_{T}(p)|^{2}\nonumber\\
        &\leq c\bigg( h_T^d\sum_{p\in\mathcal{N}_T\cap\Omega}\sum_{F\in\mathcal{V}(p)}h_{F}^{1-d}\|[\bold{v}]\|_{\bold{L}^{2}(F)}^{2}
        +h_T^d\sum_{p\in\mathcal{N}_T\cap\partial\Omega}\sum_{F\subset\partial T}h_{F}^{1-d}\|[\bold{v}]\|_{\bold{L}^{2}(F)}^{2}\bigg)\nonumber\\
        &\leq c\bigg( \sum_{p\in\mathcal{N}_T\cap\Omega}\sum_{F\in\mathcal{V}(p)}h_{F}\|[\bold{v}]\|_{\bold{L}^{2}(F)}^{2}
        +\sum_{p\in\mathcal{N}_T\cap\partial\Omega}\sum_{F\subset\partial T}h_{F}\|[\bold{v}]\|_{\bold{L}^{2}(F)}^{2}\bigg).\label{lem:est_noncon_pf1}
    \end{align}
    Then summing up the estimate \eqref{lem:est_noncon_pf1} over all $T\in\mathcal{T}_k$ gives the desired estimate.
\end{proof}

\subsection{Zero limit of estimators}\label{subsect:est->zero}
Algorithm \ref{alg_anfem_topopt-Stokes_phase-field_G1} is driven by two estimators $\eta_{k,1}$ and $\eta_{k,2}$ via suitable requirement on the marking strategy \eqref{marking_G1}. The limiting behavior of the maximal error indicators among all elements in each $\cT_k$ is central to the convergence analysis \cite{Siebert:2011}. We first analyze the property of the estimator $\eta_{k,2}$ and then $\eta_{k,1}$.
\begin{lemma}\label{lem:max-est->zero}
Let $\{(\phi_{k}^{\ast},\bold{u}_{k}^{\ast})\}_{k\geq0}$ be the sequence of discrete solutions generated by Algorithm \ref{alg_anfem_topopt-Stokes_phase-field_G1}  and $\{\mathcal{M}_k\}_{k\geq 0}$ the corresponding sequence of marked sets determined by \eqref{marking_G1}. Then there holds
\begin{equation}\label{conv:max-est->zero}
        \lim_{k\to\infty} \max_{T\in\mathcal{M}_{k}}\eta_{k,2}(\phi_{k}^{\ast},\bold{u}_{k}^{\ast};T) =0.
    \end{equation}
\end{lemma}
\begin{proof}
Let $T_{k}^2$ denote any element with the largest error indicator among $\eta_{k,2}(\phi_{k}^\ast, \bold{u}^\ast_{k_j};T)$
over $\mathcal{M}_{k}$. By \eqref{disc_Poin_ineq} and \eqref{stab_discvp-hom_adaptive}, we have
\begin{equation*}
    h_{T_{k}^2}^2\|\alpha_\eps(\phi_{k}^{\ast})\bold{u}_{k}^{\ast}-\bold{f}\|_{\bold{L}^2(T_{k}^2)}^2
    \leq 2 h_{T_{k}^2}^2 (\|\bold{f}\|^2_{\bold{L}^2(\Omega)}+\overline{\alpha}_\eps^2\|\bold{u}_{k}^{\ast}\|^2_{\bold{L}^2(\Omega)})\leq
    c h_{T_{k}^2}^2 \|\bold{f}\|^2_{\bold{L}^2(\Omega)}.
\end{equation*}
Since $[\bold{u}_{k}^{\ast}] = \bold{0}$ at the center of $F\in\mathcal{F}_k$ (cf. the definition of $\bold{Z}_k$) and $\bold{\nabla}_k \bold{u}_k^{\ast}$ is a piecewise constant tensor over $\cT_k$, by the scaled trace theorem, the local quasi-uniformity of $\cT_k$ and \eqref{stab_discvp-hom_adaptive} again, we deduce that for $d=3$,
\begin{align*}
        \sum_{F\subset\partial T_{k}^2}h_{F}\int_F|[\bold{u}^{\ast}_{k}]|^{2} \ds &\leq c\sum_{F\subset\partial T^2_{k}}h_F^{3}\int_{F}|[\boldsymbol{\nabla}_{k} \bold{u}^{\ast}_{k}\times\boldsymbol{n}_{F}]|^{2}\ds\leq c h_{T^2_{k}}^{2}\sum_{T\in\mathcal{N}_k(T_{k}^2)}\int_T|\boldsymbol{\nabla}_{k}\bold{u}^{\ast}_{k}|^{2}\dx\nonumber\\
        &\leq c h_{T^{2}_{k}}^2\|\bold{u}^{\ast}_{k_j}\|_{1,k}^{2}\leq ch_{T^2_{k}}^2 \|\bold{f}\|_{\bold{L}^2(\Omega)}^2,
    \end{align*}
    and similarly, with $\bold{t}_F$ being a unit tangent vector by rotating $\bold{n}_F$ $90$ degrees counter-clockwise, for $d=2$
\begin{equation*}
    \sum_{F\subset\partial T_{k}^2}h_{F}\int_F|[\bold{u}^{\ast}_{k}]|^{2} \ds \leq c\sum_{F\subset\partial T^2_{k}}h_F^{3}\int_{F}|[\boldsymbol{\nabla}_{k} \bold{u}^{\ast}_{k} \cdot \boldsymbol{t}_{F}]|^{2}\ds\leq c h_{T^{2}_{k}}^2\|\bold{u}^{\ast}_{k}\|_{1,k}^{2}\leq ch_{T^2_{k}}^2 \|\bold{f}\|_{\bold{L}^2(\Omega)}^2.
    \end{equation*}
 Since $T_{k}^2\in \mathcal{M}_{k}^2 \subset  \mathcal{M}_{k} \subset \cT_{k}^0$, we deduce from \eqref{mesh-size->zero} that
    \begin{equation}\label{conv:max-est_pf4}
       \lim_{k\to\infty} h_{T_{k}^2} \leq \lim_{k\to\infty}\|h_{k}\chi_{k}^0\|_{L^\infty(\Omega)}= 0.
    \end{equation}
Hence the assertion in \eqref{conv:max-est->zero} follows from the preceding estimates.
\end{proof}

\begin{theorem}\label{thm:conv_est_state&costate}
The sequence $\{\eta_{k,2}(\phi_{k}^\ast,\bold{u}_k^\ast)\}_{k\geq 0}$  generated by Algorithm \ref{alg_anfem_topopt-Stokes_phase-field_G1} converges to zero.
\end{theorem}
\begin{proof}
The proof is inspired by \cite{Siebert:2011}. Since $\cT_{l}^+ \subset \cT_k^+$ for $k>l$, one may split $\eta_{k,2}^2(\phi_k^\ast,\bold{u}_k^\ast)$ into
\begin{equation*}
\eta_{k,2}^2(\phi_k^\ast,\bold{u}_k^\ast) = \eta_{k,2}^2(\phi_k^\ast,\bold{u}_k^\ast; \cT_{k} \setminus \cT_{l}^+ ) + \eta_{k,2}^2(\phi_k^\ast,\bold{u}_k^\ast; \cT_{l}^+ ).
\end{equation*}
 The argument for Lemma \ref{lem:max-est->zero} gives
 \begin{align*}
        \eta_{k,2}^2(\phi_k^\ast,\bold{u}_k^\ast; \cT_{k} \setminus \cT_{l}^+ ) & \leq \sum_{T\in \cT_{k} \setminus \cT_{l}^+  } \bigg( h_{T}^2\|\alpha_\eps(\phi_{k}^{\ast})\bold{u}_{k}^{\ast} - \bold{f}\|_{\bold{L}^2(T)}^2 + \sum_{F\subset \partial T} h_F \| [ \bold{u}_k ]\|_{\bold{L}^2(F)}^2 \bigg) \\
        & \leq c \bigg( \max_{T\in \cT_{k}\setminus\cT_l^+}h_T^2 \|\alpha_\eps(\phi_{k}^{\ast})\bold{u}_{k}^{\ast} - \bold{f}\|_{\bold{L}^2(\Omega)}^2  + \sum_{T\in\cT_{k}\setminus\cT_l^+}h_T^2\sum_{T'\in\omega_k(T)}\int_{T'}|\bold{\nabla}_{k}\bold{u}_{k}^\ast|^2\dx\bigg) \\
        & \leq c \max_{T\in \cT_{k}\setminus\cT_l^+}h_T^2 \left( \|\alpha_\eps(\phi_{k}^{\ast})\bold{u}_{k}^{\ast} - \bold{f}\|_{\bold{L}^2(\Omega)}^2 + \|\bold{u}^\ast_k\|_{1,k} \right) \leq c \max_{T\in \cT_{k}\setminus\cT_l^+}h_T^2 \|\bold{f}\|_{\bold{L}^2(\Omega)}^2.
\end{align*}
Since $ \bigcup_{T\in\cT_{k}\setminus\cT_{l}^{+}}T = \bigcup_{T\in\cT_{l}\setminus\cT_{l}^{+}}T$, the monotonicity of $h_k$ and \eqref{mesh-size->zero} further imply
$\lim_{l\to\infty} \max_{T\in \cT_{k}\setminus\cT_l^+}h_T^2 \leq \lim_{l\to\infty}\|h_l\chi_l^0\|_{L^\infty(\Omega)}= 0.$
    So when $l$ is sufficiently large, $\eta_{k,2}^2(\phi_k^\ast,\bold{u}_k^\ast; \cT_{k} \setminus \cT_{l}^+ )$ is smaller than any given positive number. Using $\mathcal{M}_k\subset \cT_k^0$ and $\cT_l^+ \subset \cT_k^+ \subset \cT_k$ for $k>l$ and the marking assumption \eqref{marking_G1}, we obtain
    \[
        \eta_{k,2}^2(\phi_k^\ast,\bold{u}_k^\ast; \cT_{l}^+ )\leq \# \cT_l^+ \max_{T\in\cT_l^+} \eta_{k,2}^2(\phi_k^\ast,\bold{u}_k^\ast; T )\leq \# \cT_l^+ \max_{T\in\mathcal{M}_k} \eta_{k,2}^2(\phi_k^\ast,\bold{u}_k^\ast; T ).
    \]
    With the index $l$ fixed, the right hand side tends to zero as $k\to\infty$ by Lemma \ref{lem:max-est->zero}. Therefore, we may choose an $L\in\mathbb{N}$ such that $\eta_{k,2}^2(\phi_k^\ast,\bold{u}_k^\ast; \cT_{l}^+ )$ is also sufficiently small for all $k>l\geq L$. Therefore, $\lim_{k\to\infty}\eta_{k,2}(\phi_{k}^\ast,\bold{u}_k^\ast)= 0$.
\end{proof}

Next we analyze the convergence of the sequence $\{\eta_{k,1}(\phi_{k}^\ast,\bold{u}_k^\ast))\}_{k\geq0}$. First we extract a weakly convergent subsequence from $\{\bold{u}_k^\ast\}_{k\geq0}$ in $\bold{L}^2(\Omega)$.
\begin{lemma}\label{lem:vel&costate_weak_adaptive}
Let $\{\bold{u}_{k}^\ast\}_{k\geq 0}$ be the sequence of discrete velocity fields generated by Algorithm \ref{alg_anfem_topopt-Stokes_phase-field_G1}. Then there exists a subsequence $\{\bold{u}_{k_j}^\ast\}_{j \geq 0}$ and some $\overline{\bold{u}} \in \bold{V}$ such that
\begin{align}\label{vel_weak_adaptive}
   \bold{u}_{k_j}^\ast \rightarrow  \overline{\bold{u}} \quad \text{weakly in}~\bold{L}^2(\Omega)\quad\mbox{and}\quad
        \bold{\nabla}_{k_j}\bold{u}_{k_j}^\ast\rightarrow   \bold{\nabla} \overline{\bold{u}}\quad \text{weakly in}~\mathbb{L}^2(\Omega).
\end{align}
\end{lemma}
\begin{proof}
The proof uses Theorem \ref{thm:conv_est_state&costate} and consists of three steps dealing with the following three assertions (1) $\overline{\bold{u}}\in \bold{H}^1(\Omega)$; (2) $\overline{\bold{u}} = \bold{0}$ on $\partial\Omega$; (3) $\mathrm{div} \overline{\bold{u}} = 0$, separately.

\noindent \textit{Step 1.} $\overline{\bold{u}}\in \bold{H}^1(\Omega)$. Note that each $\bold{u}_k^\ast\in \bold{V}_k$ solves \eqref{disvp_stokes_div-free} over $\cT_k$. So by the estimate \eqref{stab_discvp-hom_adaptive}, the sequence $ \{\|\bold{u}^\ast_k\|_{\bold{L}^2(\Omega)} + \|\bold{u}^\ast_k\|_{1,k}\}$ is uniformly bounded. Then there exist two subsequences, denoted by $\{\bold{u}_{k_j}^\ast\}_{j\geq0}$ and $\{\bold{\nabla}_{k_j}\bold{u}_{k_j}^\ast\}_{j\geq0}$, and two weak limits $\overline{\bold{u}}\in \bold{L}^2(\Omega)$ and $\bold{\sigma} \in \mathbb{L}^{2}(\Omega)$ such that
    \begin{equation}\label{lem:vel_weak_adaptive_pf01}
        \bold{u}_{k_j}^\ast \rightarrow \overline{\bold{u}}\quad \mbox{weakly in }\bold{L}^2(\Omega)\quad\mbox{and}\quad \bold{\nabla}_{k_j}\bold{u}_{k_j}^\ast \rightarrow \bold{\sigma} \quad\mbox{weakly in }\mathbb{L}^2(\Omega).
    \end{equation}
It suffices to show $\overline{\bold{u}}\in \bold{H}_0^1(\Omega)$ with $\bold{\sigma} = \bold{\nabla}\overline{\bold{u}}$ and $\mathbf{div}\overline{\bold{u}} = 0$. First, by elementwise integration by parts, we obtain \begin{equation}\label{lem:vel_weak_adaptive_pf02}
(\bold{u}_{k_j}^\ast, \mathbf{div} \bold{\varphi})
         = - \sum_{T\in\cT_{k_j}}( \bold{\nabla} \bold{u}_{k_j}^\ast, \bold{\varphi})_{L^2(T)} + \sum_{F\in\mathcal{F}_{k_j}}( [\bold{u}_{k_j}^\ast],  \bold{\varphi}\cdot\bold{n}_F)_{L^2(F)}, \quad \forall \bold{\varphi} \in C^\infty(\overline{\Omega})^{d\times d}.
    \end{equation}
Since $\bold{u}_{k_j}^\ast \in \bold{V}_{k_j} \subset \bold{Z}_{k_j}$, $\int_F [\bold{u}_{k_j}^\ast] \ds =\bold{0}$ on $F\in\mathcal{F}_{k_j}$. Inserting $\bold{\varphi}_F = \int_F \bold{\varphi} \ds / | F |$ again for all $F\in \mathcal{F}_{k_j}$ and repeating the argument for \eqref{estimate-pf03_state_G1} with $\bold{g}=\bold{0}$, we obtain
    \begin{equation}\label{lem:vel_weak_adaptive_pf03}
        \bigg| \sum_{F\in\mathcal{F}_{k_j}} ([\bold{u}_{k_j}^\ast], \bold{\varphi}\cdot\bold{n}_F)_{L^2(F)} \bigg| = \bigg| \sum_{F\in\mathcal{F}_{k_j}} ([\bold{u}_{k_j}^\ast], (\bold{\varphi} - \bold{\varphi}_F )\cdot\bold{n}_F)_{L^2(F)} \bigg| \leq c \sum_{F\in\mathcal{F}_{k_j}} h_F \|[\bold{u}_{k_j}^\ast]\|_{\bold{L}^2(F)}^2 \|\bold{\nabla}\bold{\varphi}\|_{\mathbb{L}^2(\Omega)}.
        \end{equation}
Passing to the limit $j\to\infty$ in \eqref{lem:vel_weak_adaptive_pf02}, \eqref{lem:vel_weak_adaptive_pf01}, \eqref{lem:vel_weak_adaptive_pf03} and the relation $\lim_{k\to\infty}\eta_{k,2}(\phi_k^\ast, \bold{u}_{k}^\ast)= 0$ in Theorem \ref{thm:conv_est_state&costate} yield
\begin{equation}\label{lem:vel_weak_adaptive_pf04}
        (\overline{\bold{u}}, \mathbf{div} \bold{\varphi} ) = - ( \bold{\sigma}, \bold{\varphi}), \quad \forall \bold{\varphi}\in C^\infty(\overline{\Omega})^{d\times d}.
    \end{equation}
Restricting $\bold{\varphi}$ to $C^\infty_0(\Omega)^{d\times d}$ and the definition of the distributional gradient give
$\bold{\sigma} = \bold{\nabla}\overline{\bold{u}}$ and $\overline{\bold{u}}\in \bold{H}^1(\Omega).$

\noindent \textit{Step 2.} $\overline{\bold u}=\bold{0}$ on $\partial\Omega$. Let $\bold{H}(\mathrm{div}):=\{\boldsymbol{v}\in
\bold{L}^{2}(\Omega),\mathrm{div}\boldsymbol{v}\in L^{2}(\Omega)\}$ and $\mathbb{H}(\mathbf{div}) :=\{\bold{\tau}\in \mathbb{L}^2(\Omega), \mathbf{div}\bold{\tau}\in \bold{L}^2(\Omega)\}$. 
Note the density of $C^\infty(\overline{\Omega})^d$ in $\bold{H}(\mathrm{div})$ and the unique continuous extension of the normal
trace mapping $\gamma_n: C^\infty(\overline{\Omega})^d \to H^{-1/2}(\partial\Omega) $ with $\bold{v} \mapsto \bold{v}
\cdot \bold{n} $ to $\gamma_n: \bold{H}(\mathrm{div}) \to H^{-1/2}(\partial\Omega) $
\cite[Theorems I.2.4 and I.2.5, p.27]{GiraultRaviart:1986}. Applying the density result row-wisely to each
$\bold{\varphi} \in C^\infty(\overline{\Omega})^{d\times d}$ and using \eqref{lem:vel_weak_adaptive_pf04} give
$( \overline{\bold{u}} , \mathbf{div}\bold{\varphi} ) = - ( \bold{\nabla} \overline{\bold{u}} , \bold{\varphi})$ for all $\bold{\varphi} \in \mathbb{H}(\mathbf{div})$.
Applying  \cite[Theorems I.2.4 and I.2.5]{GiraultRaviart:1986} row-wisely also leads to Green's formula
\begin{equation*}
        ( \overline{\bold{u}}, \mathbf{div} \bold{\varphi} ) + ( \bold{\nabla} \overline{\bold{u}}, \bold{\varphi} ) = \langle \bold{\varphi}\cdot\bold{n}, \overline{\bold{u}}\rangle_{H^{-1/2}(\partial\Omega)^d \times H^{1/2}(\partial\Omega)^d}, \quad \forall \bold{\varphi}\in \mathbb{H}(\mathbf{div}).
    \end{equation*}
Consequently,
$\langle \bold{\varphi}\cdot\bold{n}, \overline{\bold{u}}\rangle_{H^{-1/2}(\partial\Omega)^d \times H^{1/2}(\partial\Omega)^d} = 0$ for any $ \bold{\varphi}\in \mathbb{H}(\mathbf{div}).$ This and the surjectivity of $\gamma_n: \bold{H}(\mathrm{div}) \to H^{-1/2}(\partial\Omega)$ \cite[Corollary I.2.8, p.28]{GiraultRaviart:1986} yield $\overline{\bold{u}}=\bold{0}$ on $\partial \Omega$.

\noindent \textit{Step 3.} $\mathrm{div}\overline{\bold{u}}=0$. Since $\bold{g}=0$ and $\bold{u}^\ast_{k_j}\in\bold{V}_{k_j}$, by elementwise integration by parts, we obtain
\begin{equation}\label{lem:vel_weak_adaptive_pf05}
(\bold{u}_{k_j}^\ast, \bold{\nabla} \varphi) = - \sum_{T\in\cT_{k_j}} (  \mathrm{div} \bold{u}_{k_j}^\ast, \varphi)_{L^2(T)} + \sum_{F\in \mathcal{F}_{k_j}} ( [\bold{u}_{k_j}^\ast]\cdot \bold{n}_F, \varphi)_{L^2(F)}, \quad \forall \varphi \in C^\infty(\overline{\Omega}).
    \end{equation}
    Then with $\varphi_{F} = \int_F \varphi \dx / |F|$ for all $F\in \mathcal{F}_{k_j}$ and $\varphi_{T} = \int_T \varphi \dx / |T|$ for all $T\in\cT_{k_j}$, the second term in \eqref{lem:vel_weak_adaptive_pf05} can be estimated by a scaled trace theorem and Poincar\'{e} inequality as
    \begin{align*}
    &\quad \bigg|
         \sum_{F\in\mathcal{F}_{k_j}}([\bold{u}_{k_j}^\ast]\cdot\bold{n}_F ,\varphi)_{L^2(F)}
    \bigg|  = \bigg| \sum_{F\in\mathcal{F}_{k_j}}( [\bold{u}_{k_j}^\ast]\cdot\bold{n}_F ,\varphi-\varphi_F)_{L^2(F)} \bigg| \leq \sum_{F\in\mathcal{F}_{k_j}} \|[\bold{u}_{k_j}^\ast]\cdot\bold{n}_F\|_{L^2(F)} \|\varphi-\varphi_F\|_{L^2(F)} \\
    & \leq \frac{1}{2}\sum_{F\in\mathcal{F}_{k_j}} \|[\bold{u}_{k_j}^\ast]\|_{L^2(F)} \sum_{T\in\mathcal{N}_{k_j}(F)}\|\varphi-\varphi_T\|_{L^2(F)} \leq c \sum_{F\in\mathcal{F}_{k_j}} h_F^{1/2}\|[\bold{u}_{k_j}^\ast]\|_{L^2(F)} \|\bold{\nabla}\varphi\|_{L^2(\omega_{k_j}(F))} \\
    &\leq c\bigg(\sum_{F\in\mathcal{F}_{k_j}} h_F\|[\bold{u}_{k_j}^\ast]\|_{L^2(F)}^2\bigg)^{1/2}\|\bold{\nabla}\varphi\|_{L^2(\Omega)}, \quad \forall \varphi \in C^\infty(\overline{\Omega}),
    \end{align*}
which, together with the identity $\lim_{k\to\infty}\eta_{k,2}(\phi_k^\ast, \bold{u}_{k}^\ast)= 0$ from Theorem \ref{thm:conv_est_state&costate}, implies
\begin{equation}\label{lem:vel_weak_adaptive_pf06}
   \lim_{j\to\infty} \sum_{F\in\mathcal{F}_{k_j}}( [\bold{u}_{k_j}^\ast]\cdot\bold{n}_F, \varphi )_{L^2(F)}= 0, \quad \forall \varphi \in C^\infty(\overline{\Omega}).
\end{equation}
So from $\mathrm{div}_{k_j}\bold{u}_{k_j}^\ast =0 $, \eqref{lem:vel_weak_adaptive_pf01} and \eqref{lem:vel_weak_adaptive_pf06}, passing to the limit $j\to\infty$ on both sides of \eqref{lem:vel_weak_adaptive_pf05} yields that $(\overline{\bold{u}}, \bold{\nabla} \varphi ) = 0$ for any $\varphi \in C_0^\infty(\Omega)$. This concludes the proof of the lemma.
\end{proof}

So far we have shown that $\{\bold{u}_{k}^\ast\}_{j\geq 0}$ and its piecewise gradient have some $L^2$ weak limits up to a subsequence.
The next lemma asserts that the $L^2$ weak limit $\overline{\bold{u}}$ of $\{\bold{u}_{k_j}^\ast\}_{j\geq 0}$ in Lemma \ref{lem:vel&costate_weak_adaptive} is actually an $L^2$ strong limit. This lifting can be viewed as the discrete compactness of nonconforming linear (the lowest-order CR) finite elements over a sequence of locally refined meshes by AFEM.

\begin{lemma}\label{lem:vel&costate_L2-strong-conv_adaptive}
The $\bold{L}^2(\Omega)$ weakly convergent subsequence $\{\bold{u}_{k_j}^\ast\}_{j\geq 0}$ in Lemma \ref{lem:vel&costate_weak_adaptive}
converges strongly to $\overline{\bold{u}}\in \bold{V}$ with respect to the $\bold{L}^2(\Omega)$-topology, i.e.,
 $     \lim_{j\to\infty}\|  \bold{u}_{k_j}^\ast -\overline{\bold{u}}\|_{\bold{L}^2(\Omega)}=0$.
\end{lemma}
\begin{proof}
By applying the connection operator $\E_{k_j}$ to $\bold{u}_{k_j}^\ast$, we deduce from Poincar\'{e} inequality, \eqref{con_operator-stab_hom} and \eqref{stab_discvp-hom_adaptive} that the sequence
$  \{\|\E_{k_j}\bold{u}_{k_j}^\ast\|_{\bold{H}^1(\Omega)}\}$  is uniformly bounded,
which, along with Sobolev compact embedding theorem, further yields a subsequence, still denoted by $\{\E_{k_j}\bold{u}_{k_j}^\ast\}_{j\geq 0}$, and some $\widehat{\bold{u}}\in \bold{H}_0^1(\Omega)$ satisfying
$\lim_{j\to\infty}\|\E_{k_j}\bold{u}_{k_j}^\ast  - \widehat{\bold{u}}\|_{\bold{L}^2(\Omega)} = 0$.
Lemma \ref{lem:est_noncon} and the relation $\lim_{k\to\infty}\eta_{k,2}(\phi_{k}^\ast,\bold{u}_k^\ast)=0$ in Theorem \ref{thm:conv_est_state&costate} imply
\begin{equation*}
      \lim_{j\to\infty}  \| \bold{u}_{k_j}^\ast - \E_{k_j}\bold{u}_{k_j}^\ast \|_{\bold{L}^2(\Omega)} \leq \lim_{j\to\infty} c\bigg(\sum_{F\in\mathcal{F}_{k_j}}h_F\|[\bold{u}_{k_j}^\ast]\|_{\bold{L}^{2}(F)}^{2}\bigg)^{1/2} = 0.
\end{equation*}
Consequently, $\lim_{j\to\infty}\|\bold{u}_{k_j}^\ast - \widehat{\bold{u}} \|_{\bold{L}^2(\Omega)} = 0$. This and the weak limit $\overline{\bold{u}}$ in Lemma \ref{lem:vel&costate_weak_adaptive} lead to the desired assertion.
\end{proof}

\begin{remark}\label{rem:vel&costate_disc-comp_adaptive}
The discrete compactness in Lemma \ref{lem:vel&costate_L2-strong-conv_adaptive} holds for a specific sequence of CR FE approximations, not each sequence of CR FE functions, over the adaptively generated meshes. This agrees with the intuition that AFEM approximates only one single function in $\bold{V}$. More importantly, a jump term associated with the discrete solution is included in Algorithm \ref{alg_anfem_topopt-Stokes_phase-field_G1} as part of an estimator and the sequence is proved to tend to zero as the adaptive algorithm proceeds in Theorem \ref{thm:conv_est_state&costate}. In the proof of Lemma \ref{lem:vel&costate_L2-strong-conv_adaptive}, the connection operator $\E_k$, Sobolev compact embedding theorem and Lemma \ref{lem:est_noncon} connect the compactness property with the null sequence in Theorem \ref{thm:conv_est_state&costate}, which plays the role of $\lim_{k\to\infty}\|h_k\|_{L^\infty(\Omega)}=0$ in the discrete compactness property for uniform mesh refinements \cite{Stummel:1980}.
\end{remark}

Now we analyze $\{\eta_{k,1}(\phi_k^\ast,\bold{u}_k^\ast)\}_{k\geq0}$. We use  the following
limiting variational problem: find $\bold{u}_\infty \in \bold{V}$ such that
\begin{equation}
    \mu (\bold{\nabla} \bold{u}_\infty , \bold{\nabla} \bold{v}) + (\alpha_\eps(\phi_\infty) \bold{u}_\infty, \bold{v})  = (\bold{f} , \bold{v} ), \quad \forall \bold{v} \in \bold{V} \label{vp_stokes_div-free_limit}
\end{equation}
with $\{\U_k\}_{k\geq 0}$ given by Algorithm \ref{alg_anfem_topopt-Stokes_phase-field_G1} and a fixed $\phi_\infty\in \U_\infty : = \overline{\bigcup_{k\geq 0}\U_k}$ in the $H^1(\Omega)$-norm. The solution maps defined by \eqref{disvp_stokes_div-free} over $\cT_k$ and \eqref{vp_stokes_div-free} respectively are denoted by $\S_k:\U_k\to\bold{V}_k$ and $\S:\U\to\bold{V}$. This motivates reformulating the solution $\bold{u}_\infty$ of problem \eqref{vp_stokes_div-free_limit} as $\bold{u}_\infty = \S(\phi_\infty)$. Since Lemmas \ref{lem:vel&costate_weak_adaptive} and \ref{lem:vel&costate_L2-strong-conv_adaptive} involve a sequence of adaptive solutions $\{\bold{u}_k^\ast = \S_k(\phi_k^\ast)\}_{k\geq 0}$, we analyze the continuity of $\S_k$ along with $\{\phi_k^\ast\}_{k\geq0}$.

\begin{lemma}\label{lem:sol-map_cont_limit}
    The sequence of discrete minimizing pairs $\{(\phi_k^\ast,\bold{u}_k^\ast)\}_{k\geq 0}$ generated by Algorithm \ref{alg_anfem_topopt-Stokes_phase-field_G1} has a subsequence converging to a  $\overline{\phi}_\infty \in \U_\infty$  in $H^1(\Omega)$ and the associated velocity field $\overline{\bold{u}}_\infty = \S(\overline{\phi}_\infty) \in \bold{V}$ in the sense that
    \begin{equation}\label{vel_strong-conv_limit}
      \lim_{j\to\infty}  \|\bold{u}^\ast_{k_j} - \overline{\bold{u}}_\infty \|_{\bold{L}^2(\Omega)} + \|\bold{\nabla}_{k_j}(\bold{u}^\ast_{k_j} - \overline{\bold{u}}_\infty) \|_{\bold{L}^2(\Omega)} = 0.
    \end{equation}
\end{lemma}
\begin{proof}
We take $\phi_k = \beta$ in \eqref{dismin_phase-field} over $\cT_k$ for each $k\geq 0$. By the estimate \eqref{stab_discvp-hom_adaptive}, the sequence $\{\|\S_k(\phi_k^\ast)\|_{L^2(\Omega)} + \|\S_k(\phi_k^\ast)\|_{1,k}\}_{k\geq0}$ is uniformly bounded. Since  $\J_{k}^\eps(\beta) \leq c_2 + \gamma \mathcal{P}_\eps(\beta) = c_2 + \tfrac{\gamma}{\eps} (f(\beta),1)$, we deduce
    \[
         -c_1 + \gamma \tfrac{\eps}{2} |\phi_k^\ast|^2_{H^1(\Omega)} \leq  -c_1 + \gamma\mathcal{P}_\eps(\phi^\ast_k)  \leq \J_{k}^\eps(\phi_k^\ast) \leq \J_{k}^\eps(\beta) \leq c_2 +  \tfrac{\gamma}  {\eps} (f(\beta),1) , \quad \forall k\geq 0.
    \]
This and the box constraint in $\U_k$ imply that $\{\|\phi_k^\ast\|_{H^1(\Omega)}\}_{k\geq0}$ is uniformly bounded. Then $\U_\infty$ is a closed and convex subset of $\U$. So there exist a non-relabeled subsequence $\{\phi_{k_j}^\ast\}_{j\geq 0}$ and some $\overline{\phi}_\infty \in \U_\infty$ such that
    \begin{equation}\label{lem:sol-map_cont_limit_pf01}
        \phi_{k_j}^\ast \rightarrow \overline{\phi}_\infty \quad \text{weakly in}~H^1(\Omega),\quad \phi_{k_j}^\ast \to \overline{\phi}_\infty \quad \text{strongly in}~L^2(\Omega), \quad  \phi_{k_j}^\ast \to \overline{\phi}_\infty \quad \text{a.e. in}~\Omega.
    \end{equation}
By Lemmas \ref{lem:vel&costate_weak_adaptive} and \ref{lem:vel&costate_L2-strong-conv_adaptive}, we may assume that the subsequence $\{\bold{u}^\ast_{k_j}\}_{j\geq0}$ converges to some $\overline{\bold{u}}\in\bold{V}$
    \begin{equation}\label{lem:sol-map_cont_limit_pf02}
        \bold{u}^\ast_{k_j} \to \overline{\bold{u}} \quad \text{strongly in}~\bold{L}^2(\Omega),\quad \bold{\nabla}_{k_j}\bold{u}^\ast_{k_j} \rightarrow \bold{\nabla}\overline{\bold{u}}\quad\text{weakly in}~\mathbb{L}^2(\Omega).
    \end{equation}
The rest of the lengthy proof is divided into three steps: (1) $\overline{\bold{u}} = \S(\overline{\phi}_\infty)$; (2) strong convergence \eqref{vel_strong-conv_limit}; (3) $H^1(\Omega)$ strong convergence of $\{\phi_{k_j}^\ast\}_{j\geq 0}$.

\noindent\textit{Step 1.} $\overline{\bold{u}} = \S(\overline{\phi}_\infty)$. The residual $\mathcal{R}_2$ in \eqref{def:res2_G1} associated with $(\phi_{k_j}^\ast,\bold{u}_{k_j}^\ast)\in \U_{k_j} \times \bold{V}_{k_j}$ is given by
\[
        \langle \mathcal{R}_{2}(\phi_{k_j}^\ast,\bold{u}_{k_j}^\ast),\bold{v} \rangle : = \mu ( \bold{\nabla}_{k_j} \bold{u}_{k_j}^\ast, \bold{\nabla} \bold{v} ) + ( \alpha_\eps(\phi^\ast_{k_j}) \bold{u}_{k_j}^\ast, \bold{v}) - ( \bold{f}, \bold{v}) , \quad \forall \bold{v} \in \bold{V}.
\]
It follows from Lebesgue dominated convergence theorem, the assumption on $\alpha_\eps$ and pointwise convergence in  \eqref{lem:sol-map_cont_limit_pf01} that $\lim_{j\to \infty}\|\alpha_\eps(\phi_{k_j}^\ast)\bold{v} - \alpha_\eps(\overline{\phi}_{\infty})\bold{v}\|_{\bold{L}^2(\Omega)}=0$. This and \eqref{lem:sol-map_cont_limit_pf02} imply
\[
   \lim_{j\to\infty} \mu ( \bold{\nabla}_{k_j} \bold{u}_{k_j}^\ast, \bold{\nabla} \bold{v}) = \mu (\bold{\nabla} \overline{\bold{u}} , \bold{\nabla} \bold{v} ),\quad
   \lim_{j\to\infty} (\alpha_\eps(\phi^\ast_{k_j}) \bold{u}_{k_j}^\ast, \bold{v} ) = ( \alpha_\eps(\overline{\phi}_\infty) \overline{\bold{u}}, \bold{v} ), \quad \forall \bold{v}\in \bold{V}.
\]
Similar to the derivation of \eqref{estimate-pf01_state_G1}, we have
$$\left| \langle \mathcal{R}_{2}(\phi_k^\ast,\bold{u}_k^\ast),\bold{v} \rangle \right| \leq c \bigg(\sum_{T\in \cT_k}h_T^2\|\alpha_\eps(\phi_k^\ast) \bold{u}_k^\ast - \bold{f}\|_{\bold{L}^2(T)}^2\bigg)^{1/2} \|\bold{\nabla}\bold{v}\|_{\bold{L}^2(\Omega)} \leq c \eta_{k,2}(\phi_k^\ast,\bold{u}_k^\ast) \|\bold{\nabla}\bold{v}\|_{\bold{L}^2(\Omega)} ,\quad \forall \bold{v}\in\bold{V}.$$
This and Theorem \ref{thm:conv_est_state&costate} imply $\lim_{j\to\infty}\langle \mathcal{R}_{2}(\phi_{k_j}^\ast,\bold{u}_{k_j}^\ast),\bold{v} \rangle = 0$. Thus, $\overline{\bold{u}} \in \bold{V}$ solves problem \eqref{vp_stokes_div-free} with $\phi=\overline{\phi}_\infty$, i.e. $\overline{\bold{u}}=\S(\overline{\phi}_\infty)=\overline{\bold{u}}_\infty$.

\noindent\textit{Step 2.} Strong convergence  \eqref{vel_strong-conv_limit}. By taking $\bold{v}_{k_j} = \bold{u}_{k_j}^\ast$ over $\cT_{k_j}$ in \eqref{opt-sys_disc_state_G1}, passing to the limit $j\to\infty$ and using \eqref{lem:sol-map_cont_limit_pf02} again, we obtain
\[
   \lim_{j\to\infty} \mu \|\bold{\nabla}_{k_j}\bold{u}_{k_j}^\ast\|_{\bold{L}^2(\Omega)}^2 + ( \alpha_\eps(\phi_{k_j}^\ast),|\bold{u}_{k_j}^\ast|^2)  = \lim_{j\to\infty} ( \bold{f}, \bold{u}_{k_j}^\ast) = (\bold{f}, \overline{\bold{u}}_\infty) = \mu \|\bold{\nabla}\overline{\bold{u}}_{\infty}\|_{\bold{L}^2(\Omega)}^2 +( \alpha_\eps(\overline{\phi}_\infty),|\overline{\bold{u}}_{\infty}|^2).
\]
Lebesgue dominated convergence theorem, the assumption on $\alpha_\eps$, \eqref{lem:sol-map_cont_limit_pf01} and \eqref{lem:sol-map_cont_limit_pf02} further imply
\[
   \lim_{j\to\infty} \left|( \alpha_\eps(\phi_{k_j}^\ast),|\bold{u}_{k_j}^\ast|^2 - |\overline{\bold{u}}_\infty|^2) \right| \leq \lim_{j\to\infty} c \| \bold{u}_{k_j}^\ast - \overline{\bold{u}}_\infty \|_{\bold{L}^2(\Omega)} = 0\quad \mbox{and} \quad \lim_{j\to\infty}( \alpha_\eps(\phi_{k_j}^\ast), |\overline{\bold{u}}_\infty|^2) = ( \alpha_\eps(\overline{\phi}_{\infty}) , |\overline{\bold{u}}_\infty|^2),
\]
which {leads to}
\begin{equation}\label{lem:sol-map_cont_limit_pf03}
   \lim_{j\to\infty} ( \alpha_\eps(\phi_{k_j}^\ast), |\bold{u}_{k_j}^\ast|^2) = ( \alpha_\eps(\overline{\phi}_{\infty}), |\overline{\bold{u}}_\infty|^2 ).
\end{equation}
Hence $\lim_{j\to\infty}\|\bold{u}_{k_j}^\ast\|_{1,k_j}^2 = |\overline{\bold{u}}_{\infty}|_{H^1(\Omega)}^2$. This, \eqref{lem:sol-map_cont_limit_pf02} and $\overline{\bold{u}}=\S(\overline{\phi}_\infty)$ yield \eqref{vel_strong-conv_limit}.

\noindent\textit{Step 3.} $H^1(\Omega)$ {strong} convergence of $\{\phi_{k_j}\}_{j\geq0}$. The pointwise convergence in \eqref{lem:sol-map_cont_limit_pf01} and Lebesgue dominated convergence theorem imply
\begin{equation}\label{lem:sol-map_cont_limit_pf04}
  \lim_{j\to\infty}  (f(\phi_{k_j}^\ast),1) = ( f(\overline{\phi}_{\infty}),1).
\end{equation}
Then \eqref{vel_strong-conv_limit}, \eqref{lem:sol-map_cont_limit_pf03}, \eqref{lem:sol-map_cont_limit_pf04}, the weak convergence in \eqref{lem:sol-map_cont_limit_pf01} and the weak lower semi-continuity of the $H^1(\Omega)$-norm imply
\begin{equation}\label{lem:sol-map_cont_limit_pf05}
    \begin{aligned}
     &\tfrac{1}{2} ( \alpha_\eps(\overline{\phi}_\infty), |\overline{\bold{u}}_\infty|^2) + \tfrac{\mu}{2}\| \bold{\nabla} \overline{\bold{u}}_\infty \|_{\bold{L}^2(\Omega)}^2 - (\bold{f},\overline{\bold{u}}_\infty) +  \tfrac{\gamma\eps}{2}|\overline{\phi}_\infty|_{H^1(\Omega)}^2 + \tfrac{\gamma}{\eps}(f(\overline{\phi}_{\infty}),1) \\
     \leq &\liminf_{j\to\infty} \left(\tfrac{1}{2} (\alpha_\eps(\phi^\ast_{k_j}), |\bold{u}^\ast_{k_j}|^2) + \tfrac{\mu}{2} \| \bold{\nabla}_{k_j} \bold{u}^\ast_{k_j} \|_{\bold{L}^2(\Omega)}^2 - (\bold{f}, \bold{u}^\ast_{k_j})  + \tfrac{\gamma\eps}{2}|\phi_{k_j}^\ast|^2_{H^1(\Omega)} + \tfrac{\gamma}{\eps}( f(\phi_{k_j}^\ast),1) \right).
    \end{aligned}
\end{equation}
Now the density of $\bigcup_{k\geq0}\U_k$ in $\U_\infty$ in the $H^1(\Omega)$-topology admits the existence of a sequence $\{\overline{\phi}_k\}_{k\geq0} \subset \bigcup_{k\geq0}\U_k$ converging to $\overline{\phi}_\infty$ strongly in $H^1(\Omega)$. We may extract a non-relabeled subsequence $\{\overline{\phi}_{k_j}\}_{j\geq 0}$ with each $\overline{\phi}_{k_j} \in \U_{k_j}$ such that
\begin{equation}\label{lem:sol-map_cont_limit_pf06}
    \overline{\phi}_{k_j} \to \overline{\phi}_\infty \quad \text{strongly in}~H^1(\Omega),\quad \overline{\phi}_{k_j} \to \overline{\phi}_\infty \quad \text{a.e. in}~\Omega.
\end{equation}
Since $\bold{u}_{k_j}^\ast = \S_{k_j}(\phi_{k_j}^\ast)\in \bold{V}_{k_j}$ and $\S_{k_j}(\overline{\phi}_{k_j}) \in \bold{V}_{k_j}$ solve \eqref{disvp_stokes_div-free}  with $\phi_{k_j} = \phi_{k_j}^\ast$ and $\phi_{k_j} = \overline{\phi}_{k_j}$ over $\cT_{k_j}$, respectively, we have
\[
    \begin{aligned}
    \mu ( \bold{\nabla}_{k_j} ( \S_{k_j}(\phi_{k_j}^\ast) - \S_{k_j}(\overline{\phi}_{k_j}) ) , \bold{\nabla}_{k_j} \bold{v}_{k_j} ) +& ( (\alpha_\eps(\phi_{k_j}^\ast) - \alpha_\eps(\overline{\phi}_{k_j})) \S_{k_j}(\phi_{k_j}^\ast), \bold{v}_{k_j}) \\
    & + ( \alpha_\eps(\overline{\phi}_{k_j})( \S_{k_j}(\phi_{k_j}^\ast) - \S_{k_j}(\overline{\phi}_{k_j})) , \bold{v}_{k_j} )= 0,\quad \forall \bold{v}_{k_j}\in \bold{V}_{k_j}.
    \end{aligned}
\]
Setting $\bold{v}_{k_j} = \S_{k_j}(\phi_{k_j}^\ast) - \S_{k_j}(\overline{\phi}_{k_j})$ in the identity yields
\[
    \begin{aligned}
    &\quad \mu \| \S_{k_j}(\phi_{k_j}^\ast) - \S_{k_j}(\overline{\phi}_{k_j}) \|_{1,k_j}^2 + ( \alpha_\eps(\overline{\phi}_{k_j}),|\S_{k_j}(\phi_{k_j}^\ast) - \S_{k_j}(\overline{\phi}_{k_j})|^2) \\
    &  =( (\alpha_\eps(\overline{\phi}_{k_j}) - \alpha_\eps(\phi^\ast_{k_j})) \S_{k_j} (\phi_{k_j}^\ast) , \S_{k_j}(\phi_{k_j}^\ast) - \S_{k_j}(\overline{\phi}_{k_j}) ) \\
& = ( (\alpha_\eps(\overline{\phi}_{k_j}) - \alpha_\eps(\phi^\ast_{k_j})) (\S_{k_j} (\phi_{k_j}^\ast) - \S(\overline{\phi}_{\infty}) ), \S_{k_j}(\phi_{k_j}^\ast) - \S_{k_j}(\overline{\phi}_{k_j}))\\
&\quad + ((\alpha_\eps(\overline{\phi}_{k_j}) - \alpha_\eps(\phi^\ast_{k_j})) \S(\overline{\phi}_{\infty}) , \S_{k_j}(\phi_{k_j}^\ast) - \S_{k_j}(\overline{\phi}_{k_j}) ).
    \end{aligned}
\]
By the assumption on $\alpha_\eps$ and the discrete Poincar\'{e} inequality, we get
\[
    \| \S_{k_j}(\phi_{k_j}^\ast) - \S_{k_j}(\overline{\phi}_{k_j}) \|_{1,k_j} \leq c \big( \| \S_{k_j} (\phi_{k_j}^\ast) - \S(\overline{\phi}_{\infty}) \|_{\bold{L}^2(\Omega)} + \| (\alpha_\eps(\overline{\phi}_{k_j}) - \alpha_\eps(\phi^\ast_{k_j})) \S(\overline{\phi}_{\infty})  \|_{\bold{L}^2(\Omega)}  \big).
\]
By \eqref{vel_strong-conv_limit} and Lebesgue dominated convergence theorem, both terms vanish in the limit, which, along with the discrete Poincar\'{e} inequality again, implies
$$\lim_{j\to\infty}\| \S_{k_j}(\phi_{k_j}^\ast) - \S_{k_j}(\overline{\phi}_{k_j})\|_{\bold{L}^2(\Omega)} + \| \S_{k_j}(\phi_{k_j}^\ast) - \S_{k_j}(\overline{\phi}_{k_j}) \|_{1,k_j} = 0.$$
Hence by the triangle inequality and \eqref{vel_strong-conv_limit}, we arrive at
\begin{equation}\label{lem:sol-map_cont_limit_pf07}
 \lim_{j\to\infty}   \| \S_{k_j} (\overline{\phi}_{k_j}) - \S(\overline{\phi}_{\infty}) \|_{\bold{L}^2(\Omega)} +
    \| \S_{k_j} (\overline{\phi}_{k_j}) - \S(\overline{\phi}_{\infty}) \|_{1,k_j} = 0.
\end{equation}
Then the pointwise convergence in  \eqref{lem:sol-map_cont_limit_pf06}, the $L^2(\Omega)$ convergence in \eqref{lem:sol-map_cont_limit_pf07} and the argument for \eqref{lem:sol-map_cont_limit_pf03} imply
\begin{equation}\label{lem:sol-map_cont_limit_pf08}
  \lim_{j\to\infty}  (\alpha_\eps(\overline{\phi}_{k_j}), |\S_{k_j}(\overline{\phi}_{k_j})|^2) = (\alpha_\eps(\overline{\phi}_{\infty}),|\S(\overline{\phi}_{\infty})|^2 ).
\end{equation}
Now using the minimizing property of $(\phi_{k_j}^\ast,\bold{u}_{k_j}^\ast)$ over $\U_{k_j}$, \eqref{lem:sol-map_cont_limit_pf06}-\eqref{lem:sol-map_cont_limit_pf08}, the assumption on $\alpha_\eps$ and Lebesgue dominated convergence theorem, we derive
\begin{align}
    &\quad \limsup_{j\to\infty} \left( \tfrac{1}{2}  ( \alpha_\eps(\phi^\ast_{k_j}), |\bold{u}^\ast_{k_j}|^2) + \tfrac{\mu}{2} \| \bold{\nabla}_{k_j} \bold{u}^\ast_{k_j} \|_{\bold{L}^2(\Omega)}^2 - (\bold{f}, \bold{u}^\ast_{k_j})  + \tfrac{\gamma\eps}{2}|\phi_{k_j}^\ast|^2_{H^1(\Omega)} + \tfrac{\gamma}{\eps}( f(\phi_{k_j}^\ast),1) \right) \nonumber \\
    &\leq \limsup_{j\to\infty} \left( \tfrac{1}{2}( \alpha_\eps(\overline{\phi}_{k_j}), |\S_{k_j}(\overline{\phi}_{k_j})|^2) + \tfrac{\mu}{2} \| \bold{\nabla}_{k_j} \S_{k_j}(\overline{\phi}_{k_j}) \|_{\bold{L}^2(\Omega)}^2 - (\bold{f}, \S_{k_j}(\overline{\phi}_{k_j})) + \tfrac{\gamma\eps}{2}|\overline{\phi}_{k_j}|^2_{H^1(\Omega)} + \tfrac{\gamma}{\eps}(f(\overline{\phi}_{k_j}),1) \right) \nonumber \\
    & = \tfrac{1}{2}( \alpha_\eps(\overline{\phi}_{\infty}), |\S(\overline{\phi}_{\infty})|^2) + \tfrac{\mu}{2} \| \bold{\nabla} \S(\overline{\phi}_{\infty}) \|_{\bold{L}^2(\Omega)}^2 - (\bold{f}, \S(\overline{\phi}_{\infty}))  + \tfrac{\gamma\eps}{2}|\overline{\phi}_{\infty}|^2_{H^1(\Omega)} + \tfrac{\gamma}{\eps}(f(\overline{\phi}_{\infty}),1).   \label{lem:sol-map_cont_limit_pf09}
    \end{align}
In view of the relation $\overline{\bold{u}}_\infty = \S(\overline{\phi}_\infty)$, \eqref{lem:sol-map_cont_limit_pf05} and \eqref{lem:sol-map_cont_limit_pf09}, the following identity holds
\[
    \begin{aligned}
    &\quad
    \lim_{j\to\infty} \left( \tfrac{1}{2} ( \alpha_\eps(\phi^\ast_{k_j}), |\bold{u}^\ast_{k_j}|^2) + \tfrac{\mu}{2} \| \bold{\nabla}_{k_j} \bold{u}^\ast_{k_j} \|_{\bold{L}^2(\Omega)}^2 - (\bold{f}, \bold{u}^\ast_{k_j})  + \tfrac{\gamma\eps}{2}|\phi_{k_j}^\ast|^2_{H^1(\Omega)} + \tfrac{\gamma}{\eps}( f(\phi_{k_j}^\ast),1) \right) \\
    & = \tfrac{1}{2}( \alpha_\eps(\overline{\phi}_{\infty}), |\overline{\bold{u}}_{\infty}|^2) + \tfrac{\mu}{2} \| \bold{\nabla} \overline{\bold{u}}_{\infty} \|_{\bold{L}^2(\Omega)}^2 - (\bold{f}, \overline{\bold{u}}_{\infty})  + \tfrac{\gamma\eps}{2}|\overline{\phi}_{\infty}|^2_{H^1(\Omega)} + \tfrac{\gamma}{\eps}( f(\overline{\phi}_{\infty}),1).
    \end{aligned}
\]
This identity, together with \eqref{vel_strong-conv_limit}, \eqref{lem:sol-map_cont_limit_pf03} and \eqref{lem:sol-map_cont_limit_pf04} again, leads to $\lim_{j\to\infty}|\phi_{k_j}^\ast|_{H^1(\Omega)}^2 = |\phi_\infty|_{H^1(\Omega)}^2.$
This, the $H^1(\Omega)$ weak convergence and the $L^2(\Omega)$ convergence in \eqref{lem:sol-map_cont_limit_pf01} conclude the proof of the lemma.
\end{proof}

To deal with the subsequence $\{\eta_{k_j,1}(\phi_{k_j}^\ast,\bold{u}_{k_j}^\ast)\}_{j\geq 0}$ of error estimators, we need three preliminary results.    The next result is analogous to Lemma \ref{lem:max-est->zero}.
\begin{lemma}\label{lem:max-est_vi->zero}
Let $\{(\phi_{k_j}^\ast, \bold{u}_{k_j}^\ast)\}_{j\geq 0}$ be the convergent subsequence in Lemma \ref{lem:sol-map_cont_limit} and $\{\mathcal{M}_{k_j}\}_{j\geq 0}$ the corresponding subsequence of marked sets determined by \eqref{marking_G1}. Then there holds
    \begin{equation}\label{max-est_vi->zero}
        \lim_{j\to\infty} \max_{T\in\mathcal{M}_{k_j}} \eta_{k_j,1}(\phi_{k_j}^\ast, \bold{u}_{k_j}^\ast;T) = 0.
    \end{equation}
\end{lemma}

\begin{proof}
Let $T_{j}^1$ be the element attaining the largest error indicator among $\eta_{k_j,1}(T)$ over $\mathcal{M}_{k_j}$. Due to the marking condition \eqref{marking_G1}, $T_{j}^1\in \mathcal{M}_{k_j}^1 \subset \mathcal{M}_{k_j} \subset \cT_{k_j}^0$. Therefore, the local quasi-uniformity of $\cT_{k_j}$ and \eqref{mesh-size->zero} imply
\begin{equation}\label{lem:max-est_vi->zero_pf01}
  \lim_{j\to\infty}     h^d_{T_{j}^1} \leq \lim_{j\to\infty} |\omega_{k_j}(T_j^1)| \leq \lim_{j\to\infty} c \|h_{k_j}\chi_{k_j}^0\|_{L^\infty(\Omega)}^d= 0.
    \end{equation}
Since $\alpha_\eps\in C^1[0,1]$, $\alpha'_\eps(\phi_{k_j}^\ast)$ is bounded. Then by H\"{o}lder inequality and a scaled trace theorem, we have
    \begin{align*}
        h_{T_{j}^1}^2\|R_{T,1}^1(\phi_{k_j}^\ast,\bold{u}^\ast_{k_j},\bold{q}^\ast_{k_j})\|_{L^2(T^1_j)}^2&
        \leq c h_{T_{j}^1}^2 \big( \|f'(\phi_{k_j}^\ast)\|_{L^2(T_{j}^1)}^2 + \|\bold{u}_{k_j}^\ast\|_{\bold{L}^4(T_{j}^1)}^4 \big),\\
        \sum_{F\subset\partial T_j^1}h_F\|J_{F,1}^2(\phi_{k_j}^\ast)\|^2_{L^2(F)}
       & \leq c \|\bold{\nabla}\phi_{k_j}^\ast\|_{L^2\left(\omega_{k_j}(T_j^1)\right)}^2.
    \end{align*}
    Combining these two estimates further leads to
    \begin{equation}\label{lem:max-est_vi->zero_pf02}
       \eta_{k_j,1}^2(T_{j}^1) \leq c ( \|\bold{\nabla}(\phi_{k_j}^\ast - \overline{\phi}_\infty)\|^2_{L^2(\Omega)} + \|\bold{\nabla}\overline{\phi}_{\infty}\|^2_{L^2\left(\omega_{k_j}(T_j^1)\right)} + \|h_{k_j}\chi_{k_j}^0\|_{L^\infty(\Omega)}^2 ( \|f'(\phi_{k_j}^\ast)\|_{L^2(\Omega)}^2 + \|\bold{u}_{k_j}^\ast\|_{\bold{L}^4(\Omega)}^4)).
    \end{equation}
In view of the boundedness of $f'$ in $[0,1]$ and the pointwise convergence $\phi_{k_j}^\ast \to \overline{\phi}_{\infty}$ a.e. in $\Omega$, cf. \eqref{lem:sol-map_cont_limit_pf01}, $$\lim_{j\to\infty}\|f'(\phi_{k_j}^\ast) - f'(\overline{\phi}_{\infty})\|_{L^2(\Omega)}^2 = 0,$$ by Lebesgue dominated convergence theorem. So $\{\|f'(\phi_{k_j}^\ast)\|_{L^2(\Omega)}\}_{j\geq 0}$ is uniformly bounded. Then by Lemmas \ref{lem:sol-map_cont_limit}, \begin{equation*}
\lim_{j\to\infty}\|\bold{u}_{k_j}^\ast - \overline{\bold{u}}_\infty\|_{1,k_j}=0.
\end{equation*}
This and the discrete Sobolev inequality \eqref{disc_Sob_ineq_hom} imply that $\{\|\bold{u}_{k_j}^\ast\|_{\bold{L}^4(\Omega)}\}_{j\geq 0}$ is uniformly bounded. Then we may deduce from \eqref{mesh-size->zero} that the third term in \eqref{lem:max-est_vi->zero_pf02} vanishes in the limit. By Lemma \ref{lem:sol-map_cont_limit}, $\lim_{j\to\infty}\|\bold{\nabla}(\phi_{k_j}^\ast - \overline{\phi}_\infty)\|^2_{L^2(\Omega)}= 0$. The second term also tends to zero by \eqref{lem:max-est_vi->zero_pf01} and the absolute continuity of the norm $\|\cdot\|_{L^2(\Omega)}$ with respect to the Lebesgue measure.
\end{proof}

The stability in non-refined elements and the reduction on refined elements given in the next two lemmas form the basis to prove Theorem \ref{thm:conv_est_control}. These two properties of an a posteriori error estimator first appeared in \cite{CKNS:2008} and were developed in \cite{GantnerPraetorius:2022} for the plain convergence of the estimator to zero in adaptive conforming methods for direct problems. We adapt the properties for a nonconforming FEM that leads to $\eta_{k,1}(\phi_k^\ast, \bold{u}_{k}^\ast)$.
\begin{lemma}\label{lem:est_perturb_nonrefined}
    Let $\{(\phi_{k_j}^\ast, \bold{u}_{k_j}^\ast)\}_{j\geq 0}$ be the convergent subsequence in Lemma \ref{lem:sol-map_cont_limit}. Then for $j<l$, there holds
\begin{equation}\label{est_perturb_nonrefined}
        \begin{aligned}
        \eta_{k_l,1}(\cT_{k_j}\cap \cT_{k_l}) \leq
        \eta_{k_j,1}(\cT_{k_j}\cap \cT_{k_l})&  + c\left(\|\bold{\nabla}(\phi_{k_j}^\ast - \phi_{k_l}^\ast)\|_{L^2(\Omega)} + \|f'(\phi_{k_j}^\ast) - f'(\phi_{k_l}^\ast)\|_{L^2(\Omega)} \right. \\
        & \qquad \left. + \|\bold{u}_{k_j}^\ast - \bold{u}_{k_l}^\ast\|_{\bold{L}^2(\Omega)} + \| \alpha'_\eps(\phi^\ast_{k_j}) - \alpha'_\eps(\phi^\ast_{k_l}) \|_{L^6(\Omega)}\right).
        \end{aligned}
    \end{equation}
\end{lemma}

\begin{proof}
    By the triangle inequality, we have
    \begin{equation}\label{lem:est_perturb_nonrefined_pf01}
               \eta_{k_l,1}(\phi_{k_l}^\ast,\bold{u}_{k_l}^\ast;\cT_{k_j}\cap \cT_{k_l}) \leq \eta_{k_j,1}(\phi_{k_j}^\ast,\bold{u}_{k_j}^\ast;\cT_{k_j}\cap \cT_{k_l}) + \bigg(\sum_{T\in \cT_{k_j}\cap\cT_{k_l}} \left( \mathrm{I} + \mathrm{II} \right) \bigg)^{1/2}
    \end{equation}
     with the two estimators $\mathrm{I}$ and $\mathrm{II}$ given respectively by
     \begin{align*}
        \mathrm{I} &: = h_T^2\left\| \tfrac{\gamma}{\eps} f'(\phi^\ast_{k_j}) + \tfrac{1}{2}\alpha'_\eps(\phi^\ast_{k_j}) |\bold{u}_{k_j}^\ast|^2  - \tfrac{\gamma}{\eps} f'(\phi^\ast_{k_l}) - \tfrac{1}{2}\alpha'_\eps(\phi^\ast_{k_l}) |\bold{u}_{k_l}^\ast|^2   \right\|_{L^2(T)}^2,\\
        \mathrm{II} &: = \sum_{F \subset \partial T}h_F\|[\bold{\nabla}\phi_{k_j}^\ast - \bold{\nabla}\phi_{k_l}^\ast] \cdot \bold{n}_{F}\|_{L^2(F)}^2.
     \end{align*}
By the scaled trace theorem and the finite overlapping property of $\omega_{k_l}(T)$, we have
     \begin{equation}\label{lem:est_perturb_nonrefined_pf02}
        \sum_{T\in \cT_{k_j}\cap\cT_{k_l}}\mathrm{II} \leq c \sum_{T\in \cT_{k_j}\cap\cT_{k_l}}\|\bold{\nabla}(\phi_{k_j}^\ast - \phi_{k_l}^\ast)\|_{L^2(\omega_{k_l}(T))}^2  \leq c \| \bold{\nabla}(\phi_{k_j}^\ast - \phi_{k_l}^\ast)\|^2_{L^2(\Omega)}.
     \end{equation}
By the boundedness of $\alpha'_\eps$ in $[0,1]$, the inverse estimate and H\"{o}lder inequality, we obtain on any $T\in \cT_{k_j}\cap\cT_{k_l}$,
\begin{align*}
    h_T^2\|\alpha'_\eps(\phi^\ast_{k_j}) (|\bold{u}_{k_j}^\ast|^2 -  |\bold{u}_{k_l}^\ast|^2)\|_{L^2(T)}^2
    & \leq c h_T^2\|\bold{u}_{k_j}^\ast - \bold{u}_{k_l}^\ast\|_{\bold{L}^4(T)}^2 \|\bold{u}_{k_j}^\ast + \bold{u}_{k_l}^\ast\|_{\bold{L}^4(T)}^2  \\
    &\leq c h_T^{2-\frac{d}{2}}\| \bold{u}_{k_j}^\ast - \bold{u}_{k_l}^\ast \|^2_{\bold{L}^2(T)} \|\bold{u}_{k_j}^\ast + \bold{u}_{k_l}^\ast\|_{\bold{L}^4(\Omega)}^2,\\
        h_T^2 \| (\alpha'_\eps(\phi^\ast_{k_j}) - \alpha'_\eps(\phi^\ast_{k_l})) |\bold{u}_{k_l}^\ast|^2 \|_{L^2(T)}^2 & \leq h_T^2 \| \alpha'_\eps(\phi^\ast_{k_j}) - \alpha'_\eps(\phi^\ast_{k_l}) \|_{L^6(T)}^2 \| \bold{u}_{k_l}^\ast \|_{\bold{L}^6(T)}^4  \\
        & \leq  c h_T^{2-\frac{d}{3}} \| \alpha'_\eps(\phi^\ast_{k_j}) - \alpha'_\eps(\phi^\ast_{k_l}) \|_{L^6(\Omega)}^2 \| \bold{u}_{k_l}^\ast \|_{\bold{L}^4(T)}^4.
    \end{align*}
Summing up these two estimates over all $T\in \cT_{k_j}\cap\cT_{k_l}$ and the Cauchy-Schwarz inequality yield
\begin{align}
   &\quad \sum_{T\in\cT_{k_j}\cap\cT_{k_l}} h_T^2 \| \alpha'_\eps(\phi^\ast_{k_j}) |\bold{u}_{k_j}^\ast|^2 - \alpha'_\eps(\phi^\ast_{k_l}) |\bold{u}_{k_l}^\ast|^2\|_{L^2(T)}^2 \nonumber \\
            & \leq c \big(\| \bold{u}_{k_j}^\ast - \bold{u}_{k_l}^\ast \|^2_{\bold{L}^2(\Omega)} \| \bold{u}_{k_j}^\ast + \bold{u}_{k_l}^\ast\|_{\bold{L}^4(\Omega)}^2 + \| \alpha'_\eps(\phi^\ast_{k_j}) - \alpha'_\eps(\phi^\ast_{k_l}) \|_{L^6(\Omega)}^2 \| \bold{u}_{k_l}^\ast \|_{\bold{L}^4(\Omega)}^4  \big) \nonumber \\
            &\leq c \big(\| \bold{u}_{k_j}^\ast - \bold{u}_{k_l}^\ast \|^2_{\bold{L}^2(\Omega)} + \| \alpha'_\eps(\phi^\ast_{k_j}) - \alpha'_\eps(\phi^\ast_{k_l}) \|_{L^6(\Omega)}^2 \big) , \label{lem:est_perturb_nonrefined_pf04}
    \end{align}
since the sequence $\{\|\bold{u}_{k_j}^\ast\|_{\bold{L}^4(\Omega)}\}_{j\geq0}$ is uniformly bounded, in view of the discrete Sobolev inequality \eqref{disc_Sob_ineq_hom} and the strong convergence in Lemma \ref{lem:sol-map_cont_limit}. The preceding estimates complete the proof of the lemma.
\end{proof}

\begin{lemma}\label{lem:est_reduction_refined}
    Let $\{(\phi_{k_j}^\ast, \bold{u}_{k_j}^\ast)\}_{j\geq 0}$ be the convergent subsequence in Lemma \ref{lem:sol-map_cont_limit}. Then for $j<l$, there exists a $q\in (0,1)$ such that
    \begin{equation}\label{est_reduction_refined}
        \eta_{k_l,1}^2(\cT_{k_l}\setminus\cT_{k_j}) \leq q \eta_{k_j ,1}^2(\cT_{k_j}\setminus\cT_{k_l})
          + c\epsilon_1(k_j,k_l)
    \end{equation}
    with $\epsilon_1(k_j,k_l)= \|\bold{\nabla}(\phi_{k_j}^\ast - \phi_{k_l}^\ast)\|^2_{L^2(\Omega)} + \|f'(\phi_{k_j}^\ast) - f'(\phi_{k_l}^\ast)\|^2_{L^2(\Omega)}  + \|\bold{u}_{k_j}^\ast - \bold{u}_{k_l}^\ast\|^2_{\bold{L}^2(\Omega)} + \| \alpha'_\eps(\phi^\ast_{k_j}) - \alpha'_\eps(\phi^\ast_{k_l}) \|_{L^6(\Omega)}^2$.
\end{lemma}

\begin{proof}
Similar to Lemma \ref{lem:est_perturb_nonrefined}, we use the triangle inequality for each $T\in \cT_{k_l}\setminus \cT_{k_j}$
    \begin{equation}\label{lem:est_reduction_refined_pf01}
        \eta_{k_l, 1}(\phi_{k_l}^\ast,\bold{u}_{k_l}^\ast;T) \leq \eta_{k_l,1}(\phi_{k_j}^\ast,\bold{u}_{k_j}^\ast;T) + (\mathrm{I} + \mathrm{II})^{1/2},
    \end{equation}
    with $\mathrm{I}$ and $\mathrm{II}$ given in \eqref{lem:est_perturb_nonrefined_pf01}. Since $\bold{u}_{k_j}^\ast$ is linear in each $T\in \cT_{k_l}\setminus \cT_{k_j}$,  the argument for \eqref{lem:est_perturb_nonrefined_pf02}-\eqref{lem:est_perturb_nonrefined_pf04} yields
$\sum_{T\in \cT_{k_l}\setminus \cT_{k_j}}(\mathrm{I} + \mathrm{II}) \leq c \epsilon_1(k_j,k_l)$.
By applying Young's inequality with $\delta > 0$ to \eqref{lem:est_reduction_refined_pf01} and summing the resulting estimate over $T\in \cT_{k_l}\setminus \cT_{k_j}$, we obtain
        \begin{align*}
        \eta_{k_l,1}^2(\phi_{k_l}^\ast,\bold{u}_{k_l}^\ast;&\cT_{k_l}\setminus\cT_{k_j}) \leq
        (1+\delta)\eta_{k_l,1}^2(\phi_{k_j}^\ast,\bold{u}_{k_j}^\ast;\cT_{k_l}\setminus\cT_{k_j}) + (1+\delta^{-1})c\epsilon_1(k_j,k_l).
        \end{align*}
Since $\cT_{k_l}$ is a refinement of $\cT_{k_j}$ by bisection, each $T\in \cT_{k_l} \setminus \cT_{k_j}$ is a successor of an element $T'\in \cT_{k_j}\setminus \cT_{k_l}$. Thus $J_{F,1}^1(\phi_{k_j}^\ast)$ is zero across $F\in \mathcal{F}_{k_l}$ in the interior of $T'$ since $\phi_{k_j}^\ast \in \U_{k_j}\subset S_{k_j}$ and $\bold{\nabla}\phi_{k_j}^\ast$ is a constant vector in $T'$. For any $F=\partial T \cap F'$ with $F'\subset \partial T'$ and $F'\in \mathcal{F}_{k_j}$, by the bisection rule, we have
$$h_T^d = |T| \leq \tfrac{1}{2} |T'| = \tfrac{1}{2} h_{T'}^d\quad\mbox{and}\quad h_{F}^{d-1}  = |F| \leq \tfrac{1}{2} |F'| = \tfrac{1}{2} h_{F'}^{d-1}.$$
Then we get
    \[
        \eta^2_{k_l,1}(\phi_{k_j}^\ast,\bold{u}_{k_j}^\ast;\cT_{k_l}\setminus\cT_{k_j}) \leq  2^{-1/(d-1)}\eta^2_{k_j,1}(\phi_{k_j}^\ast,\bold{u}_{k_j}^\ast;\cT_{k_j}\setminus\cT_{k_l}),
    \]
which completes the proof with a sufficiently small $\delta >0$ such that $q=(1+\delta)2^{-1/(d-1)}\in (0,1)$.
\end{proof}

The next result gives the vanishing limit of the estimator subsequence $\{\eta_{k_j,1}\}_{j\geq 0}$.
\begin{theorem}\label{thm:conv_est_control}
    The subsequence $\{\eta_{k_j,1}\}_{j\geq 0}$ associated with the convergent subsequence $\{(\phi_{k_j}^\ast, \bold{u}_{k_j}^\ast)\}_{j\geq 0}$ in Lemma \ref{lem:sol-map_cont_limit} converges to zero.
\end{theorem}
\begin{proof}
The proof is inspired by the work \cite{GantnerPraetorius:2022}, and uses also the ideas in \cite{MorinSiebert:2008} and \cite{Siebert:2011}. The main change consists in an enlarged set for $j<l$,
$$\cT_{k_{j}\to k_{l}}^+:=\{ T \in \cT_{k_j} \cap \cT_{k_l}~|~ T \cap \Omega^+_{k_j} \neq \emptyset \}$$
and the geometric observation \eqref{thm:conv_est_control_pf04}. By definition, $\cT_{k_j}^+ = \bigcap_{m\geq k_j} \cT_{m}^+ \subset \cT_{k_j} \cap \cT_{k_l}$. Hence, the set $\cT_{k_{j}\to k_{l}}^+$ is well-defined and consists of all elements in $\cT_{k_j}$ neighboring $\Omega^+_{k_j}$ and not refined until at least the $k_l$-th adaptive loop. Moreover, the uniform shape regularity of $\{\cT_k\}_{k\geq0}$ by Algorithm \ref{alg_anfem_topopt-Stokes_phase-field_G1} entails
    \begin{equation}\label{thm:conv_est_control_pf01}
        \# \cT_{k_{j}\to k_{l}}^+ \leq c \# \cT_{k_{j}}^+,
    \end{equation}
    with the constant $c$ depending only on $\cT_0$. Now we decompose $\cT_{k_l}$ into \cite{GantnerPraetorius:2022} \cite{LiXuZhu:2023}
        \begin{align*}
        \cT_{k_l}&=\big(\cT_{k_l}\setminus\cT_{k_j}\big)\cup \big(\cT_{k_l}\cap\cT_{k_j}\big)
				= \big(\cT_{k_l}\setminus\cT_{k_j}\big)\cup \big(\cT_{k_l}\cap((\cT_{k_j}\setminus\cT^+_{k_j\to k_l})\cup \cT^+_{k_j \to k_l})\big)\\
				& = \big(\cT_{k_l}\setminus\cT_{k_j}\big) \cup \big(\cT_{k_l}\cap (\cT_{k_j}\setminus\cT^+_{k_j\to k_l}) \big)  \cup \big(\cT_{k_l}\cap\cT^+_{k_j\to k_l}\big) \\
       &=
\big(\cT_{k_l}\setminus\cT_{k_j}\big)\cup \big(\cT_{k_l}\cap (\cT_{k_j}\setminus\cT^+_{k_j\to k_l})\big)  \cup \cT^+_{k_j\to k_l}.
\end{align*}
The set $\cT_{k_l}\setminus\cT_{k_j}$ contains elements in $\cT_{k_l}$ generated by refinements after the ${k_j}$-th adaptive loop, and elements in $\cT_{k_l}\cap (\cT_{k_j}\setminus\cT^+_{k_j\to k_l})$ will eventually be refined after the $k_l$-th adaptive loop. These properties are important to deriving \eqref{thm:conv_est_control_pf04} in Step 2 below.   
For $j<l$, $\eta_{k_l,1}^2(\phi_{k_l}^\ast, \bold{u}_{k_l}^\ast)$ is split into
\begin{equation}\label{thm:conv_est_control_pf02}
    \eta_{k_l,1}^2(\phi_{k_l}^\ast, \bold{u}_{k_l}^\ast) =  \eta_{k_l,1}^2( \cT_{k_l}\setminus\cT_{k_j}) +  \eta_{k_l,1}^2(\cT_{k_l}\cap (\cT_{k_j}\setminus\cT^+_{k_j\to k_l})) + \eta_{k_l,1}^2( \cT^+_{k_j\to k_l}).
\end{equation}
The rest of the proof is to bound the three terms in \eqref{thm:conv_est_control_pf02} separately, and proceeds in three steps.

\noindent\textit{Step 1.} By Lemma \ref{lem:est_reduction_refined}, we have
\begin{equation}\label{thm:conv_est_control_pf03}
    \eta_{k_l,1}^2( \cT_{k_l}\setminus\cT_{k_j}) \leq q \eta_{k_j,1}^2 ( \cT_{k_j}\setminus\cT_{k_l}) + \epsilon_{1}(k_j,k_l).
\end{equation}
with $\epsilon_{1}(k_j,k_l)\to 0$ as $j,l\to\infty$, which is guaranteed by Lemma \ref{lem:sol-map_cont_limit} for $\|\bold{\nabla}(\phi_{k_j}^\ast - \phi_{k_l}^\ast)\|^2_{L^2(\Omega)} +\|\bold{u}_{k_j}^\ast - \bold{u}_{k_l}^\ast\|^2_{\bold{L}^2(\Omega)}$ and the pointwise convergence in \eqref{lem:sol-map_cont_limit_pf01}, the continuity of $f'$, $\alpha_{\eps}'$ in $[0,1]$ and Lebesgue dominated convergence theorem for $\|f'(\phi_{k_j}^\ast) - f'(\phi_{k_l}^\ast)\|^2_{L^2(\Omega)} + \| \alpha'_\eps(\phi^\ast_{k_j}) - \alpha'_\eps(\phi^\ast_{k_l}) \|_{L^6(\Omega)}^2$ respectively.

\noindent\textit{Step 2.} From $\cT_{k_j}^+\subset \cT_{k_j\to k_l}^+$ for $j<l$, we deduce that $\cT_{k_j}\setminus \cT_{k_j\to k_l}^+
\subset \cT_{k_j}\setminus \cT_{k_j}^+ = \cT_{k_j}^0$ and $\cT_{k_l}\cap (\cT_{k_j}\setminus\cT^+_{k_j \to k_l})\subset \cT_{k_l}\cap\cT_{k_j}^0$.
Since any element $T\in \cT_{k_l}\cap\cT_{k_j}^0$ is not refined until after the $k_l$-th iteration, we have $\cT_{k_l}\cap\cT_{k_j}^0=\cT_{k_l}^0\cap\cT_{k_j}^0 \supset \cT_{k_l}\cap (\cT_{k_j}\setminus\cT^+_{k_j \to k_l})$. When $l$ is sufficiently large, due to the property \eqref{mesh-size->zero}, the mesh size of each element in $\cT_{k_l}^0$ is smaller than that of each one in $\cT_{k_j}^0$ for each fixed $j\in\mathbb{N}_0$. Thus, $\cT_{k_j}^0\cap\cT_{k_l}^0=\emptyset$ for sufficiently large $l$ and by letting $\Omega(\cT^0_{k_j}\cap\cT^0_{k_l}):=\bigcup_{T\in\cT^0_{k_j}\cap\cT^0_{k_l}}T$, there holds that for each fixed $j\in\mathbb{N}_0$,
\begin{equation}\label{thm:conv_est_control_pf04}
\lim_{l\to\infty}\big|\Omega(\cT^0_{k_j}\cap\cT^0_{k_l})\big|=0.
\end{equation}
The following local stability for $\eta_{k_l,1}(\phi_{k_l}^\ast,\bold{u}_{k_l}^\ast;T)$  on $T\in \cT_{k_l}$ in Lemma \ref{lem:max-est_vi->zero} and the inverse estimate give
    \begin{align*}
    \eta^2_{k_l,1}(\phi_{k_l}^\ast,\bold{u}_{k_l}^\ast;T) & \leq c ( \|\bold{\nabla}\phi_{k_l}^\ast\|_{L^2(\omega_{k_l}(T))}^2+ \|f'(\phi_{k_l}^\ast)\|_{L^2(T)}^2 + h_T^2\|\bold{u}_{k_l}^\ast\|_{\bold{L}^4(T)}^4   )  \\
    & \leq c ( \|\bold{\nabla}\phi_{k_l}^\ast\|_{L^2(\omega_{k_l}(T))}^2+ \|f'(\phi_{k_l}^\ast)\|_{L^2(T)}^2 + h_T^{2-\frac{d}{2}}\|\bold{u}_{k_l}^\ast\|_{\bold{L}^2(T)}^2\|\bold{u}_{k_l}^\ast\|_{\bold{L}^4(\Omega)}^2 ).
    \end{align*}
By summing the inequality over $T \in \cT_{k_l}\cap (\cT_{k_j}\setminus\cT^+_{k_j\to k_l})$, noting the relation $\cT_{k_l}\cap (\cT_{k_j}\setminus\cT^+_{k_j\to k_l}) \subset \cT_{k_l}^0\cap\cT_{k_j}^0 $ and using the Cauchy-Schwarz inequality and the discrete Sobolev inequality \eqref{disc_Sob_ineq_hom}, we arrive at
    \begin{align*}
        &\eta_{k_l,1}^2(\cT_{k_l}\cap (\cT_{k_j}\setminus\cT^+_{k_j\to k_l})) \\
        \leq& c \sum_{T\in \cT_{k_l}\cap (\cT_{k_j}\setminus\cT^+_{k_j\to k_l})} \big( \|\bold{\nabla}\phi_{k_l}^\ast\|_{L^2(\omega_{k_l}(T))}^2+ \|f'(\phi_{k_l}^\ast)\|_{L^2(T)}^2 + \|\bold{u}_{k_l}^\ast\|_{\bold{L}^2(T)}^2\|\bold{u}_{k_l}^\ast\|_{\bold{L}^4(\Omega)}^2  \big) \\
        \leq& c \big( \|\bold{\nabla}\phi_{k_l}^\ast\|_{L^2(\Omega(\cT^0_{k_j}\cap\cT^0_{k_l}))}^2+ \|f'(\phi_{k_l}^\ast)\|_{L^2(\Omega(\cT^0_{k_j}\cap\cT^0_{k_l}))}^2 + \|\bold{u}_{k_l}^\ast\|_{\bold{L}^2(\Omega(\cT^0_{k_j}\cap\cT^0_{k_l}))}^2 \|\bold{u}_{k_l}^\ast\|_{\bold{L}^4(\Omega)}^2  \big)\\
        \leq &c \Big( \|\bold{\nabla}(\phi_{k_l}^\ast-\overline{\phi}_{\infty})\|_{L^2(\Omega)}^2+ \|\bold{\nabla}\overline{\phi}_{\infty}\|_{L^2(\Omega(\cT^0_{k_j}\cap\cT^0_{k_l}))}^2+ \|f'(\phi_{k_l}^\ast) - f'(\overline{\phi}_\infty)\|_{L^2(\Omega)}^2 + \|f'(\overline{\phi}_{\infty})\|_{L^2(\Omega(\cT^0_{k_j}\cap\cT^0_{k_l}))}^2 \\
        & \qquad  + (\|\bold{u}_{k_l}^\ast - \overline{\bold{u}}_{\infty}\|_{\bold{L}^2(\Omega)}^2 + \|\overline{\bold{u}}_{\infty}\|_{\bold{L}^2(\Omega(\cT^0_{k_j}\cap\cT^0_{k_l}))}^2) \|\bold{u}_{k_l}^\ast\|_{1,k_l}^2    \Big).
    \end{align*}
Similar to Step 1, $\|\bold{\nabla}(\phi_{k_l}^\ast-\overline{\phi}_{\infty})\|_{L^2(\Omega)}^2$, $\|f'(\phi_{k_l}^\ast) - f'(\overline{\phi}_\infty)\|_{L^2(\Omega)}^2$ and $\|\bold{u}_{k_l}^\ast - \overline{\bold{u}}_{\infty}\|_{\bold{L}^2(\Omega)}^2$ all tend to zero as $l\to\infty$. By \eqref{thm:conv_est_control_pf04} and the absolute continuity of the $L^2$ norm with respect to the Lebesgue measure, $\|\bold{\nabla}\overline{\phi}_{\infty}\|_{L^2(\Omega(\cT^0_{k_j}\cap\cT^0_{k_l}))}^2$, $\|f'(\phi_{k_l}^\ast)\|_{L^2(\Omega(\cT^0_{k_j}\cap\cT^0_{k_l}))}^2$ and $\|\overline{\bold{u}}_{\infty}\|_{\bold{L}^2(\Omega(\cT^0_{k_j}\cap\cT^0_{k_l}))}^2$ also tend to zero as $l\to\infty$, for each fixed $j\in \mathbb{N}_0$. These vanishing limits and the uniform boundedness of the sequence $\{\|\bold{u}_{k_j}^\ast\|_{1,k_j}\}_{j\geq0}$ (cf. Lemma \ref{lem:sol-map_cont_limit}) imply that for each fixed $j\in \mathbb{N}_0$,
\begin{equation}\label{thm:conv_est_control_pf05}
    \lim_{l\to\infty} \eta_{k_l,1}^2(\cT_{k_l}\cap (\cT_{k_j}\setminus\cT^+_{k_j\to k_l})) = 0.
\end{equation}

\noindent \textit{Step 3.} For the term $\eta_{k_l,1}^2(\cT^+_{k_j\to k_l})$, we resort to \eqref{marking_G1} in the module MARK of Algorithm \ref{alg_anfem_topopt-Stokes_phase-field_G1}. In fact, it follows from \eqref{thm:conv_est_control_pf01}, \eqref{marking_G1} and Lemma \ref{lem:max-est_vi->zero} that for each fixed $j\in\mathbb{N}_0$,
\begin{equation}\label{thm:conv_est_control_pf06}
    \eta_{k_l,1}^2(\cT^+_{k_j\to k_l}) \leq \#\cT^+_{k_j\to k_l} \max_{T\in \cT^+_{k_j\to k_l}} \eta_{k_l,1}(T) \leq c \#\cT_{k_j}^+ \max_{T\in \mathcal{M}_{k_l}}\eta_{k_l,1}(T) \to 0 \quad \text{as}~l\to\infty.
\end{equation}
Now combining \eqref{thm:conv_est_control_pf05} and \eqref{thm:conv_est_control_pf06} imply that for each fixed $j\in\mathbb{N}_0$, there exists $N(j)>j$ such that for each $l>N(j)$,
\begin{equation}\label{thm:conv_est_control_pf07}
    \eta^2_{k_l,1}(\cT_{k_l}\cap\cT_{k_j}) \leq (q'-q)\eta^2_{k_j,1}
\end{equation}
with some $q'\in (q,1)$ and $q$ as given in \eqref{thm:conv_est_control_pf03}. Then \eqref{thm:conv_est_control_pf03} and \eqref{thm:conv_est_control_pf07} allow extracting a further non-relabeled subsequence $\{\eta_{k_{j_n},1}\}_{n\geq0}$ satisfying $\eta^2_{k_{j_{n+1}},1} \leq q' \eta^2_{k_{j_{n}},1} + \epsilon_{1}(k_{j_n},k_{j_{n+1}})$ with $\epsilon_{1}(k_{j_n},k_{j_{n+1}})\to 0$ as $n\to\infty$ as asserted in Step 1. Hence, an elementary argument in calculus (cf. Lemma 3.2 in \cite{GantnerPraetorius:2022}) yields $\eta_{k_{j_n},1}\to 0$ as $n\to\infty$. Finally, by Lemmas \ref{lem:est_perturb_nonrefined}-\ref{lem:est_reduction_refined} and the argument in Step 1 again, $\eta_{k_j,1}^2 \leq 2 \eta^2_{k_{j_n},1} + \epsilon_{1}(k_{j_n},k_{j})$ with $\epsilon_{1}(k_{j_n},k_{j}) \to 0$ as $k_j > k_{j_n}\to\infty$ and $n\to\infty$. The desired conclusion follows from the last two results.
\end{proof}

\begin{remark}\label{rem:conv_est_control}
Using an alternative decomposition of the error estimator, Morin et al \cite{MorinSiebert:2008} prove the plain convergence of AFEM for a class of linear elliptic problems. The three parts in the decomposition are associated with three classes of elements that are always refined, never refined and between the two extreme cases, respectively, and they are shown to tend to zero separately as the adaptive algorithm proceeds. The contractive property up to a Cauchy sequence in \eqref{thm:conv_est_control_pf03} indicates that the decomposition \eqref{thm:conv_est_control_pf02} is completely different from that in \cite{MorinSiebert:2008}. However, the proof of \eqref{thm:conv_est_control_pf05} in Step 2, especially the geometric observation \eqref{thm:conv_est_control_pf04}, is somewhat motivated by \cite[Proposition 4.2]{MorinSiebert:2008} for the elements in the intermediate state while we follow the argument in \cite{Siebert:2011} for \eqref{thm:conv_est_control_pf06}. Finally, the subsequence argument at the end of the proof is borrowed from \cite{GantnerPraetorius:2022}. This analysis strategy has been applied to the plain convergence of an adaptive conforming FEM for an eigenvalue optimization problem \cite{LiXuZhu:2023} and an inverse problem from cardiac electrophysiology \cite{JinWangXu:2025}. In this work, we generalize the strategy to a conforming-nonconforming coupled method.
\end{remark}

\subsection{Convergence of the discrete minimizing pair}\label{subsect:conv_phase&velocity}

Now we can state the main result of the work.
\begin{theorem}\label{thm:conv_adaptive}
The sequence of discrete minimizing pairs $\{(\phi_{k}^\ast,\bold{u}_k^\ast)\}_{k\geq0}$ generated by Algorithm \ref{alg_anfem_topopt-Stokes_phase-field_G1} contains a subsequence $\{ (\phi_{k_j}^\ast,\bold{u}_{k_j}^\ast) \}_{j\geq 0}$ converging to a solution $(\phi^\ast,\bold{u}^\ast)\in \U\times \bold{V}$ of problem \eqref{opt-sys_G1} in the sense that
\begin{equation}\label{conv_adaptive}
   \lim_{j\to\infty} \|\phi_{k_j}^\ast - \phi^\ast\|_{H^1(\Omega)}+ \|\bold{u}^\ast_{k_j} - \bold{u}^\ast \|_{\bold{L}^2(\Omega)} + \|\bold{\nabla}_{k_j}(\bold{u}^\ast_{k_j} - \bold{u}^\ast) \|_{\bold{L}^2(\Omega)} = 0.
\end{equation}
\end{theorem}
\begin{proof}
     By Lemma \ref{lem:sol-map_cont_limit}, we have a convergent subsequence $\{(\phi_{k_j}^\ast,\bold{u}_{k_j}^\ast)\}_{j\geq0}$ with the limit $(\overline{\phi}_\infty, \overline{\bold{u}}_\infty)$ solving \eqref{opt-sys_state_G1}. It suffices to prove that \eqref{opt-sys_control_G1} also holds with $(\phi^\ast,\bold{u}^\ast)=(\overline{\phi}_\infty, \overline{\bold{u}}_\infty)$. We invoke the residual operator $\mathcal{R}_{1}(\phi_{k}^\ast,\bold{u}_k^\ast)$ defined in \eqref{def:res1_G1} and argue as in \eqref{estimate-pf01_control_G1}, \eqref{estimate-pf02_control_G1} and \eqref{estimate-pf03_control_G1} to get
     \begin{align}\label{thm:conv_adaptive_pf01}
        \langle \mathcal{R}_{1}(\phi_{k_j}^\ast,\bold{u}_{k_j}^\ast), \phi - \phi_{k_j}^\ast \rangle &\geq \langle \mathcal{R}_{1}(\phi_{k_j}^\ast,\bold{u}_{k_j}^\ast), \phi - \Pi_{k_j}\phi \rangle, \quad \forall \phi\in \U,\\
\label{thm:conv_adaptive_pf02}
        \left| \langle \mathcal{R}_{1}(\phi_{k_j}^\ast,\bold{u}_{k_j}^\ast), \phi - \Pi_{k_j}\phi \rangle \right| &\leq c \eta_{k_{j},1}(\phi_{k_j}^\ast,\bold{u}^\ast_{k_j})\|\bold{\nabla}\phi\|_{L^2(\Omega)}, \quad \forall \phi \in \U.
     \end{align}
 By passing the limit $j\to\infty$ in \eqref{thm:conv_adaptive_pf01} and noting  $\lim_{j\to\infty}\eta_{k_{j},1}(\phi_{k_j}^\ast,\bold{u}^\ast_{k_j})= 0$ (cf. Theorem \ref{thm:conv_est_control}), we obtain
     \begin{equation}\label{thm:conv_adaptive_pf03}
        \liminf_{j\to\infty}\langle \mathcal{R}_{1}(\phi_{k_j}^\ast,\bold{u}_{k_j}^\ast), \phi - \phi_{k_j}^\ast \rangle\geq 0, \quad \forall \phi\in \U.
     \end{equation}
     Due to the convergence of $\{\bold{u}_{k_j}^\ast\}_{j\geq0}$ in Lemma \ref{lem:sol-map_cont_limit} and the discrete Sobolev inequality \eqref{disc_Sob_ineq_hom}, $\{\|\bold{u}_{k_j}^\ast\|_{\bold{L}^4(\Omega)}\}_{j\geq0}$ is uniformly bounded. Since $\phi_{k_j}$, $\overline{\phi}_\infty\in [0,1]$ a.e. in $\Omega$, pointwise convergence of $\phi_{k_j}\to \overline{\phi}_{\infty}$ a.e. in $\Omega$ (see \eqref{lem:sol-map_cont_limit_pf01}), the continuity of $\alpha'_{\eps}$ in $[0,1]$ and Lebesgue dominated convergence theorem imply
        \begin{align*}
       \lim_{j\to\infty} \left| ((\alpha_\eps'(\phi_{k_j}^\ast)\phi_{k_j}^\ast - \alpha_\eps'(\overline{\phi}_{\infty})\overline{\phi}_{\infty} ), |\bold{u}_{k_j}^\ast|^2) \right| & \leq \lim_{j\to\infty}\|\alpha_\eps'(\phi_{k_j}^\ast)\phi_{k_j}^\ast - \alpha_\eps'(\overline{\phi}_{\infty})\overline{\phi}_{\infty}\|_{L^2(\Omega)}\|\bold{u}_{k_j}^\ast\|_{\bold{L}^4(\Omega)}^2 =0.
        \end{align*}
By the $\bold{L}^2(\Omega)$ strong convergence of $\{\bold{u}_{k_j}^\ast\}_{j\geq0}$ again,
\[
    \lim_{j\to\infty} (\alpha_\eps'(\overline{\phi}_{\infty})\overline{\phi}_{\infty} , |\bold{u}_{k_j}^\ast|^2 ) = (\alpha_\eps'(\overline{\phi}_{\infty})\overline{\phi}_{\infty} , |\overline{\bold{u}}_\infty|^2 ).
\]
Thus we arrive at
\begin{equation*}
        \lim_{j\to\infty} ( \alpha_\eps'(\phi_{k_j}^\ast) |\bold{u}_{k_j}^\ast|^2,  \phi_{k_j}^\ast) = ( \alpha_\eps'(\overline{\phi}_{\infty}) |\overline{\bold{u}}_\infty|^2 , \overline{\phi}_{\infty}).
    \end{equation*}
Similarly, we deduce
$\lim_{j\to\infty} ( \alpha_\eps'(\phi_{k_j}^\ast) |\bold{u}_{k_j}^\ast|^2 ,  \phi ) = (\alpha_\eps'(\overline{\phi}_{\infty}) |\overline{\bold{u}}_\infty|^2 ,  \phi)$ for any $\phi \in \U$.
Then the $H^1(\Omega)$ convergence and the pointwise convergence of $\{\phi_{k_j}^\ast\}_{j\geq0}$ in Lemma \ref{lem:sol-map_cont_limit}
and Lebesgue dominated convergence theorem imply
\begin{align*}
        \lim_{j\to\infty}(\bold{\nabla}\phi_{k_j}^\ast, \bold{\nabla} (\phi - \phi_{k_j}^\ast)) &= \lim_{j\to\infty}(\bold{\nabla}\overline{\phi}_{\infty},
        \bold{\nabla}(\phi - \overline{\phi}_{\infty})), \quad \forall \phi \in \U,\\
        \lim_{j\to\infty}( f'(\phi_{k_j}^\ast) ,\phi - \phi_{k_j}^\ast) &= (f'(\overline{\phi}_{\infty}) ,\phi - \overline{\phi}_{\infty}), \quad \forall \phi \in \U.
    \end{align*}
    The proof of the theorem is completed by collecting the preceding estimates {and \eqref{thm:conv_adaptive_pf03}}.
\end{proof}

In the practical numerical treatment of problem \eqref{vp_stokes_div-free}, one approximates the pressure field $p$ by piecewise constants. This motivates the convergence analysis of the discrete pressure fields $\{p_k^\ast\in Q_k\}_{k\geq0}$ associated with minimizing pairs $\{(\phi_k^\ast,\bold{u}_k^\ast)\}_{k\geq 0}$ for problem \eqref{dismin_phase-field}-\eqref{disvp_stokes_div-free}.
Given the minimizing pair $(\phi^\ast,\bold{u}^\ast)\in \U\times \bold{U}$ to problem \eqref{vp_stokes_div-free}, a linear functional in  $\bold{V}^{0}:=\{\bold{l}\in \bold{H}^{-1}(\Omega) ~|~\langle \bold{l},\bold{v} \rangle = 0~\forall \bold{v}\in \bold{V}\}$ is given by $(\phi^\ast, \bold{u}^\ast) \in \U \times \bold{U}$ through \eqref{vp_stokes_div-free}. By the
inf-sup condition \cite[Corollary I.2.4 or (I.5.14) on p. 81]{GiraultRaviart:1986} and \cite[Lemma I.4.1]{GiraultRaviart:1986}, there exists a unique $p^\ast \in L_0^2(\Omega)$ such that
\begin{subequations}
    \begin{align}
    \mu (\bold{\nabla} \bold{u}^\ast,\bold{\nabla} \bold{v}) + (\alpha_\eps(\phi^\ast) \bold{u}^\ast, \bold{v}) - ( \mathrm{div} \bold{v}, p^\ast)&= (\bold{f}, \bold{v}), \quad \forall \bold{v}  \in \bold{H}_0^1(\Omega), \label{vp_stokes_motion}\\
    \label{vp_stokes_incompress}
       (\mathrm{div} \bold{u}^\ast ,q ) &= 0 ,\quad \forall q \in L_0^2(\Omega).
    \end{align}
\end{subequations}
Likewise under the discrete inf-sup condition \cite{CrouzeixRaviart:1973}, there exists a unique $p_k^\ast \in Q_k$ associated with the minimizing pair $(\phi_k^\ast,\bold{u}_k^\ast)\in \U_k \times \bold{U}_k$ of problem \eqref{dismin_phase-field}-\eqref{disvp_stokes_div-free} such that
\begin{subequations}\label{disvp_stokes}
    \begin{align}\label{disvp_stokes_motion}
    \mu ( \bold{\nabla} \bold{u}_k^\ast , \bold{\nabla} \bold{v}_k) + (\alpha_\eps(\phi^\ast) \bold{u}_k^\ast , \bold{v}_k) - (\mathrm{div}_k \bold{v}_k, p_k^\ast) &= ( \bold{f}, \bold{v}_k) ,\quad \forall \bold{v}_k  \in \bold{Z}_k,\\
   \label{disvp_stokes_incompress}
        (\mathrm{div}_k \bold{u}^\ast_k, q_k) &= 0 ,\quad \forall q_k \in Q_k.
    \end{align}
\end{subequations}

Then we have the following convergence result for discrete pressure fields generated by Algorithm \ref{alg_anfem_topopt-Stokes_phase-field_G1}.

\begin{theorem}\label{thm:conv_pre_adaptive}
    Let $\{p_k^\ast\in Q_k\}_{k \geq 0}$ be the sequence of discrete pressure fields associated with $\{(\phi_{k}^\ast,\bold{u}_{k}^\ast)\}_{k\geq0}$ generated by Algorithm \ref{alg_anfem_topopt-Stokes_phase-field_G1} through \eqref{disvp_stokes_motion}-\eqref{disvp_stokes_incompress}. Then the subsequence $\{p_{k_j}^\ast\}_{j\geq 0}$ associated with the convergent subsequence $\{(\phi_{k_j}^\ast,\bold{u}_{k_j}^\ast)\}_{j\geq0}$ in Theorem \ref{thm:conv_adaptive} converges strongly to $p^\ast \in L_0^2(\Omega)$ in \eqref{vp_stokes_motion} with respect to the strong $L^2(\Omega)$-topology.
\end{theorem}

\begin{proof}
The proof is identical with that of \cite[Theorem 3.3]{JinLiXu:2025}, except the reasoning for
$(\sum_{T\in \cT_{k_j}}h_T^2\|\alpha_\eps(\phi_{k_j}^*\bold{u}_{k_j}^*-\bold{f}\|^2_{L^2(T)})^{1/2}$ at \cite[(3.43)]{JinLiXu:2025}, which is replaced with Theorem \ref{thm:conv_est_state&costate}. In fact, the argument in \cite[Theorem 3.3]{JinLiXu:2025} yields
\begin{equation*}
    \|p^\ast - p^\ast_{k_j}\|_{L^2(\Omega)} \leq c \Big(\Big(\sum_{T\in \cT_k}h_T^2\|\alpha_\eps(\phi_{k_j}^\ast) \bold{u}_{k_j}^\ast - \bold{f}\|_{\bold{L}^2(T)}^2\Big)^{1/2} + \| \bold{u}^\ast - \bold{u}_{k_j}^\ast \|_{1,k_j} + \| \alpha_\eps(\phi^\ast) \bold{u}^\ast - \alpha_\eps(\phi_{k_j}^\ast)\bold{u}_{k_j}^\ast \|_{\bold{L}^2(\Omega)} \Big).
\end{equation*}
Then Theorem \ref{thm:conv_est_state&costate} ensures a vanishing limit of the first term, and the rest also tend to zero by Theorem \ref{thm:conv_adaptive}.
\end{proof}

\section{Numerical results and discussions}\label{sect:numerics}
In this section, we present four numerical experiments to illustrate the performance of the adaptive algorithm. All numerical simulations are conducted using a desktop with a 13th Gen Intel(R) Core(TM) i7-13700 (24 CPUs) and 32GB RAM, and a Stokes solver, coded using MATLAB as well as the built-in PDE Toolbox, in \cite{LarsonBengzon:2013} and the package iFEM \cite{Chen:2009} are utilized in the 2d and 3d examples, respectively. In problem \eqref{min_phase-field}-\eqref{vp_stokes_div-free}, we set the viscosity coefficient $\mu=1$, $\alpha_\eps(\phi)=10000(1-\phi)^2$, $\bold{f}=\bold{0}$ and both the relaxation parameter $\eps$ and the regularization parameter $\gamma$ to $1\times 10^{-2}$ except for Example \ref{exp:bypass}, for which $\eps=5\times 10^{-3}$ and $\gamma=1\times 10^{-1}$. Instead of the pure Dirichlet condition in \eqref{vp_stokes_div-free}, we impose a traction-free condition, namely $(\mu\bold{\nabla}\bold{u} - p \bold{I}_d) \cdot \bold{n}=\bold{0}$, at the outlet, since the case of a mixed boundary condition is practical. The jump terms $J_{F,2}(\bold{u}_k^\ast)$ in $\eta_{k,2}$ for boundary edges/faces $F$ at the outlet are not required in computation accordingly. The inlet boundary data is given below and $\bold{g}=\bold{0}$ on the remaining part. In the adaptive algorithm, we take $K=4$ (2d) and $K=5$ (3d) refinement steps and fix $(\theta_{1},\theta_{2})$ in the module \textsf{MARK} at $(0.4,0.6)$ for the 2d cases and $(0.8,0.8)$ for the 3d case, respectively:
\begin{equation*}
	\eta_{k,1}^2(\phi_k^\ast,\bold{u}_k^\ast,\mathcal{M}_k^1) \geq \theta_1 \eta_{k,1}^2 (\phi_k^\ast,\bold{u}_k^\ast)\quad\mbox{and}\quad
    		\eta_{k,2}^2(\phi_k^\ast,\bold{u}_k^\ast,\mathcal{M}_k^2) \geq \theta_2 \eta_{k,2}^2 (\phi_k^\ast,\bold{u}_k^\ast).
\end{equation*}
The $\textsf{OPTIMIZE}$ module in Algorithm \ref{alg_anfem_topopt-Stokes_phase-field_G1} employs an inner-outer iterative scheme \cite[Algorithm 2.1]{JinLiXu:2025} to minimize the augmented Lagrangian functional $\mathcal{L}(\phi_{k},\bold{u}_k):=\J^\eps(\phi_k) + \ell W(\phi_k) + \frac{\zeta}{2} W(\phi_k)^2$ (with $ W(\phi_k):=\int_\Omega \phi_k \dx - \beta |\Omega|$) of \eqref{dismin_phase-field}-\eqref{disvp_stokes_div-free} on each mesh $\cT_k$ (subject to \eqref{disvp_stokes_div-free}) with $\phi_k$ in the convex subset of $S_k$ restricted into $[0,1]$. The number $N$ of outer loops is set to $50$ for the 2d examples and $10$ for the 3d example for \eqref{disvp_stokes_div-free}. In each outer loop, we employ 10 inner loops to update the discrete phase field using an $L^2$-gradient flow based scheme due to \cite{LiYang:2022} with a stabilized coefficient $\tilde{S}=0.25$ or $1$ (for Example \ref{exp:bypass}). The initial Lagrangian multiplier $\ell_0$ and the initial penalty parameter $\zeta_0$ (both for the volume constraint) are taken to be zero and $100$ with an increasing factor $\kappa=1.1$, unless otherwise specified.

\begin{enumerate}[label=(\alph*)]
  \item (Left Inflow) The domain $\Omega$ is $(0,1)^2$. A left inflow with the inlet boundary condition $\bold{g}=\left(4y\left(1-y\right),0\right)^{\mathrm{T}}$ prescribed on the left side $x=0$, and the outlet boundary lies at $\{(x,y)~|~x=1, 0.3\leq y \leq 0.7\}$. Set $\Delta t=1\times 10^{-4}$ and $\beta=0.5$. \label{exp:leftinflow}

  \item (Three Inflows) The domain $\Omega$ is $(0,1)^2$. The three inflows are positioned in the middle of the top, bottom, and left edges of the square, each spanning a width of 0.2. The flow directions are downward, upward, and leftward, respectively, with a unit speed. The outlet is located on in the middle of the right edge with a spanning width of 0.2. Set $\Delta t= 5\times 10^{-5}$ and $\beta=0.36$.  \label{exp:threeinflows}

  \item (Bypass) The design domain $\Omega$ for the bypass is a rectangle $(0,1.5) \times(-0.5,0.5)$ with two inlets located at $\{(x,y)~|~x=0, 0.15 \leq y  \leq 0.35\} $ and $ \{(x,y)~|~ x=0,-0.35 \leq y \leq -0.15\} $ and two outlets at the same vertical positions along $x=1.5$. The flow through the inlet is given by $\boldsymbol{g}=\left(-100\left(y^2-0.35^2\right)\left(y^2-0.15^2\right), 0\right)^{\mathrm{T}}$. Set $\Delta t=5\times 10^{-3}$, $\eps=5\times 10^{-3}$, $\gamma=1\times 10^{-1}$, $\tilde{S}=1$ and $\zeta_{0}=50$. And the target volume is $0.7$. \label{exp:bypass}

\item (Pipe in 3d) The domain $\Omega$ is $(0,1)^3$. Two inlets represent circular regions at the both ends of $\Omega$ in the $x$-direction, positioned along the $x=0$ and $x=1$ planes, respectively. Each of these regions has a radius of 0.2, centered at the point $(y=0.5, z=0.7)$ in the $x z$-plane with $\bold{g}=\left(\pm 1,0,0\right)^{\mathrm{T}}$. Similarly, the remaining two inlets are located at the boundaries in the $y$-direction, positioned along the $y=0$ and $y=1$ planes, with their centers at $(x=0.5, z=0.7)$ in the $y z$-plane given $\bold{g}=\left(0,\pm 1,0\right)^{\mathrm{T}}$. The outlet is defined as a circular region in the $x y$-plane at $z=0$ with the same radius, centered at the point $(x=0.5, y=0.5)$. These regions collectively define the inlet and outlet boundaries for the fluid flow. Set $\beta=0.2$, $\zeta_{0}=200$ and $\Delta t=5\times10^{-4}$. \label{exp:pipe3d}
\end{enumerate}

\begin{figure}[hbt!]
\centering\setlength{\tabcolsep}{0pt}
\begin{tabular}{ccccc}
\includegraphics[width=0.2\textwidth]{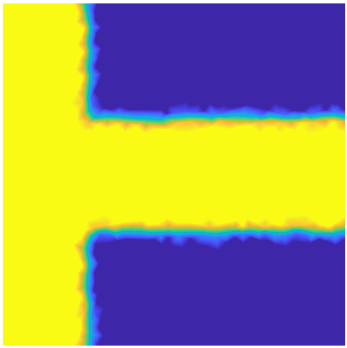} & \includegraphics[width=0.2\textwidth]{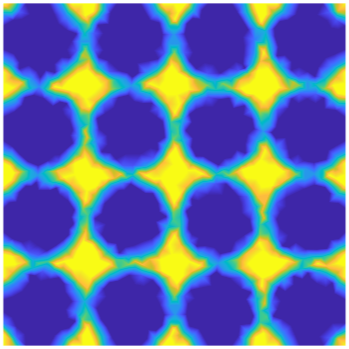}  &
\includegraphics[width=0.25\textwidth]{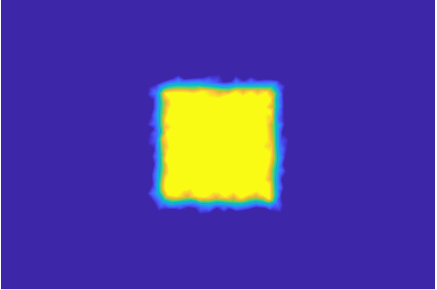}&
\includegraphics[width=0.31in]{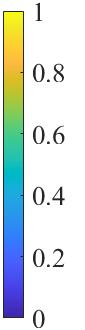} & \includegraphics[width=0.25\textwidth]{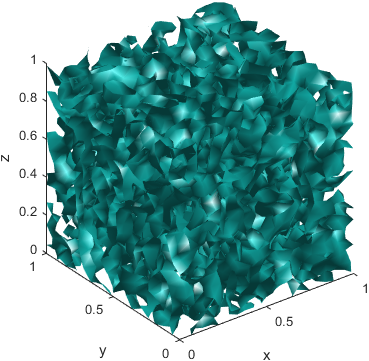}\\
\ref{exp:leftinflow} & \ref{exp:threeinflows} & \ref{exp:bypass} & & \ref{exp:pipe3d}
\end{tabular}
\caption{Initial phase-field functions for {Examples \ref{exp:leftinflow}-\ref{exp:pipe3d}}.
}\label{fig: initialshapes}
\end{figure}

Fig. \ref{fig: initialshapes} presents the initial phase-field functions over the initial meshes for the examples. The initial shapes, represented by discrete phase-field functions taking 1, are indicated in yellow for 2d examples. The computation of the 3d example begins with the one generated randomly with a uniform
distribution in the interval $[0, 1]$, whose iso-surface with a value of $0.5$ is shown in Fig. \ref{fig: initialshapes}(d).

\begin{figure}[hbt!]
	\centering\setlength{\tabcolsep}{0pt}
	\begin{tabular}{cccccccc}
		\includegraphics[width=1in]{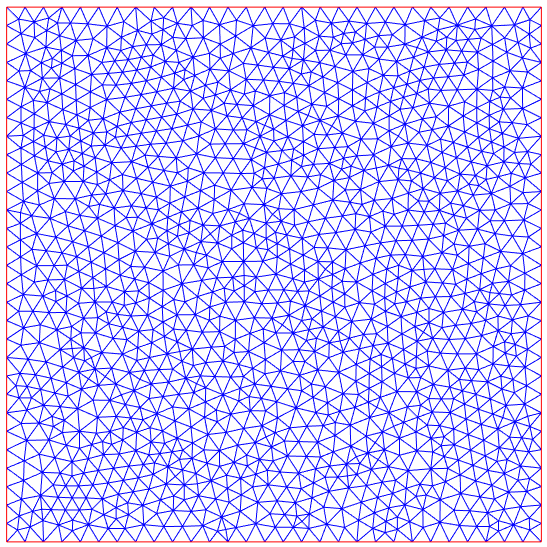}
        &\includegraphics[width=1in]{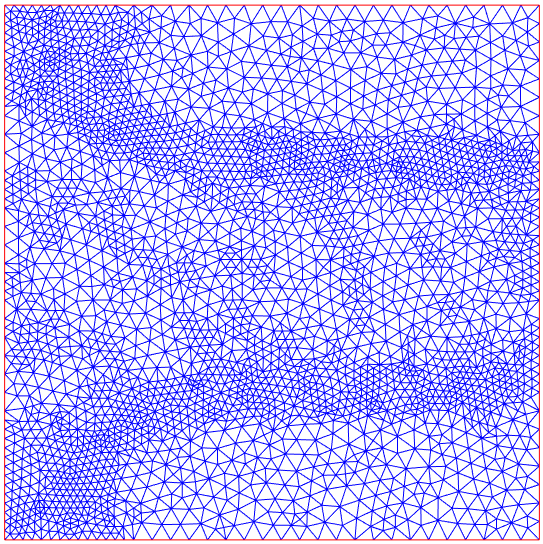}
		&\includegraphics[width=1in]{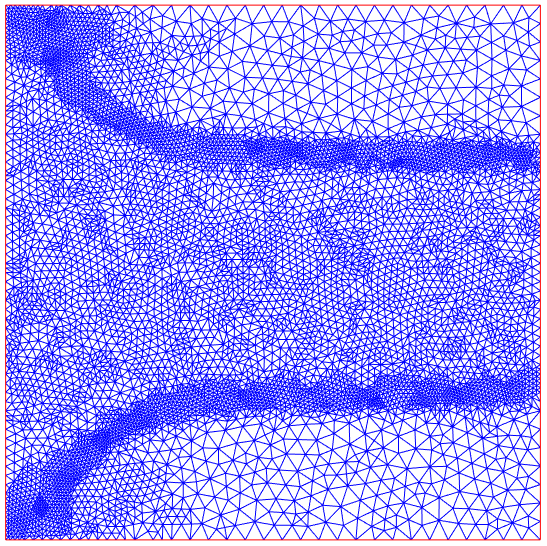}
        &\includegraphics[width=1in]{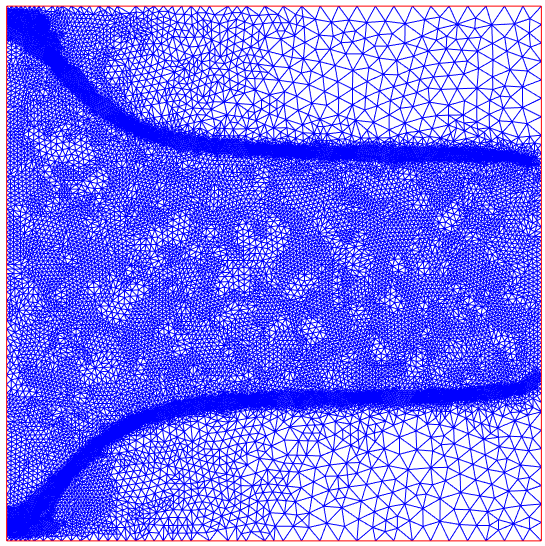}&\\
        \includegraphics[width=1in]{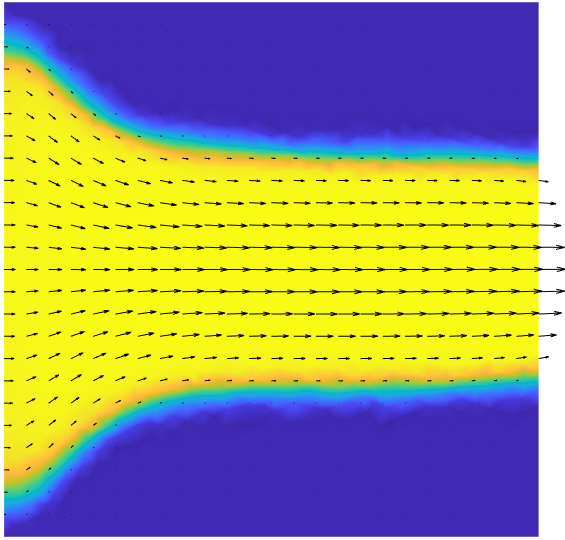}
        &\includegraphics[width=1in]{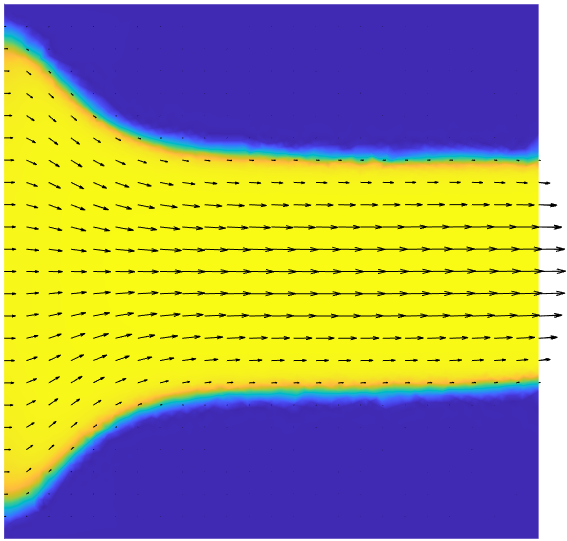}
        &\includegraphics[width=1in]{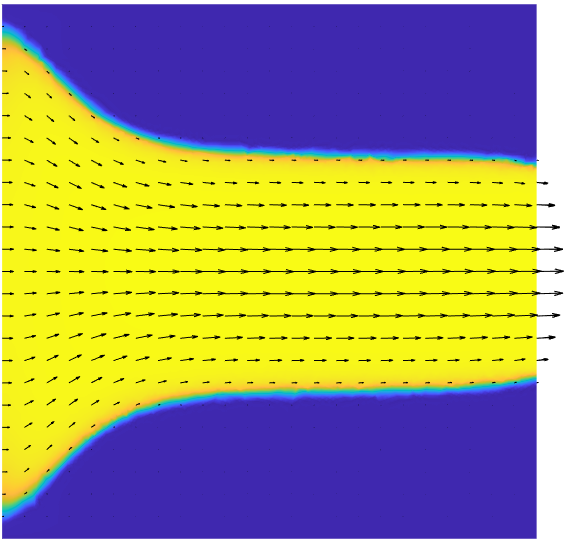}
        &\includegraphics[width=1in]{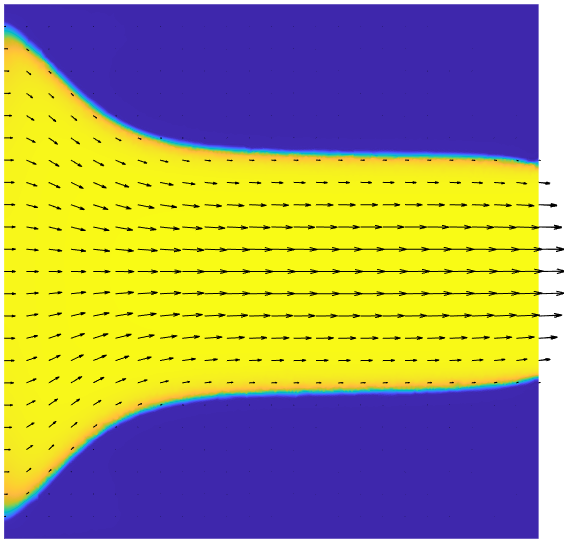}
        &\includegraphics[width=0.245in]{DesignColorbar.png}\\
        \includegraphics[width=1in]{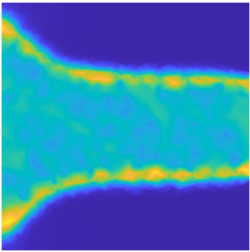}
        &\includegraphics[width=1in]{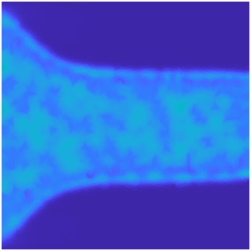}
        &\includegraphics[width=1in]{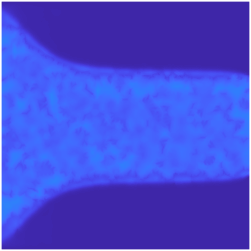}
        &\includegraphics[width=1in]{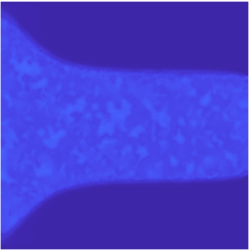}
        &\includegraphics[width=0.185in]{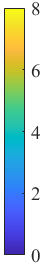}\\
        \includegraphics[width=1in]{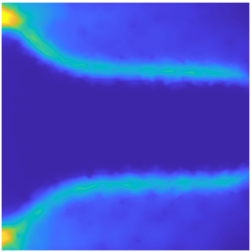}
        &\includegraphics[width=1in]{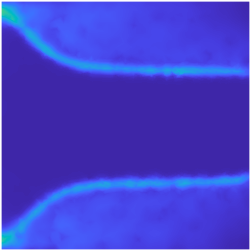}
        &\includegraphics[width=1in]{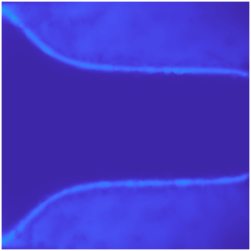}
        &\includegraphics[width=1in]{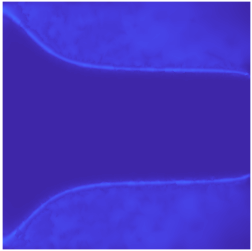}
        &\includegraphics[width=0.19in]{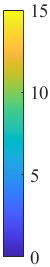}\\

        $k=0$ & $k=1$ & $k=2$ & $k=3$
	\end{tabular}
\caption{Numerical results for Example \ref{exp:leftinflow} from top to bottom: mesh, optimized design $\phi_k^*$, estimators $\eta_{k,1}$ and $\eta_{k,2}$. The number of vertices on each mesh is 1441, 2954, 6946 and 17354.
}\label{fig: LeftInflowAdapProcess}
\end{figure}

\begin{figure}[hbt!]
\centering\setlength{\tabcolsep}{0pt}
	\begin{tabular}{cccccccc}
		\includegraphics[width=1in]{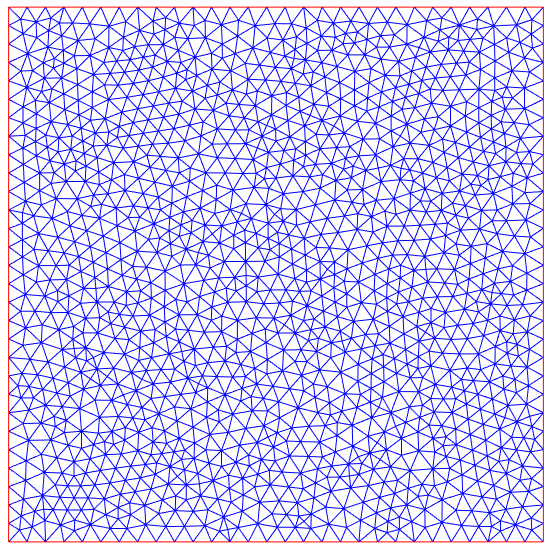}
        &\includegraphics[width=1in]{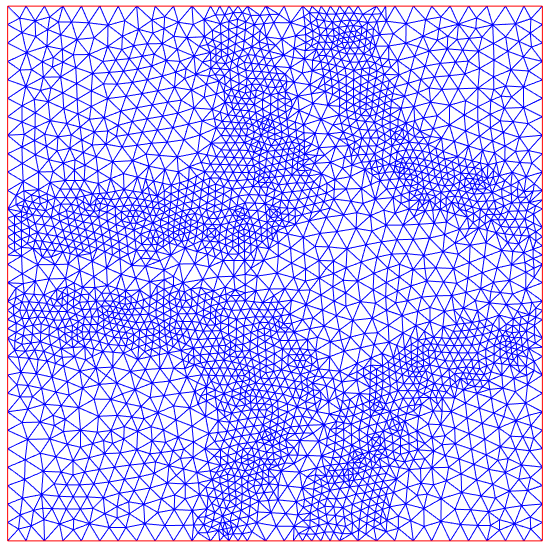}
		&\includegraphics[width=1in]{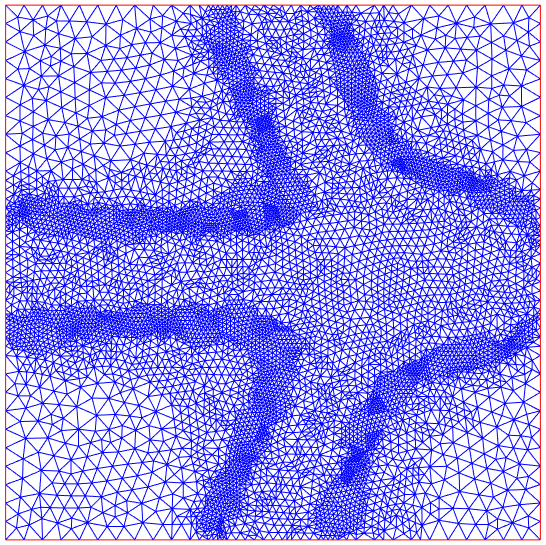}
        &\includegraphics[width=1in]{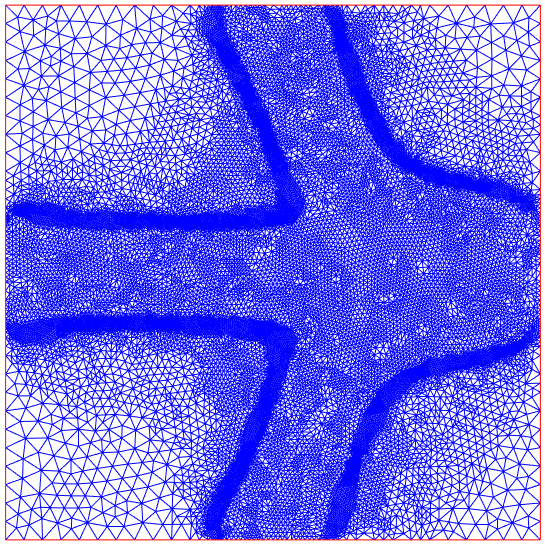}
        &
        \\
        \includegraphics[width=1in]{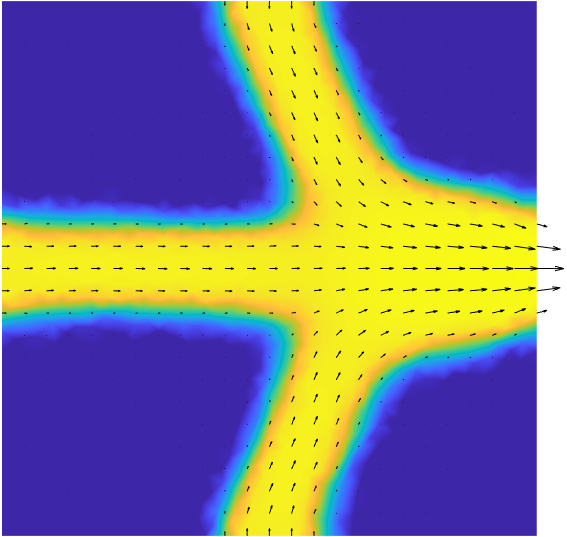}
        &\includegraphics[width=1in]{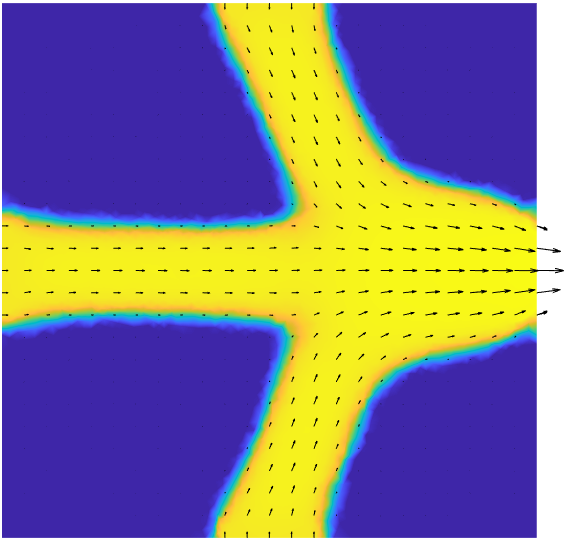}
        &\includegraphics[width=1in]{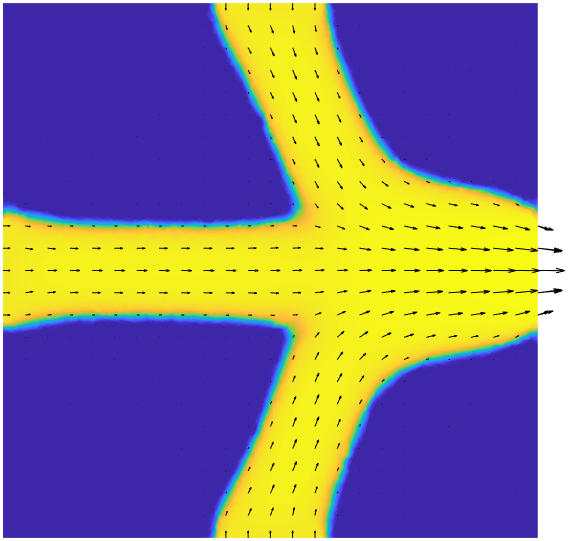}
        &\includegraphics[width=1in]{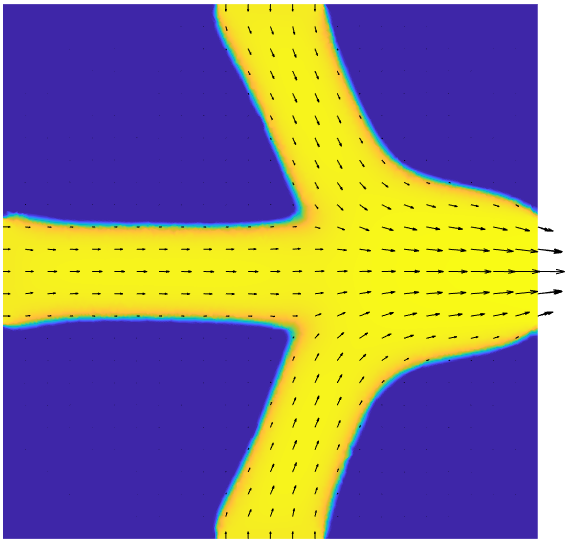}
        &\includegraphics[width=0.245in]{DesignColorbar.png}
        \\
        \includegraphics[width=1in]{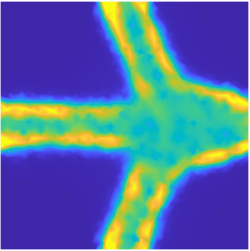}
        &\includegraphics[width=1in]{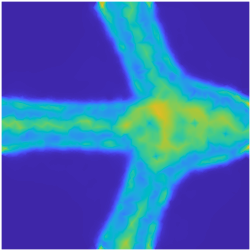}
        &\includegraphics[width=1in]{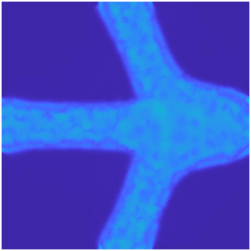}
        &\includegraphics[width=1in]{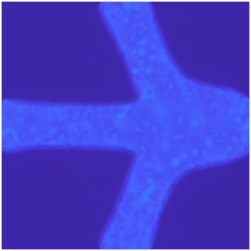}
        &\includegraphics[width=0.195in]{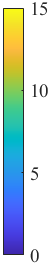}
        \\
        \includegraphics[width=1in]{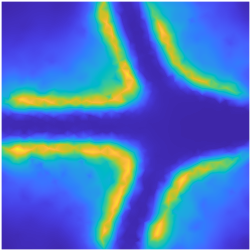}
        &\includegraphics[width=1in]{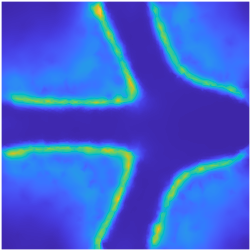}
        &\includegraphics[width=1in]{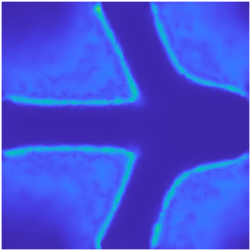}
        &\includegraphics[width=1in]{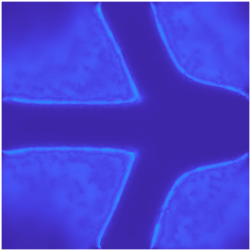}
        &\includegraphics[width=0.195in]{ThreeInflowsAdap1EstimatorColorbar.png}\\
        $k=0$ & $k=1$ & $k=2$ & $k=3$
	\end{tabular}
\caption{Numerical results for Example \ref{exp:threeinflows} from top to bottom: mesh, optimized designs $\phi_k^\ast$, and the estimators $\eta_{k,1}$ and $\eta_{k,2}$. The number of vertices on each mesh is 1441, 3070, 7442 and 19364.
}\label{fig: ThreeInflowsAdapProcess}
\end{figure}

\begin{figure}[hbt!]
	\centering\setlength{\tabcolsep}{0pt}
	\begin{tabular}{cccccccc}
		\includegraphics[width=1.3in]{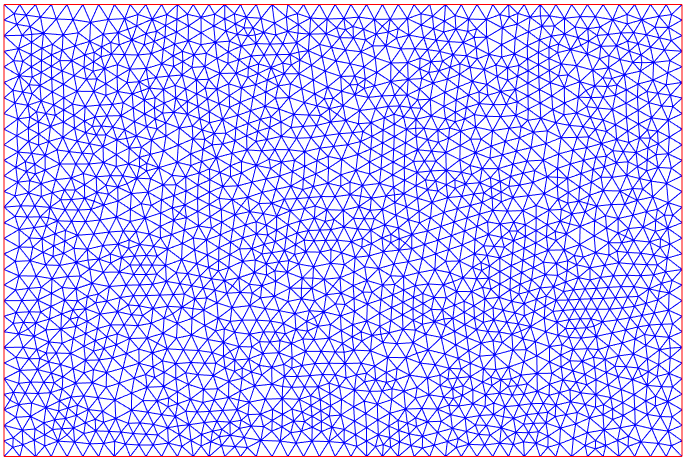}  \hspace{0.1cm}
        &\includegraphics[width=1.3in]{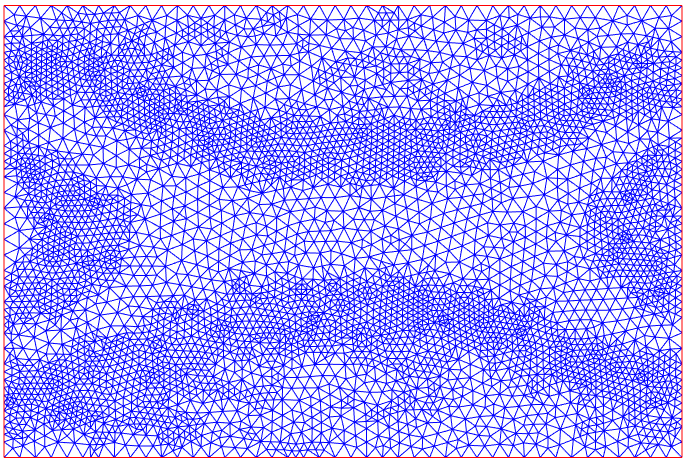} \hspace{0.1cm}
		&\includegraphics[width=1.3in]{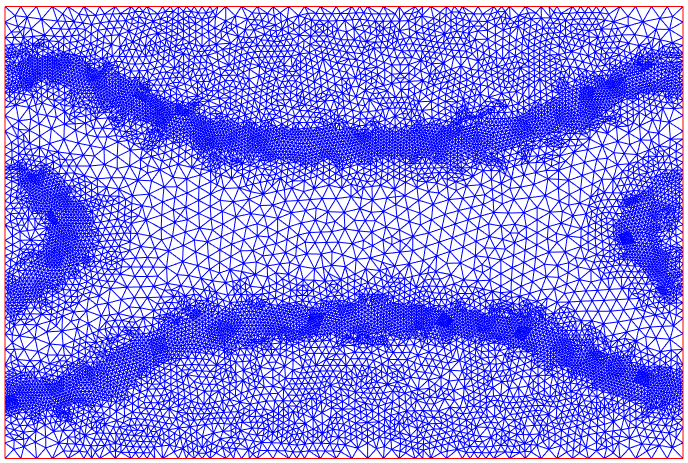} \hspace{0.1cm}
        &\includegraphics[width=1.3in]{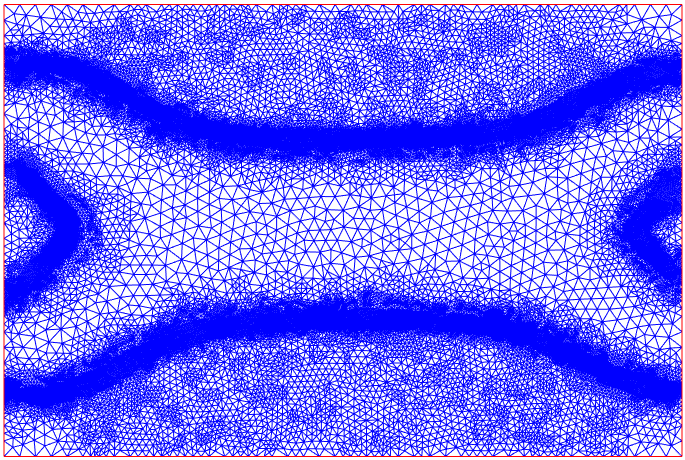}&  \\
        \includegraphics[width=1.3in]{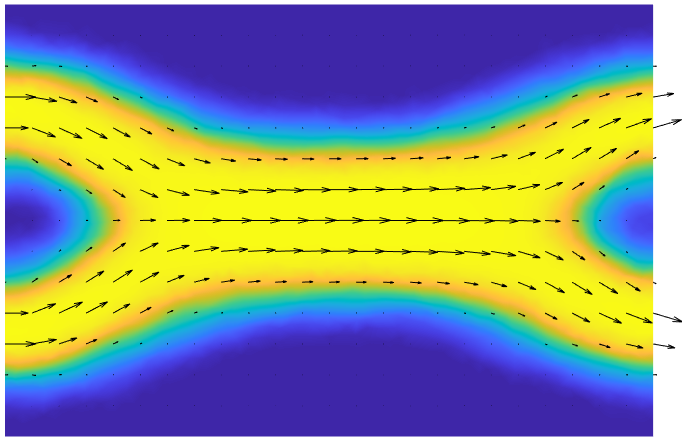}  \hspace{0.1cm}
        &\includegraphics[width=1.3in]{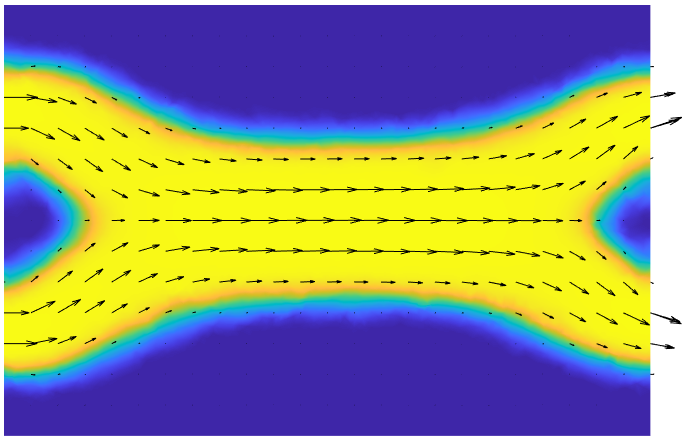} \hspace{0.1cm}
        &\includegraphics[width=1.3in]{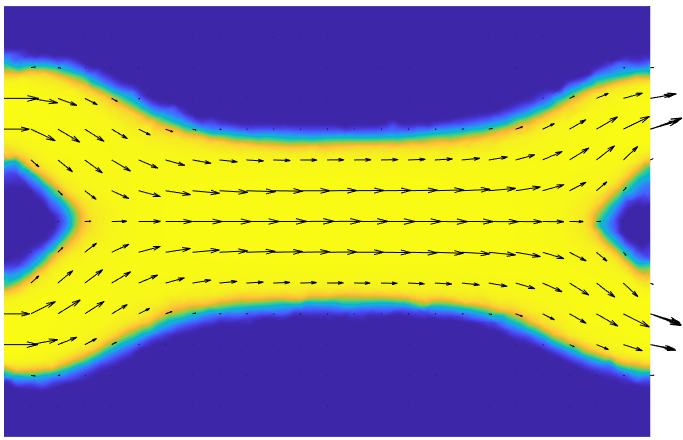} \hspace{0.1cm}
        &\includegraphics[width=1.3in]{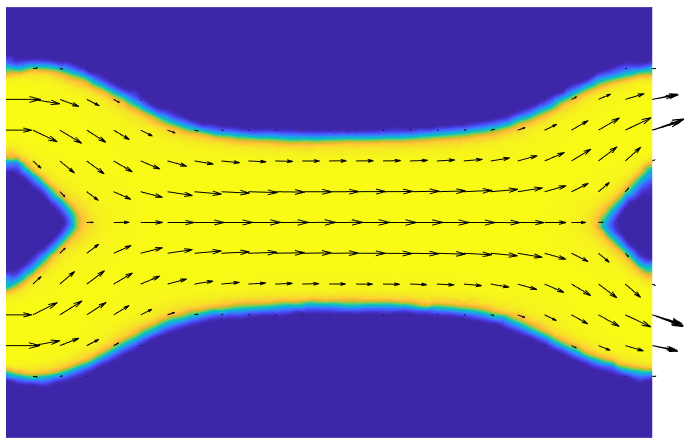}
        &\includegraphics[width=0.215in]{DesignColorbar.png}
        \\
        \includegraphics[width=1.3in]{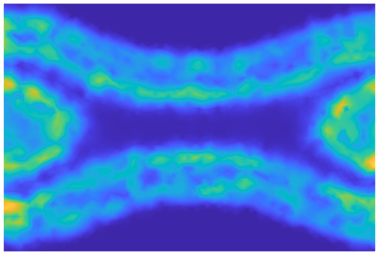} \hspace{0.1cm}
        &\includegraphics[width=1.3in]{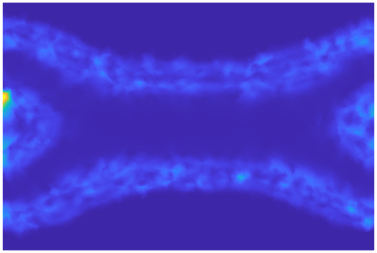} \hspace{0.1cm}
        &\includegraphics[width=1.3in]{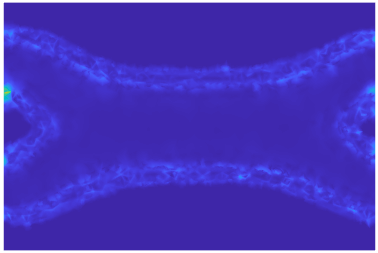} \hspace{0.1cm}
        &\includegraphics[width=1.3in]{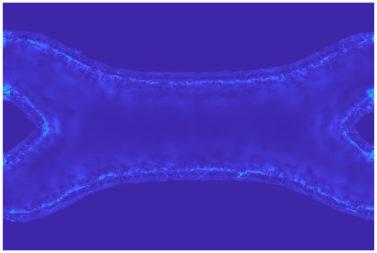}
        \hspace{0.0cm}
        &\includegraphics[width=0.21in]{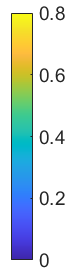}
        \\
        \includegraphics[width=1.3in]{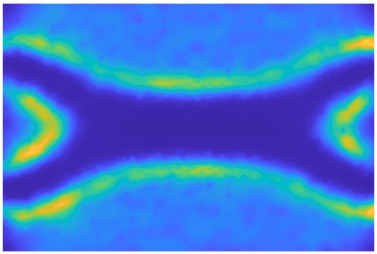} \hspace{0.1cm}
        &\includegraphics[width=1.3in]{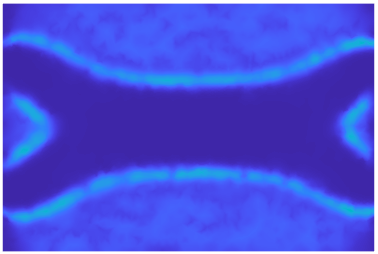} \hspace{0.1cm}
        &\includegraphics[width=1.3in]{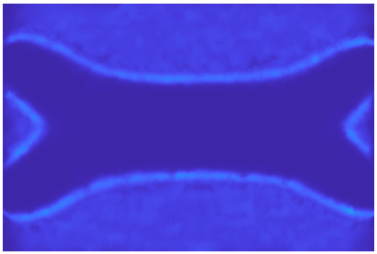} \hspace{0.1cm}
        &\includegraphics[width=1.3in]{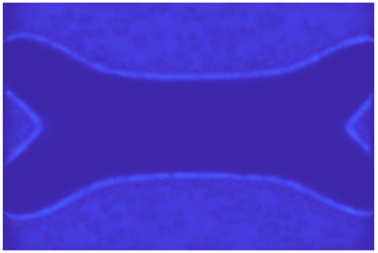}
        \hspace{0.0cm}
        &\includegraphics[width=0.2in]{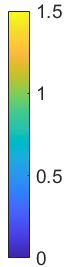}\\
        $k=0$ & $k=1$ & $k=2$ & $k=3$
	\end{tabular}
\caption{Numerical results for Example \ref{exp:bypass} from top to bottom: mesh, optimized designs $\phi_k^\ast$, and the estimators $\eta_{k,1}$ and $\eta_{k,2}$. The number of vertices on each mesh is 2174, 4973, 11669 and 28468. }\label{fig: ByPasssAdapProcess}
\end{figure}

Figs. \ref{fig: LeftInflowAdapProcess}, \ref{fig: ThreeInflowsAdapProcess}, and \ref{fig: ByPasssAdapProcess} present the results for the 2d examples, including the evolution of meshes $\cT_k$, optimal designs $\phi_k^*$ and the error estimators $\eta_{k,1}(\phi_k^*,\bm{u}_k^*,T)$ and $\eta_{k,2}(\phi_k^*,\bm{u}_k^*,T)$ during the adaptive refinement process. The black arrows indicate the direction of fluid flow. It is observed that as the adaptive refinement progresses, Algorithm \ref{alg_anfem_topopt-Stokes_phase-field_G1} can accurately capture the interfaces between the fluid and void regions, around which local refinements are concentrated. More precisely, the first adaptive iteration $(k=0)$ yields a rough shape of the design, and with subsequent adaptive refinements, the algorithm produces improved boundary interfaces of the fluid region (in yellow), due to the higher resolution of the velocity field $\bm{u}$. The additional refinements mainly take place around the interfaces during the adaptive refinement process, which clearly shows the effectiveness of the adaptive algorithm. The estimators $\eta_{k,1}$ and $\eta_{k,2}$ play different roles: $\eta_{k,1}$ indicates the presence of the fluid region, whereas $\eta_{k,2}$ concentrates in the void region due to the penalty term $\alpha_{\epsilon}(\phi)$. The magnitudes of the two error estimators are relatively large at the interface between the fluid and non-fluid regions, and decrease steadily to zero throughout the adaptive process.

\begin{figure}[hbt!]
\centering
\setlength{\tabcolsep}{0pt}
\begin{tabular}{cccc}
\includegraphics[width=0.25\textwidth]{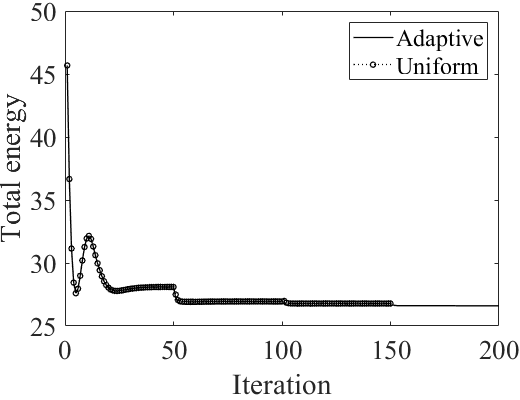}
&\includegraphics[width=0.25\textwidth]{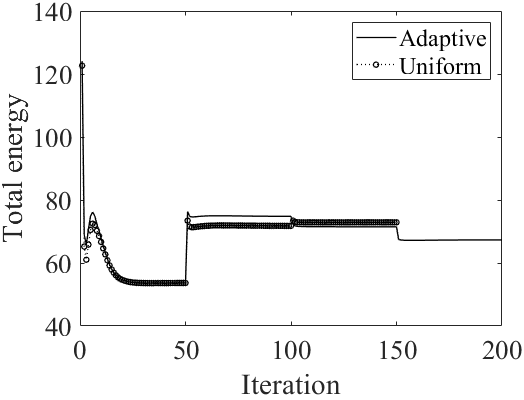}
&\includegraphics[width=0.25\textwidth]{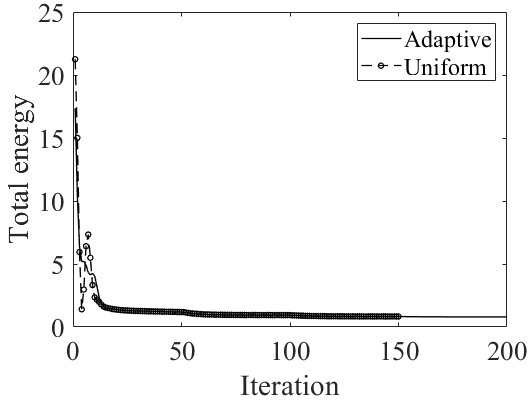}
&\includegraphics[width=0.25\textwidth]{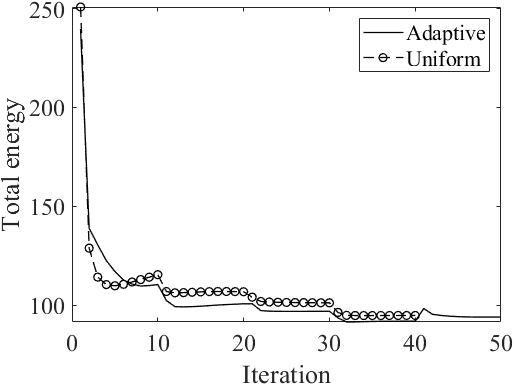}\\
\includegraphics[width=0.25\textwidth]{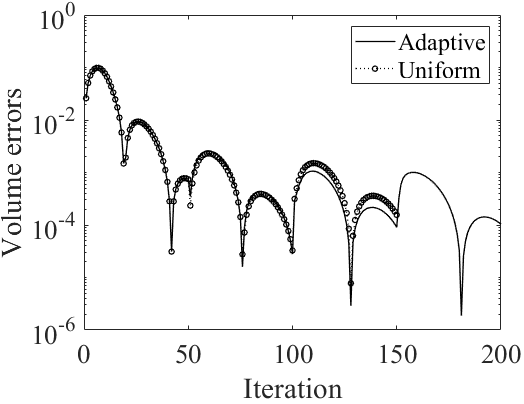} &
\includegraphics[width=0.25\textwidth]{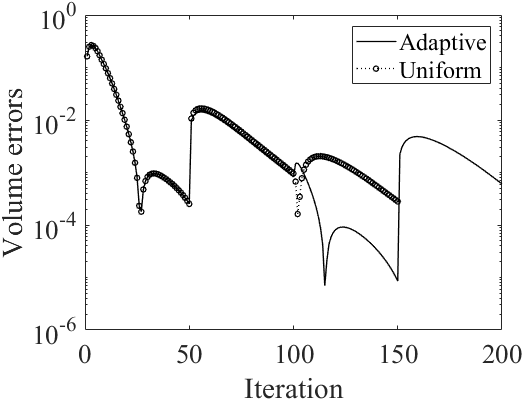}
&\includegraphics[width=0.25\textwidth]{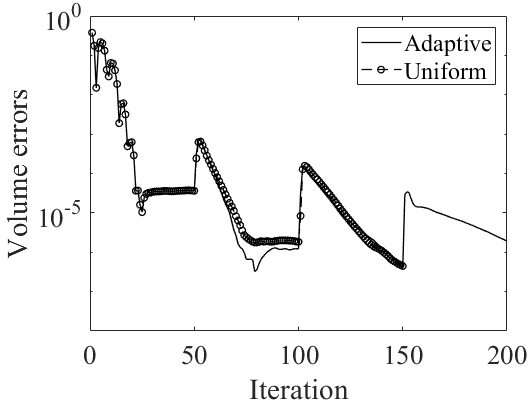} &\includegraphics[width=0.25\textwidth]{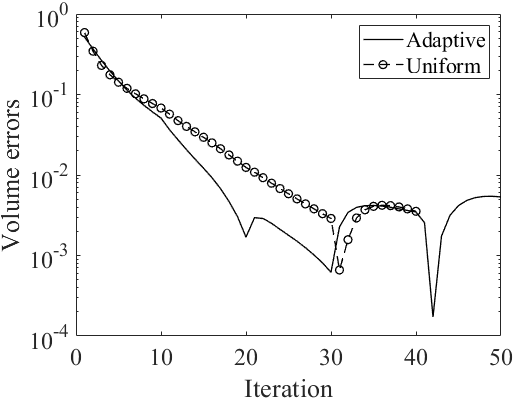}\\
\ref{exp:leftinflow} & \ref{exp:threeinflows} & \ref{exp:bypass} & \ref{exp:pipe3d}
\end{tabular}
\caption{The convergence history of the total energy (top) and the volume constraint error (bottom) versus the total number ($K \times N$) of outer iterations performed.}\label{fig:conv2D}
\end{figure}

\begin{figure}[hbt!]
	\centering
	\begin{tabular}{ccc}
        \includegraphics[width=1in]{LeftInflowAdap4thOptimizedDeisgnFlow.png}
        \includegraphics[width=1in]{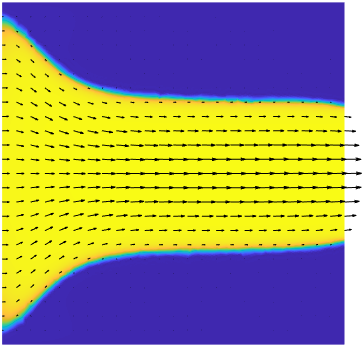}
        &
        \includegraphics[width=1in]{ThreeInflowsAdap4thOptimizedDeisgnFlow.png}
        \includegraphics[width=1in]{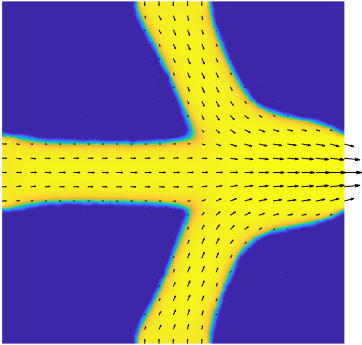} & \includegraphics[width=1.1in]{ByPassAdap4thOptimizedDeisgnFlow2.png}
        \includegraphics[width=1.1in]{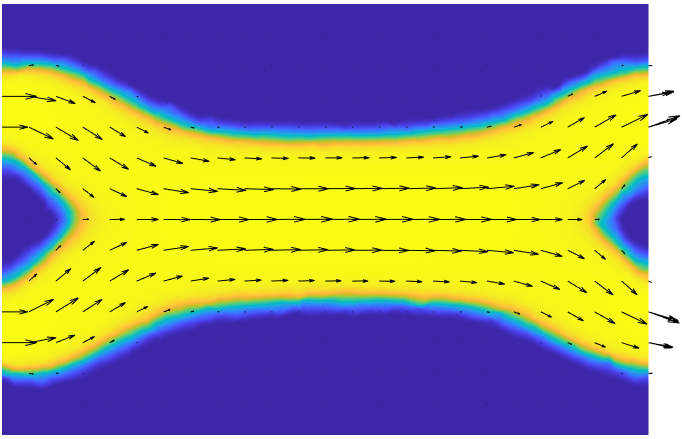} \\
          \ref{exp:leftinflow} & \ref{exp:threeinflows} &  \ref{exp:bypass}

	\end{tabular}
\caption{Comparison of optimal designs $\phi_{K-1}^\ast$ by adaptive (left within each pair) and uniform (right within each pair) refinements with optimal shapes indicated in yellow for Examples \ref{exp:leftinflow}-\ref{exp:bypass}.}\label{fig:designCompare}
\end{figure}

In Fig. \ref{fig:conv2D}, we compare the decay of the augmented Lagrangian functional $\mathcal{L}$ and the volume error for adaptive and uniform refinements, where the refinement steps and iteration numbers of the augmented Lagrangian method for adaptive and uniform refinements are $\left(K,N\right)=\left(4,50\right)$ and $\left(K,N\right)=\left(3,50\right)$, respectively. For both strategies, the functional values perform equally well, and the volume constraint is well preserved. The functional value decreases rapidly at the beginning and then stagnates. The magnitude of two final functional values are close to each other, which agrees with the fact that the two optimal designs are visually similar (see Fig. \ref{fig:designCompare}). The total energy experiences slight spikes when changing from one mesh to the other. This is due to the interpolation and projection errors when initializing the velocity and phase-field functions from the coarse mesh to the fine mesh.

Table \ref{tab:compare_results} presents quantitative results: the values of the Lagrangian functional augmented $\mathcal{L}(\phi_{K-1}^\ast,\bold{u}_{K-1}^\ast)$ (achieved in the final meshes generated by uniform and adaptive refinements), the number of vertices in the final mesh and the computing time (in seconds). The objective functional by the adaptive refinement strategy is very close to or slightly smaller than that for the uniform one with much less computing time. Despite the similarity between the optimized designs by two refinement strategies in Fig. \ref{fig:designCompare}, careful examination indicates that the profile given by the adaptive algorithm is smoother, with fewer small oscillations along the interface, than that by the uniform refinement. This is attributed to the {higher} resolution of the interface and the velocity field $\boldsymbol{u}$ by the adaptive strategy.

\begin{table}[hbt!]
\centering
\begin{threeparttable}
\caption{The comparison of quantitative results by adaptive and uniform refinements for Examples \ref{exp:leftinflow}-\ref{exp:pipe3d}: the number of vertices on the final mesh $\cT_{K-1}$, objective value and total computing time. \label{tab:compare_results}}
\begin{tabular}{|c|rrr|rrr|}
    \hline
    \multicolumn{1}{|c|}{\multirow{2}{*}{example}}
    & \multicolumn{3}{c}{adaptive}&\multicolumn{3}{|c|}{uniform}\\
    \cmidrule(lr){2-4} \cmidrule(lr){5-7}&
     vertices & objective & time (sec) & vertices & objective & time (sec) \\
    \hline
     \ref{exp:leftinflow}   & 17354  & 26.61   & 657  & 22345   & 26.79     & 19839  \\ \hline
     \ref{exp:threeinflows} & 19364  & 67.34   & 433  & 22345   & 72.91     & 1924 \\ \hline
     \ref{exp:bypass}       & 28468  & 0.80    & 6097 & 33905   & 0.84     & 11943  \\
     \hline
     \ref{exp:pipe3d}       & 24371  & 94.00    & 26246 & 42461   & 94.78     & 33140  \\
    \bottomrule
    \end{tabular}
\end{threeparttable}
\end{table}

\begin{figure}[hbt!]
	\centering\setlength{\tabcolsep}{0pt}
	\begin{tabular}{ccccccccc}
		\includegraphics[width=1.in]{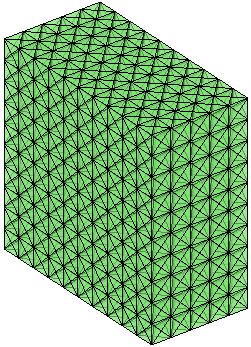}  \hspace{0.0cm}
        &\includegraphics[width=1.in]{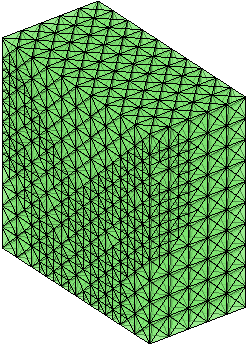} \hspace{0.0cm}
		&\includegraphics[width=1.in]{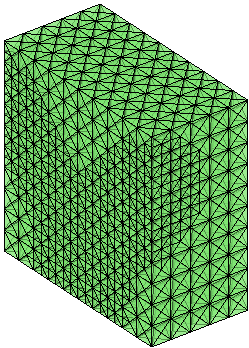} \hspace{0.0cm}
        &\includegraphics[width=1.in]{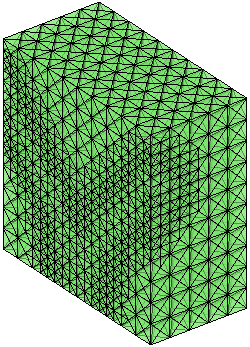}
        &\includegraphics[width=1.in]{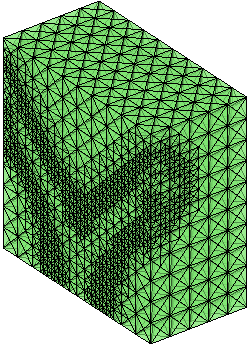}
        &
        \\
        \includegraphics[width=1.4in]{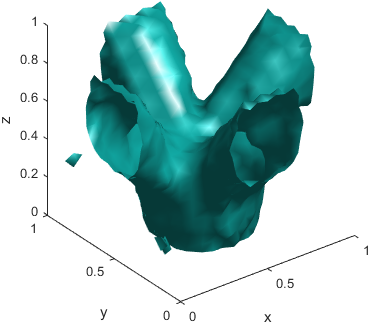}  \hspace{0.0cm}
        &\includegraphics[width=1.4in]{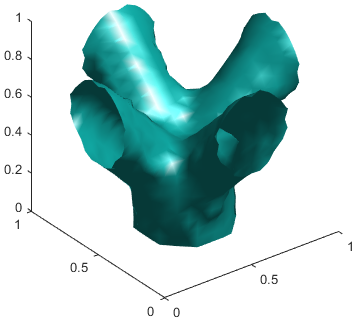} \hspace{0.0cm}
        &\includegraphics[width=1.4in]{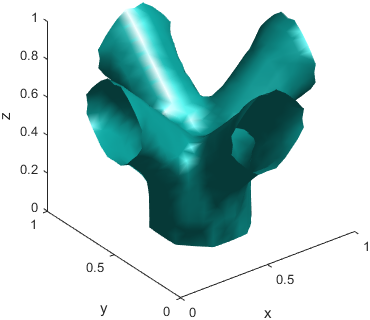} \hspace{0.0cm}
        &\includegraphics[width=1.4in]{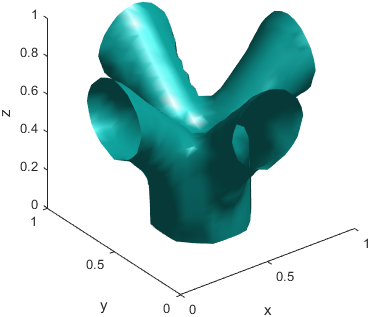}
        &\includegraphics[width=1.4in]{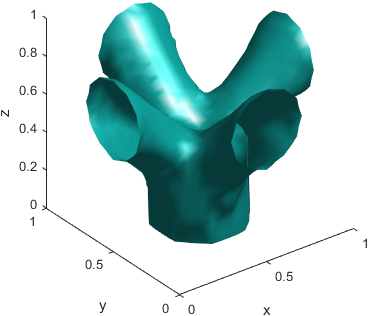}
        &
        \\
        \includegraphics[width=1.in]{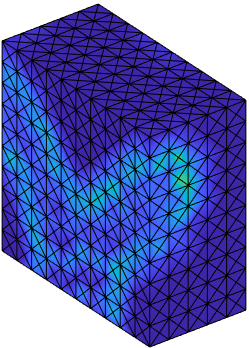} \hspace{0.0cm}
        &\includegraphics[width=1.in]{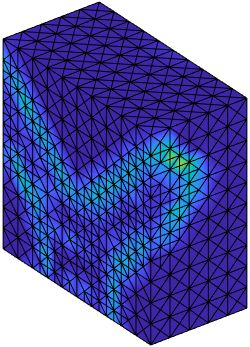} \hspace{0.0cm}
        &\includegraphics[width=1.in]{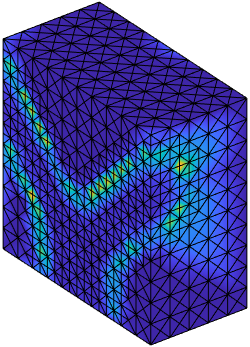} \hspace{0.0cm}
        &\includegraphics[width=1.in]{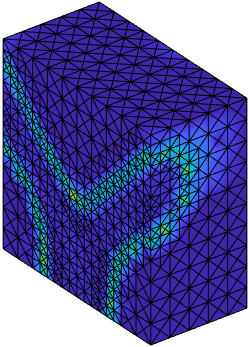}
        \hspace{0.0cm}
        &\includegraphics[width=1.in]{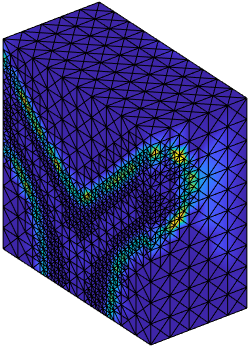}
        &\includegraphics[width=0.26in]{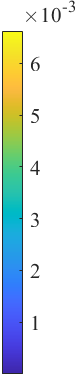}
        \\
        \includegraphics[width=1.in]{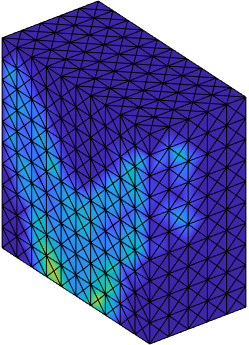} \hspace{0.0cm}
        &\includegraphics[width=1.in]{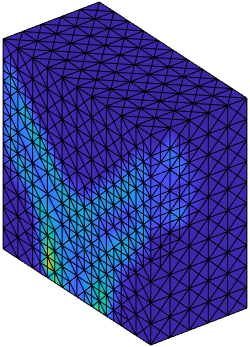} \hspace{0.0cm}
        &\includegraphics[width=1.in]{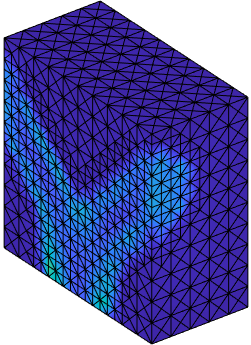} \hspace{0.0cm}
        &\includegraphics[width=1.in]{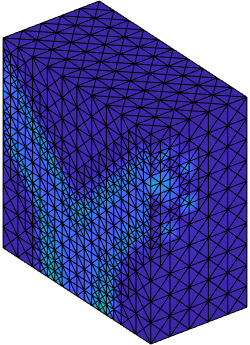}
        \hspace{0.0cm}
        &\includegraphics[width=1.in]{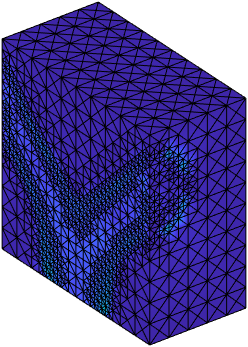}
        &\includegraphics[width=0.26in]{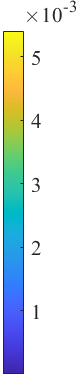}\\
        $k=0$ & $k=1$ & $k=2$ & $k=3$ & $k=4$
	\end{tabular}
\caption{Numerical results for Example \ref{exp:pipe3d}
from top to bottom: mesh, optimized designs $\phi_k^\ast$, and the estimators $\eta_{k,1}$ and $\eta_{k,2}$. The number of vertices on each mesh is 5631, 6929, 9631, 16052 and 24371.
}\label{fig: Pipe3DAdapProcess}
\end{figure}

\begin{figure}[hbt!]
	\centering
	\begin{tabular}{c}
        \includegraphics[width=1.5in]{isosurface4pipe3D.png}
        \includegraphics[width=1.5in]{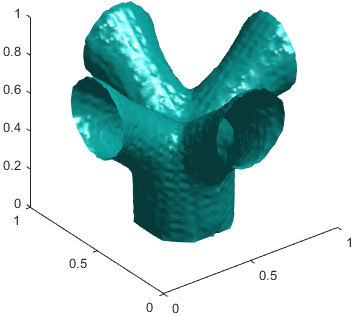}
	\end{tabular}
\caption{Comparison of optimal iso-surfaces with a value of $0.5$ by adaptive (left) and uniform (right) refinements with optimal shapes indicated in green for Example \ref{exp:pipe3d}.}\label{fig:designCompare3D}
\end{figure}

To further illustrate the performance of the adaptive method, we present one 3d example: the design of a pipe with four inlets and one outlet.
Fig. \ref{fig: Pipe3DAdapProcess} shows the evolution of mesh levels, optimized designs and  two estimators in the cross section $x>0.5$, in view of the symmetry during the adaptive refinement process. The adaptive strategy effectively identifies the diffuse interface of the phase-field functions, thereby allowing for precise local refinements. Meanwhile, as the mesh refinement process goes on, the optimized shapes become increasingly detailed and smooth, and the magnitude of the error estimators decreases stably. The convergence plots in Fig. \ref{fig:conv2D}(d) indicate the decay of the total energy and the good preservation of the volume constraint, where $(K,N)=(4,10)$ in the uniform refinement strategy. Similar to the 2d examples, the adaptive algorithm yields a better design represented by the iso-surface of $0.5$ (see Fig. \ref{fig:designCompare3D}) with less computing time (cf. Table \ref{tab:compare_results}) when compared with the uniform refinement strategy. 

\paragraph{Funding} The work of B. Jin is supported by Hong Kong RGC General Research Fund (Project 14306824), and a start-up fund from The Chinese University of Hong Kong. The work of Y. Xu is supported by National Natural Science Foundation of China under grants 12250013, 12261160361 and 12271367. The work of S. Zhu is supported by National Natural Science Foundation of China (12471377) and Science and the Technology Commission of Shanghai Municipality (22ZR1421900 and 22DZ2229014).

\bibliographystyle{abbrv}
\bibliography{top_fluid}

\begin{thebibliography}{10}

\bibitem{Age:2008}
N.~Aage, T.~H. Poulsen, A.~Gersborg-Hansen, and O.~Sigmund.
\newblock Topology optimization of large scale {Stokes} flow problems.
\newblock {\em Struct. Multidiscip. Optim.}, 35(2):175-- 180, 2008.

\bibitem{AinsworthOden:2000}
M.~Ainsworth and J.~T. Oden.
\newblock {\em A {P}osteriori {E}rror {E}stimation in {F}inite {E}lement
  {A}nalysis}.
\newblock Wiley-Interscience, John Wiley \& Sons, New York, 2000.

\bibitem{Aleandersen:2020}
J.~Alexandersen and C.~S. Andreasen.
\newblock A review of topology optimisation for fluid-based problems.
\newblock {\em Fluids}, 5(1):29, 2020.

\bibitem{Alonso:2018}
D.~H. Alonso, L.~F.~N. de~S\'{a}, J.~S.~R. Saenz, and E.~C.~N. Silva.
\newblock Topology optimization applied to the design of {2D} swirl flow
  devices.
\newblock {\em Struct. Multidiscip. Optim.}, 58(6):2341--2364, 2018.

\bibitem{Alonso:2020}
D.~H. Alonso, J.~S.~R. Saenz, and E.~C.~N. Silva.
\newblock Non-newtonian laminar {2D} swirl flow design by the topology
  optimization method.
\newblock {\em Struct. Multidiscip. Optim.}, 62(1):299--321, 2020.

\bibitem{BeckerMao:2011}
R.~Becker and S.~Mao.
\newblock Quasi-optimality of adaptive nonconforming finite element methods for
  the {S}tokes equations.
\newblock {\em SIAM J. Numer. Anal.}, 49(3):970--991, 2011.

\bibitem{BonitoCanutoNochettoVeeser:2024}
A.~Bonito, C.~Canuto, R.~H. Nochetto, and A.~Veeser.
\newblock Adaptive finite element methods.
\newblock {\em Acta Numer.}, 33:163--485, 2024.

\bibitem{BorrvallPetersson:2003}
T.~Borrvall and J.~Petersson.
\newblock Topology optimization of fluids in {S}tokes flow.
\newblock {\em Int. J. Numer. Meth. Fluids}, 41(1):77--107, 2003.

\bibitem{Brenner:2003}
S.~C. Brenner.
\newblock Poincar\'{e}-{F}riedrichs inequalities for piecewise {$H^{1}$}
  functions.
\newblock {\em SIAM J. Numer. Anal.}, 41(1):306--324, 2003.

\bibitem{BrennerScott:2008}
S.~C. Brenner and L.~R. Scott.
\newblock {\em The {M}athematical {T}heory of {F}inite {E}lement {M}ethods}.
\newblock Springer, New York, third edition, 2008.

\bibitem{BuffaOrtner:2009}
A.~Buffa and C.~Ortner.
\newblock Compact embeddings of broken {S}obolev spaces and applications.
\newblock {\em IMA J. Numer. Anal.}, 29(4):827--855, 2009.

\bibitem{CarstensenFeischlPagePraetorius:2014}
C.~Carstensen, M.~Feischl, M.~Page, and D.~Praetorius.
\newblock Axioms of adaptivity.
\newblock {\em Comp. Math. Appl.}, 67(6):1195--1253, 2014.

\bibitem{CarstensenPeterseimRabus:2013}
C.~Carstensen, D.~Peterseim, and H.~Rabus.
\newblock Optimal adaptive nonconforming fem for the {S}tokes problem.
\newblock {\em Numer. Math.}, 123(2):291--308, 2013.

\bibitem{CKNS:2008}
J.~M. Cascon, C.~Kreuzer, R.~H. Nochetto, and K.~G. Siebert.
\newblock Quasi-optimal convergence rate for an adaptive finite element method.
\newblock {\em SIAM J. Numer. Anal.}, 46(5):2524--2550, 2008.

\bibitem{Chen:2009}
L.~Chen.
\newblock {\textit{iFEM}: an integrated finite element method package in
  MATLAB}.
\newblock Technical report, University of California at Irvine, 2009.
\newblock Available at \url{https://github.com/lyc102/ifem}.

\bibitem{Ciarlet:2002}
P.~G. Ciarlet.
\newblock {\em The {F}inite {E}lement {M}ethod for {E}lliptic {P}roblems}.
\newblock SIAM, Philadelphia, PA, 2002.

\bibitem{CrouzeixRaviart:1973}
M.~Crouzeix and P.-A. Raviart.
\newblock Conforming and nonconforming finite element methods for solving the
  stationary {S}tokes equations.
\newblock {\em RAIRO Anal. Num\'{e}r.}, 7(3):33--76, 1973.

\bibitem{DariDuranPadra:1995}
E.~Dari, R.~Dur\'{a}n, and C.~Padra.
\newblock Error estimators for nonconforming finite element approximations of
  the {S}tokes problem.
\newblock {\em Math. Comp.}, 64(211):1017--1033, 1995.

\bibitem{Deng:2018}
Y.~Deng, Y.~Wu, and Z.~Liu.
\newblock {\em {Topology Optimization Theory for Laminar Flow: Applications in
  Inverse Design of Microfluidics}}.
\newblock Springer, Singapore, 2018.

\bibitem{EvansGariepy:2015}
L.~C. Evans and R.~F. Gariepy.
\newblock {\em Measure Theory and Fine Properties of Functions}.
\newblock CRC Press, Boca Raton, FL, revised edition, 2015.

\bibitem{Evgrafov:2005}
A.~Evgrafov.
\newblock The limits of porous materials in the topology optimization of
  {S}tokes flows.
\newblock {\em Appl. Math. Optim.}, 52(3):263--277, 2005.

\bibitem{GantnerPraetorius:2022}
G.~Gantner and D.~Praetorius.
\newblock Plain convergence of adaptive algorithms without exploiting
  reliability and efficiency.
\newblock {\em IMA J. Numer. Anal.}, 42(2):1434--1453, 2022.

\bibitem{GarckeHecht:2015}
H.~Garcke and C.~Hecht.
\newblock A phase field approach for shape and topology optimization in
  {S}tokes flow.
\newblock In A.~Pratelli and G.~Leugering, editors, {\em New Trends in Shape
  Optimization}, pages 103--116. Birkh\"{a}user, New York, 2015.

\bibitem{GarckeHecht:2016a}
H.~Garcke and C.~Hecht.
\newblock Applying a phase field approach for shape optimization of a
  stationary {N}avier-{S}tokes flow.
\newblock {\em ESAIM Control Optim. Calc. Var.}, 22(2):309--337, 2016.

\bibitem{GarckeHecht:2016b}
H.~Garcke and C.~Hecht.
\newblock Shape and topology optimization in {S}tokes flow with a phase field
  approach.
\newblock {\em Appl. Math. Optim.}, 73(1):23--70, 2016.

\bibitem{GarckeHechtHinzeKahle:2015}
H.~Garcke, C.~Hecht, M.~Hinze, and C.~Kahle.
\newblock Numerical approximation of phase field based shape and topology
  optimization for fluids.
\newblock {\em SIAM J. Sci. Comput.}, 37(4):A1846--A1871, 2015.

\bibitem{Gersborg-Hansen:2005}
A.~Gersborg-Hansen, O.~Sigmund, and R.~B. Haber.
\newblock Topology optimization of channel flow problems.
\newblock {\em Struct. Multidiscip. Optim.}, 30(3):181--192, 2005.

\bibitem{GiraultRaviart:1986}
V.~Girault and P.~A. Raviart.
\newblock {\em Finite Element Methods for {N}avier-{S}tokes Equations}.
\newblock Springer-Verlag, Berlin, 1986.

\bibitem{GuestPrevost:2006}
J.~K. Guest and J.~H. Pr\'{e}vost.
\newblock Topology optimization of creeping fluid flows using a
  {D}arcy-{S}tokes finite element.
\newblock {\em Int. J. Numer. Meth. Engng.}, 66(3):461--484, 2006.

\bibitem{HuXu:2013}
J.~Hu and J.~Xu.
\newblock Convergence and optimality of the adaptive nonconforming linear
  element method for the {S}tokes problem.
\newblock {\em J. Sci. Comput.}, 55:125--148, 2013.

\bibitem{JinLiXuZhu:2024}
B.~Jin, J.~Li, Y.~Xu, and S.~Zhu.
\newblock An adaptive phase-field method for structural topology optimization.
\newblock {\em J. Comput. Phys.}, 206:112932, 27 pp., 2024.

\bibitem{JinLiXu:2025}
B.~Jin, J.~Li, Y.~Xu, and S.~Zhu.
\newblock On the {Crouzeix-Raviart} finite element approximation of phase-field
  dependent topology optimization in {Stokes} flow.
\newblock Preprint, arXiv:2505.04120, 2025.

\bibitem{JinWangXu:2025}
B.~Jin, F.~Wang, and Y.~Xu.
\newblock Adaptive approximations of inclusions in a semilinear elliptic
  problem related to cardiac electrophysiology.
\newblock {\em IMA J. Numer. Anal.}, pages in press, available as
  arXiv:2504.04483, 2025.

\bibitem{JinXu:2019}
B.~Jin and Y.~Xu.
\newblock Adaptive reconstruction for electrical impedance tomography with a
  piecewise constant conductivity.
\newblock {\em Inverse Problems}, 36(1):014003, 28 pp., 2020.

\bibitem{Kossaczky:1995}
I.~Kossaczky.
\newblock A recursive approach to local mesh refinement in two and three
  dimensions.
\newblock {\em J. Comp. Appl. Math.}, 55(3):275--288, 1994.

\bibitem{Kreissl:2011}
S.~Kreissl, G.~Pingen, and K.~Maute.
\newblock Topology optimization for unsteady flow.
\newblock {\em Int. J. Numer. Methods Engrg.}, 87(13):1229--1253, 2011.

\bibitem{LarsonBengzon:2013}
M.~G. Larson and F.~Bengzon.
\newblock {\em The Finite Element Method: Theory, Implementation, and
  Applications}.
\newblock Springer, Heidelberg, 2013.

\bibitem{LiYang:2022}
F.~Li and J.~Yang.
\newblock A provably efficient monotonic-decreasing algorithm for shape
  optimization in {S}tokes flows by phase-field approaches.
\newblock {\em Comput. Methods Appl. Mech. Engrg.}, 398:115195, 24 pp., 2022.

\bibitem{LiXuZhu:2023}
J.~Li, Y.~Xu, and S.~Zhu.
\newblock Adaptive computation of elliptic eigenvalue topology optimization
  with a phase-field approach.
\newblock {\em J. Numer. Math.}, 2025.
\newblock \url{https://doi.org/10.1515/jnma-2024-0060}.

\bibitem{MorinNochettoSiebert:2000}
P.~Morin, P.~H. Nochetto, and K.~G. Siebert.
\newblock Data oscillation and convergence of adaptive fem.
\newblock {\em SIAM J. Numer. Anal.}, 38(2):466--488, 2000.

\bibitem{MorinSiebert:2008}
P.~Morin, K.~G. Siebert, and A.~Veeser.
\newblock A basic convergence result for conforming adaptive finite elements.
\newblock {\em Math. Models Methods Appl. Sci.}, 18(5):707--737, 2008.

\bibitem{NieldBejan:2017}
D.~A. Nield and A.~Bejan.
\newblock {\em {Convection in Porous Media}}.
\newblock Springer, Cham., fifth edition, 2017.

\bibitem{NochettoSiebertVeeser:2009}
R.~H. Nochetto, K.~G. Siebert, and A.~Veeser.
\newblock Theory of adaptive finite element methods: an introduction.
\newblock In {\em {Multiscale, Nonlinear and Adaptive Approximation}}, pages
  409--542. Springer, Berlin, 2009.

\bibitem{Olesen:2006}
L.~H. Olesen, F.~Okkels, and H.~Bruus.
\newblock A high-level programming-language implementation of topology
  optimization applied to steady-state {Navier--Stokes} flow.
\newblock {\em Int. J. Numer. Methods Engrg.}, 65(7):975--1001, 2006.

\bibitem{OrtnerPraetorius:2011}
C.~Ortner and D.~Praetorius.
\newblock On the convergence of adaptive nonconforming finite element methods
  for a class of convex variational problems.
\newblock {\em SIAM J. Numer. Anal.}, 49(1):346--367, 2011.

\bibitem{Papadopoulos:2022}
I.~P.~A. Papadopoulos.
\newblock Numerical analysis of a discontinuous {G}alerkin method for the
  {B}orrvall-{P}etersson topology optimization problem.
\newblock {\em SIAM J. Numer. Anal.}, 60(5):2538--2564, 2022.

\bibitem{PapadopoulosSuli:2022}
I.~P.~A. Papadopoulos and E.~S\"{u}li.
\newblock Numerical analysis of a topology optimization problem for {S}tokes
  flow.
\newblock {\em J. Comput. Appl. Math.}, 412:114295, 21 pp., 2022.

\bibitem{Sa:2018}
L.~F.~N. S\'{a}, J.~S. Romero, O.~Horikawa, and E.~C.~N. Silva.
\newblock Topology optimization applied to the development of small scale pump.
\newblock {\em Struct. Multidiscip. Optim.}, 57(5):2045--2059, 2018.

\bibitem{ShiWang:2013}
Z.~Shi and M.~Wang.
\newblock {\em Finite Element Methods}.
\newblock Science Press, Beijing, 2013.

\bibitem{Siebert:2011}
K.~G. Siebert.
\newblock A convergence proof for adaptive finite elements without lower bound.
\newblock {\em IMA J. Numer. Anal.}, 31(3):947--970, 2011.

\bibitem{Stummel:1980}
F.~Stummel.
\newblock Basic compactness properties of nonconforming and hybrid finite
  element spaces.
\newblock {\em RAIRO Anal. Num\'{e}r.}, 14(1):81--115, 1980.

\bibitem{Thore:2021}
C.-J. Thore.
\newblock Topology optimization of {S}tokes flow with traction boundary
  conditions using low-order finite elements.
\newblock {\em Comput. Methods Appl. Mech. Engrg.}, 386:114177, 14 pp., 2021.

\bibitem{Traxler:1997}
C.~Traxler.
\newblock An algorithm for adaptive mesh refinement in $n$ dimensions.
\newblock {\em Computing}, 59(2):115--137, 1997.

\bibitem{Verfurth:2013}
R.~Verf\"{u}rth.
\newblock {\em A {P}osteriori {E}rror {E}stimation {T}echniques for {F}inite
  {E}lement {M}ethods}.
\newblock Oxford University Press, Oxford, 2013.

\bibitem{WikerKlarbringBorrvall:2007}
N.~Wiker, A.~Klarbring, and T.~Borrvall.
\newblock Topology optimization of regions of {Darcy and Stokes} flow.
\newblock {\em Int. J. Numer. Meth. Engng.}, 69(7):1374--1404, 2007.

\end{thebibliography}

\end{document}